\magnification=1200
\hoffset=.13in  
\voffset=.2in
\input amssym.tex 
\input amssym.def
\outer\def\beginsection#1\par{\filbreak\bigskip\leftline{\bf#1}\nobreak
\smallskip\vskip-\parskip\noindent}
 
\def\geq{\geqslant}
\def\leq{\leqslant}
\def\le{\leq}

\footline{\hss\oldstyle\number\pageno\hss}
\headline{\hss\it\sectitle} 
\def\sectitle{\firstmark}
 
\outer\def\proclaim #1. #2\par{\medbreak
\noindent{\rm#1.\enspace}{\sl#2}\par
\ifdim\lastskip<\medskipamount \removelastskip\penalty55\medskip\fi}
 
\def\spec#1#2{\hbox{\dimen0=#1truein\vrule
\vbox{\tabskip= \dimen0\baselineskip= 27pt
\halign{\kern -.9\dimen0\rlap{$\scriptstyle##$}\hfil&&\rlap
{$\scriptstyle##$}\hfil\cr#2\crcr}\smallskip\hrule}}}
 
\def\bigspec#1#2{\hbox{\dimen0=#1truein\vrule
\vbox{\tabskip= \dimen0\baselineskip= 25pt
\halign{\kern -.9\dimen0\rlap{$\textstyle##$}\hfil&&\rlap
{$\textstyle##$}\hfil\cr#2\crcr}\smallskip\hrule}}}
 
\def\sentry#1{\null\,\vcenter{\normalbaselines
\ialign{$\scriptstyle##$\hfil&&\quad$\scriptstyle##$\hfil\crcr
\mathstrut\crcr\noalign{\kern-\baselineskip}
#1\crcr\mathstrut\crcr\noalign{\kern-\baselineskip}}}\,}
 
\def\figure #1/#2/#3/{$$\hbox{\it#1\cddot#2.\ #3.}$$}
 
\def\eee #1,#2,#3,{E_{#1}^{#2,#3}}
\def\mapright#1{\smash{\mathop{\longrightarrow}\limits^{#1}}}
\def\mapdown#1{\Big\downarrow\rlap{$\vcenter{\hbox{$\scriptstyle#1$}}$}}
\def\proof{\medbreak \noindent{\it Proof. }}
\def\sone {\ifmmode {{\rm S}^1} \else$\sone$\fi}
\def\beginsubsection #1. {\medbreak \noindent{\bf #1.}\enspace }
\def\widebar{\overline}
\def\aut {{\rm Aut}}
\def\deg {{\rm deg}}
\def\res {{\rm Res}}
\def\cor {{\rm Cor}}
\def\hom {{\rm Hom}}
\def\inf {{\rm Inf}}
\def\pra{\par}
\def\ptwo{{\widetilde P}}
\def\tilb{{\widetilde B}}
\def\tilp{{\ptwo}}
\def\bee{{\rm B}}
\def\cohz #1#2{H^{#1}(\bee#2;\Bbb Z)}
\def\cohp #1#2{H^{#1}(\bee#2;\Bbb F_p)}
\def\cohr #1#2{H^{#1}(\bee#2;R)}
\def\cosr #1#2{H^{#1}(#2;R)}
\def\coht #1#2{H^{#1}(\bee#2;\Bbb F_3)}
\def\cod #1#2{H^{#1}(\bee#2)}
\def\ch#1{{\rm Ch}(#1)}
\def\chbar#1{{\rm \widebar {Ch}}(#1)}
\def\cher#1#2{{\rm Ch}^{#1}(#2)}
\def\cohev#1{H^{\rm ev}(\bee#1)}
\def\cohod#1{H^{\rm od}(\bee#1)}
\def\co#1#2#3{H^{#1}(\bee#2;#3)}
\def\coh#1{H^*(\bee#1)}
 
\def\held{{H\kern -0.1667em e}}
\def\steen{{\cal A}_p}
\def\chain{\dot{\hbox{\kern0.3em}}}
\def\pone{{\rm P^1}}
\def\qed{\ifmmode \eqno{\scriptstyle\blacksquare}\else
{{\unskip\nobreak\hfil\penalty50\hskip 2em\hbox{}\nobreak\hfil
${\scriptstyle\blacksquare}$\parfillskip=0pt\par}}\fi}
\def\ker{{\rm Ker}}
\def\im{{\rm Im}}

\def\phi {\varphi}
\def\epsilon {\varepsilon}
\def\nilc {nilpotency class\ }
\def\tilg {{\widetilde G}}
\def\cddot {$\cdot$}
\def\Lemma {L{\sevenrm EMMA\ }}
\def\Theorem {T{\sevenrm HEOREM\ }}
\def\Corollary {C{\sevenrm OROLLARY\ }}
\def\prop{P{\sevenrm ROPOSITION\ }}
\def\Definition {D{\sevenrm EFINITION\ }}
\def\book#1/#2/#3/#4/#5/#6/ {\item{[#1]} #2{\sevenrm #3}. {\it #4}.\
#5, (#6).\par\smallskip}
\def\paper#1/#2/#3/#4/#5/#6/#7/#8/ {\item{[#1]} #2{\sevenrm #3}. #4.
{\it #5} {\bf #6} (#7), #8.\par\smallskip}
\def\prepaper#1/#2/#3/#4/#5/ {\item{[#1]} #2{\sevenrm #3}. #4.
#5\par\smallskip}


\font\frmm=cmr10 scaled  1700
\font\trmm=cmr10 scaled  1400
\font\bd=cmr10 scaled 1200
\font\wwd=cmr10 scaled 1000
\hoffset=.13in
\
\vskip 1 in
\baselineskip 30 pt
\centerline{\frmm T{\trmm HE} C{\trmm OHOMOLOGY} {\trmm OF}}
\centerline{\frmm C{\trmm ERTAIN} G{\trmm ROUPS}}
\vskip 1.42 in
\centerline{\bd Ian James Leary}
\baselineskip 23pt
\vskip 1.42 in
\baselineskip 20pt
 
\centerline{\wwd January, 1992}
 
\vfill\eject

\baselineskip= 20 pt
 
\baselineskip 35 pt
\def\leaderfil{\leaders\hbox to 1em{\hss.\hss}\hfill}
\def\ooo{\oldstyle}
\beginsection Contents.  \par
\line{Introduction\leaderfil \ooo2}
\line{The Circle Construction\leaderfil \ooo6}
\line{The Cohomology Rings of Various $p$-Groups\leaderfil \ooo20}
\line{The Size of the Chern Subring and its Closure\leaderfil \ooo67}
\line{Further Calculations at the Prime Three\leaderfil \ooo79}
\line{The Integral Cohomology of the Held Group\leaderfil \ooo94}
\line{Yagita's Invariant\leaderfil \ooo109}
\line{The Davis Complex\leaderfil \ooo115}
\line{Appendix\leaderfil \ooo122}
\line{References\leaderfil \ooo125}

\baselineskip 20 pt
 
\beginsection Introduction. \par
This dissertation comprises six chapters which are a improved version
of the author's PhD thesis (submitted in November 1990), together with a short
chapter describing more other recent work not directly related to the rest of
the thesis.
Most of the work contained in this dissertation consists of calculations in the
integral and mod-$p$ cohomology rings of specific finite groups. A number of
questions concerning group cohomology are resolved by means of these examples,
as I shall explain below. I have tried to describe my results in sufficient
detail that they may provide a tool for further calculations, which accounts
for the length of the statements of some of the theorems. Of the seven chapters
of this dissertation, all except the first and sixth involve, to some degree,
calculations for specific groups.
\par
Chapter one consists mainly of an explanation of a technique used throughout
the calculational chapters. The technique was suggested by J.~Huebschmann and
P.~H.~Kropholler, but its subsequent development in this chapter, including its
close connection with questions concerning Chern classes and their
corestrictions, is my own work.
The technique involves embedding a finite group in a compact Lie group
whose identity component is a circle. The reason why I have chosen the
`topologist's notation' for the cohomology of a group (as the cohomology of
its classifying space) for this dissertation is to avoid confusion between the
cases when I consider the cohomology of the circle as a topological space
and as a Lie group. The remainder of chapter one consists of a generalisation
whereby a finite group is embedded in a Lie group with identity component
any torus, and a classification of Lie groups containing the circle as a
subgroup of index $p^n$ for $n\leq 3$. The difficult cases of this
classification are when $p=2$ (so that the circle subgroup need not be central)
and are not used in later chapters, but are included for completeness.
\par
Chapter six of this dissertation is not very closely related to the other chapters,
and might better be described as an appendix, except that I have chosen to
reserve that title for the brief discussion of $p$-groups of small order
(which is almost entirely unoriginal). Chapter six concerns an integer
invariant of finite groups, due to Yagita [Ya], which is related to free
actions of the groups on products of spheres. The original invariant is defined
in terms of images of the integral cohomology of the group, and the aim of the
chapter is to investigate a similar definition using images of the Chern
subring. Since the new invariant may be calculated directly from the complex
representation ring of the group, it is not surprising that the invariant
bears a relation to linear actions analogous to the relation of the original
invariant to free actions.
\par
Chapter two is the first of three chapters devoted to $p$-groups for odd
primes~$p$. The groups considered are the non-abelian groups expressible
as central extensions of a cyclic group by the elementary abelian group of
order~$p^2$. There are two isomorphism classes of such group of given order.
I determine, for each such group, its integral cohomology ring and the action
of its outer automorphism group upon this ring. These calculations are referred
to frequently throughout chapters three, four and five (often without comment).
For the non-metacyclic groups I also make similar calculations for the
mod-$p$ cohomology rings and determine the action thereon of the Steenrod
algebra. Similar calculations could be made for the metacyclic groups, which
may also be expressed as split extensions of a cyclic group by a group of
order~$p$, but I have not included these. Diethelm has shown that there are
(for odd~$p$) only two isomorphism types of ring that occur as mod~$p$
cohomology rings of split metacyclic $p$-groups [Die], so one really ought to
examine the actions of the automorphism group and the Steenrod algebra for
all split metacyclics, rather than just those with a cyclic subgroup of
index~$p$. Neither the method used here, nor Diethelm's method (which involves
consideration of the spectral sequence for the `defining extension') seems
suited to this task. I feel that the best approach to this problem will
involve the spectral sequence for these groups considered as extensions of
their derived subgroups by their abelianisations. Some of the results of
chapter two contain as special cases results of other authors; I verify
Lewis' calculations of the integral cohomology rings of non-abelian groups of
order~$p^3$, Thomas' result that Chern classes generate the even degree
cohomology of each of the groups considered in chapter two, and Thomas'
calculations of the integral cohomology rings of the metacyclic $p$-groups
with a cyclic subgroup of index~$p$.  Various preprints of Huynh-Mui dated
from {\oldstyle 1981} onwards concern the calculation of the mod-$p$
cohomology ring of the group of order~$p^3$ and exponent~$p$, but the methods
used are different from those used here, and I believe that Huynh-Mui does not
give a complete description of the relations between the generators that he
uses.
\par
Chapters three and four consist mainly of calculations in the integral
cohomology of the $p$-groups of $p$-rank two and nilpotency class three. It
is these groups (in the case $p\geq 5$) that are shown to provide
counterexamples to the conjecture of C.~B.~Thomas, that for all $p$-groups of
$p$-rank two Chern classes generate the even degree cohomology, and
$p$-group counterexamples to the conjecture of M.~F.~Atiyah, that Chern classes
generate the $E_\infty$-page of the Atiyah-Hirzebruch spectral sequence. These
results have also been obtained by N.~Yagita using Brown-Peterson cohomology.
I also show that for $p=3$ Chern classes do generate the even degree cohomology
of these groups, and that for $p\geq 5$ corestrictions of Chern classes suffice.
These results do not seem to be obtainable from the work of Yagita.
In chapter four I calculate the cohomology rings of some of these groups
for $p=3$, and thus exhibit non-isomorphic groups of order~$3^5$ with
isomorphic integral cohomology rings. The best result of a similar nature
for $p\geq 5$ in chapter four is the existence of non-isomorphic groups of
order~$p^6$ with isomorphic integral cohomology groups.
\par
Chapter five consists of an extended application of the results of chapter two
to the study of the Held group, one of the sporadic simple groups [He].
The $p$-torsion in the integral cohomology of the Held group is determined
for all odd primes~$p$, by calculating the image of the restriction map to
a Sylow-$p$ subgroup. This provides enough information to determine the
$E_2$-page of the Rector spectral sequence [Re].
I intend to examine this spectral sequence in the near future.
\par
Chapter seven of this dissertation contains a description of a method
due to Davis [Da] for
constructing contractible simplicial complexes with groups acting with finite
stabilisers.  Among other things this method can yield important cohomological
information concerning Coxeter groups.  This chapter does not contain any
original results however.
\par
I hope that the notation used in this dissertation is not unusual. My notation
for cohomology is largely similar to that of [Br] and [Th3], except that I
refer to $H^*({\rm B}G)$ instead of $H^*(G)$ for reasons explained above.
Throughout chapter five much of the notation I use is an approximation to
that of [Co]. I follow other authors in using Greek letters to stand for
integral cohomology elements and Roman letters to stand for mod-$p$ elements.
If $g$, $h$ are elements of a group~$G$, then $h^g=g^{-1}hg$ and
$[g,h]=g^{-1}g^h$. If $H$ is a subgroup of $G$, then $c_g$ is the map from
$H$ to $H^g$ given by $h\mapsto h^g$, and $c_g^*$ is the induced map in
cohomology. I speak of extensions of normal subgroups by quotient groups, and
say that these are central if the normal subgroup is contained in the centre
of the big group.
\par
\vfill
\eject

\def\hatg{{\widehat G}}
\def\hatd{{\widehat D}}
\def\hatq{{\widehat Q}}
\def\gl{{\rm GL}}
\def\aut{{\rm Aut}}
\def\out{{\rm Out}}
\def\inn{{\rm Inn}}
\def\tilq{{\widetilde Q}}
\eject
\beginsection{1.} The Circle Construction.\par
In this section we describe a method for examining the cohomology of a finite
group by embedding it in a compact Lie group of dimension one. This method was
suggested by P.~Kropholler and J.~Huebschmann, and has also been used by
B.~Moselle [Hu2], [Hu3], [Mo]. We derive general properties of the method,
which demonstrate why it is particularly relevant to the study of the Chern
subring. We also describe a generalisation of the original method, and classify
the compact Lie groups of dimension one and $p^n$ components, where~$n\leq 3$.
\par
\bigskip
\mark{\ }
 
\beginsubsection Definition and Properties of $\widetilde G$.
\par \noindent
Given a finite group $G$ and a central cyclic
subgroup $C$, we fix an embedding of $C$
into \sone, and define
$$\widetilde G = \sone \times _C G.$$
Then we have a commutative diagram:
$$\matrix{C & \mapright{} & G & \mapright{} & Q \cr \mapdown{} && \mapdown{} &
& \mapdown{} \cr
\sone & \mapright{} & \widetilde G & \mapright{} & Q} $$
If $M$ is a $G$-module on which $C$ acts trivially, we may consider $M$ as a
$\widetilde G$-module by letting \sone\ act trivially, and the
Lyndon-Hochschild-Serre spectral sequence for the second extension
is often simpler than that for the first. To find $H^*(\bee G;M)$, given
$H^*(\bee\widetilde G;M)$, we use the Serre spectral sequence of the fibration
$$\sone /C \cong {\widetilde G / G} \mapright{} {\rm B}G \mapright{}
{\rm B}\widetilde G. \eqno(1)$$
This spectral sequence has $E_2^{i,j}=0$ for $j>1$, so the only possible
non-zero differential is $d_2$.
The above was first suggested to the author by P.~Kropholler.
A similar idea occurs in J.~Huebschmann's papers [Hu2], [Hu3].
\par We shall examine the cases in which $M$ is $\Bbb Z$ or $\Bbb F_p$, with
the trivial $G$-action. If we let $R$ stand for either $\Bbb Z$ or $\Bbb F_p$,
then
$$H^*(\sone;R)=\Lambda [\xi],$$
the exterior algebra over $R$ on one generator $\xi$ of degree 1, and
$$\cohr * \sone =R[\tau],$$
the polynomial algebra over $R$ on one generator $\tau$ of degree 2. The
$E_2^{*,*}$ page of our spectral sequence is given by
$$E_2^{*,*}\cong \cohr * {\widetilde G} \otimes \Lambda [\xi],$$
except that the product structure is changed by a sign. The effect of this sign
change is to make $\xi$ anticommute with elements of $\cohr * {\widetilde G}$.
We shall call this ring the anticommutative ring generated by
$\cohr * {\widetilde G}$ and $\xi$.
It is clear that
$$d_2(\xi)=\widebar \tau,$$
where $\widebar \tau$ is the inflation to $\widetilde G$ of a generator $\tau$
for $\cohr 2 {\widetilde G/G}\cong\cohr 2 \sone$. This follows by naturality of
the spectral sequence with respect to maps of fibrations, in particular for the
following commutative diagram.
$$\matrix{{\sone/C}&\mapright{}&{\rm B}G&\mapright{}&{\rm B}\widetilde G \cr
\mapdown{}&&\mapdown{}&&\mapdown{}\cr
{\widetilde G/G}&\mapright{}&{\rm B}\{1\}&\mapright{}&{\rm B}{\widetilde G/G}}
$$
\mark{The Circle Construction}
Here, all the vertical maps are induced by the quotient map
$$\widetilde G\longrightarrow{\widetilde G/G},$$
and the lower fibration is, up to homotopy, the path-loop fibration over
${\rm B}\sone$. We may also think of $\widebar\tau$ as the first Chern class of
${\rm B}G$ considered as a circle or ${\rm U(1)}$ bundle over ${\rm B}\tilg$.
This follows because the first Chern class of the standard $\sone$ bundle over
${\rm B}\sone$ generates $\cohr 2 \sone$. With integer coefficients there is
another interpretation for $\widebar\tau$. Recall that there is for any compact
Lie group $K$ a
natural isomorphism between $\cohz 2 K$ and $\hom(K,\sone)$. If we
use this isomorphism to regard $\widebar \tau$ as a homomorphism from $\tilg$
to $\sone$, its kernel is the subgroup $G$.

\par
It appears that we may obtain another filtration of $\cohr * G$ by examining
the Eilenberg-Moore spectral sequence for the pullback square:
$$\matrix{{\rm B}G&\mapright{}&{\rm B}\widetilde G\cr
\mapdown{}&&\mapdown{}\cr
\{*\}&\mapright{}&{\rm B}{\widetilde G/G}}$$
The $E_2^{*,*}$ page of this spectral sequence is
$$E_2^{*,*}\cong {\rm Tor}^{*,*}_{R[\tau]}(R,\cohr * {\widetilde G} ).$$
To calculate $E_2^{*,*}$, we may construct a 2-stage free $R[\tau]$ resolution
of $R$ by
$$0\mapright{} R[\tau]\mapright{\times\tau} R[\tau]\mapright{}R$$
and we see that
$${\rm Tor}^{i,*}_{R[\tau]}(R,\cohr * {\widetilde G})\cong
\cases{\cohr * {\widetilde G}/{\widebar\tau}\cohr * {\widetilde G}
&for $i=0$ \cr
\ker \times{\widebar\tau}:\cohr * {\widetilde G}\rightarrow\cohr *
{\widetilde G} &for $i=1$}$$
 It is now easy to see that this spectral sequence collapses, and that its
$E_2^{*,*}$ page gives the same filtration of $\cohr * G$ as the
Lyndon-Hochschild-Serre spectral sequence for the fibration (1). Of course,
these spectral sequences are just alternative ways to view the Gysin sequence
for the \sone\ bundle ${\rm B}G$.
We now prove various results concerning $\tilg$. We assume in the statements
that $G$ is expressed as a central extension of $C$ by $Q$, and that $\tilg$ is
the corresponding extension of $\sone$ by $Q$.
\par

\proclaim \Lemma 1\cddot 1.  Given $G$ as above, consider the spectral sequence
of the fibration (1) with coefficients in the trivial $G$-module $R$.
\pra a) For $R=\Bbb Z$, the spectral sequence does not collapse.
\pra b) For $R=\Bbb F_p$, the spectral sequence collapses iff
$C$ is not contained in $G^pG'$.\par
\proof By $G^p$ we mean the subgroup of $G$ generated by $p$th
powers of elements of $G$.
If we let $\xi$ be a generator of $H^1(\sone;R)$, then as
above, in each case $E_2^{*,*}$ is the anticommutative ring generated by
$\cohr * {\widetilde G}$ and $\xi$,
and the spectral sequence collapses if and only if $d_2(\xi)=0$.
$G$ is finite, so $H^1(\bee G;\Bbb Z)=0$, so in case a)
$d_2:E_2^{0,1}\rightarrow E_2^{2,0}$ must be injective.
\par
$$\hom(G,\Bbb F_p)\cong H^1(\bee G;\Bbb F_p) \cong E_3^{0,1}\oplus E_3^{1,0},$$
$$\hom(G,\Bbb F_p)\cong H^1(\bee G;\Bbb F_p) \cong E_2^{1,0},$$
and the map $E_2^{1,0}\twoheadrightarrow E_3^{1,0}$ corresponds to
restriction from $\widetilde G$ to $G$. $E_3^{0,1}$ is either $C_p$ or 0,
so the spectral sequence collapses iff
there are homomorphisms from $G$ to
$\Bbb F_p$ which cannot be extended to $\widetilde G$. A homomorphism
$f:G\rightarrow \Bbb F_p$ can be extended to $\widetilde G$ iff
$f(C)=0$, and this will be true for all $f$ iff $CG'/G'$ is contained in
$(G/G')^p$ iff $C$ is contained in $G^pG'$.\qed\par
 
\proclaim \Lemma 1\cddot 2. For $n \geq 2$, the inclusion of $G$ in $\tilg$
restricts to an isomorphism from $\Gamma_n(G)$ to $\Gamma_n(\tilg)$, and from
$G^{(n-1)}$ to $\tilg^{(n-1)}$.
\par
\proof Any element of $\tilg$ can be expressed as a product $xg$, where $x$ is
central, and $g$ is in the image of $G$. In any group, if $x$ and $x'$ are
central, then for all $g,g'$ we have $[xg,x'g']=[g,g']$.\qed
 
\proclaim \Corollary 1\cddot 3. $G$ is nilpotent (resp. soluble) if and only if
$\tilg$ is, and if so they have the same nilpotency class (resp. soluble
length). \qed\par
 
\proclaim \Lemma 1\cddot 4. The following conditions on $G$ and $\tilg$ are
equivalent:
\pra a) $\sone\rightarrowtail\tilg\twoheadrightarrow Q$ is split;
\pra b) The extension class of $G$ in $H^2({\rm B}Q;C)$ is in the kernel of the
Bockstein from $H^2({\rm B}Q;C)$ to $H^3({\rm B}Q;\Bbb Z)$ associated
to the short exact sequence $\Bbb Z\rightarrowtail\Bbb Z\twoheadrightarrow C$;
\pra c) There is an extension $\widebar G$ of $\Bbb Z$ by $Q$ and a map
$\widebar G\twoheadrightarrow G$ such that the following diagram commutes:
$$\matrix{\relax\Bbb Z & \mapright{} & \widebar G & \mapright{} & Q \cr
\mapdown{} && \mapdown{} && \mapdown{} \cr
C & \mapright{} & G & \mapright{} & Q.} $$
\par
\proof Using the classification of extensions up to equivalence, see for
example [Th3], we see that a) holds if and only if the extension class of $G$
in $H^2({\rm B}Q;C)$ maps to zero under the map to $H^2({\rm B}Q;\sone)$
induced by the inclusion of $C$ in \sone. To verify that this is equivalent to
b), we consider the following map of short exact sequences:
$$\matrix{\relax\Bbb Z & \mapright{} & \relax\Bbb Z & \mapright{} & C\cr
\mapdown{} && \mapdown{} && \mapdown{} \cr
\relax\Bbb Z & \mapright{} & \relax\Bbb R & \mapright{} & \sone.}$$
The Bockstein is natural for maps of short exact sequences, so we obtain the
following commutative diagram:
$$\matrix{H^2({\rm B}Q;C) & \mapright{} & H^3({\rm B}Q;\Bbb Z) \cr
\mapdown{} && \mapdown{} \cr
H^2({\rm B}Q;\sone) & \mapright{} & H^3({\rm B}Q;\Bbb Z).}$$
The group $Q$ is finite, so the lower map is an isomorphism, and hence a) and
b) are equivalent. The equivalence of b) and c) follows from the cohomology
long exact sequence for the coefficient sequence $\Bbb Z\rightarrowtail
\Bbb Z\twoheadrightarrow C$.\qed
\proclaim \Lemma 1\cddot 5. Any complex representation of $G$ extends to one of
$\tilg$.\par
\proof Given $\rho: G \rightarrow {\rm Aut}(V)$, an irreducible
representation of $G$, $\rho$ must restrict to $C$ as scalar multiplication,
because an eigenspace for $C$ in $V$ would be a $G$ submodule.
Hence $\rho$ may be extended to \sone, then to $\widetilde G$ by defining:
$$\rho (g)= \rho (x) \rho (h)
{\rm \qquad where}\
g = xh {\rm \ for \ }x \in \sone {\rm \ and} \  h \in i(G).\eqno{{\scriptstyle
\blacksquare}}$$
\par
The close connection between the representation ring of a group and the
$K$~theory of its classifying space suggest the following :
\proclaim \Lemma 1\cddot 6. $K^*({\rm B}\tilg)$ restricts onto $K^*({\rm B}G)$.
Similarly, if $K(n)^*$ is the $n$th Morava $K$~theory, and $K(n)^{\rm od}({\rm
B}G)$ is trivial, then $K(n)^*({\rm B}\tilg)$ restricts onto $K(n)^*
({\rm B}G)$. \par
\proof The generalised Atiyah-Hirzebruch spectral sequence for the fibration
$${\rm B}G \rightarrow {\rm B}\tilg \rightarrow {\rm B}\sone$$
has $E_2^{*,*}$ page $\co i {\sone} {{\cal H}^j({\rm B}G)}$, converging to a
filtration of ${\cal H}^{i+j}({\rm B}\tilg)$, and the fibre map is the
restriction. If ${\cal H}$ is $K$~theory, or if ${\cal H}$ is Morava $K$~theory
and $G$ satisfies the condition above, then $E_2^{i,j}$ is trivial if $i$ or
$j$ is odd, so the spectral sequence collapses. The statement that this
condition holds for all finite $G$ is known as Ravenel's conjecture. \qed
 
\proclaim \Lemma 1\cddot 7. The subring ${\rm Ch}(\tilg)$ of the
integral cohomology ring of $\tilg$ generated by Chern classes
maps onto ${\rm Ch}(G)$ under the restriction. Similarly $\chbar \tilg$ maps
onto $\chbar G$, where $\chbar H$ is the subring of the integral cohomology
generated by corestrictions of Chern classes.
\par
\proof The first statement follows immediately from lemma 1\cddot 5. For the
second statement, we note that a subgroup of $\tilg$ of finite index must
contain $\sone$, so is of the form
\def\tilh{{\widetilde H}}
$\tilh$ for some subgroup $H$ of $G$, so the image of $\chbar \tilg$ is
contained in $\chbar G$. To show that the image is all of $\chbar G$, we note
that the following diagram commutes:
$$\matrix{\ch \tilh&\mapright{\cor}&\cohev \tilg \cr
\mapdown{\res}&&\mapdown{\res}\cr
\ch H &\mapright{\cor}&\cohev G \cr}$$
This is because $\tilh \cap G = H$, and $\tilh G =\tilg$, so we only need one
double coset in the restriction-corestriction formula (see [Br]).
\qed \par
\proclaim \Lemma 1\cddot 8. $\cohz {n+1} \tilg$ restricts onto $\cohz {n+1} G$
if and only if multiplication by $c(G)$, the first Chern class of ${\rm B}G$
as an \sone\ bundle over ${\rm B}\tilg$, is injective on $\cohz n \tilg$.\par
\proof Consider the spectral sequence for the fibration
$$ {\widetilde G / G} \mapright{} {\rm B}G \mapright{}{\rm B}\widetilde G. $$
The cokernel of the restriction map from $\cohz {n+1} \tilg$ to $\cohz {n+1} G$
is isomorphic to $E_\infty^{n,1}$, which is isomorphic to the kernel of
multiplication by $c(G)$ in $\cohz n \tilg$. \qed
\par
\proclaim \Theorem 1\cddot 9. $\cohev G = \ch G$ (resp. $\cohev G = \chbar G$) if
and only if $\cohev \tilg = \ch \tilg$ (resp. $\cohev \tilg = \chbar \tilg$)
and multiplication by $c(G)$ is injective on $\cohod \tilg$.\par
\proof By lemma 1\cddot 8 the condition on $c(G)$ is equivalent to $\cohev
\tilg$ maps onto $\cohev G$ by the restriction. The two conditions on $\tilg$
therefore imply that $\ch \tilg$ (resp. $\chbar \tilg$) restricts onto $\cohev
G)$, but by lemma 1\cddot 7 the image is $\ch G$ (resp. $\chbar G$). Conversely,
by lemma 1\cddot 7, the condition on $G$ implies that $\cohev \tilg$ restricts
onto $\cohev G$.\pra
If $\cohev \tilg$ restricts onto $\cohev G = \ch G$, we prove
by induction that $\cher {2n} \tilg = \cohz {2n} \tilg$. For any compact Lie
group $H$, $\cohz 2 H = \cher 2 H$, so without loss of generality we may assume
that $n$ is at least 2. Lemma 1\cddot 7 tells us that $\cher {2n} \tilg$ maps
onto $\cher {2n} G$, so $\cher {2n} \tilg$ maps onto $E_3^{2n,0}$ in the
spectral sequence for the \sone \ bundle ${\rm B}G$. Therefore, given $x \in
\cohz {2n} \tilg$, we may pick $y\in \cohz {2n-2} \tilg$ and $z\in \cher {2n}
\tilg$ such that $x = z +c(G) y$. Hence $x\in \cher {2n} \tilg$. The proof for
the case of $\chbar \tilg$ is similar. \qed
\par
\beginsubsection A Generalisation of the Circle Construction.
\pra\noindent
Let $G$ an extension of $A\cong (C_m)^n$ by $Q$ such that the `action map'
$\phi :
Q\rightarrow \gl_n(\Bbb Z_{/m})$ lifts to a map $\widebar\phi$ from $Q$ to
$\gl_n(\Bbb Z)$. We may use this map to define an action of $Q$ on the $n$
dimensional torus, $\Bbb T^n$, which restricts to the original action on a
subgroup of $\Bbb T^n$ isomorphic to $A$, which we will identify with $A$.
We may then form the semidirect product $\Bbb T^n\rtimes G$, that is the
manifold $\Bbb T^n\times G$ with multiplication defined by
$$(x,g)(y,h)=(xy^{\bar\phi\pi(g^{-1})}, gh).$$
The set $B=\{(a^{-1},a)|a\in A\}$ is a normal subgroup, and we define
$$\hatg=\Bbb T^n\rtimes_A G = {{\Bbb T^n\rtimes G}/B}.$$
$G$ is a subgroup of $\hatg$, so there is a fibration
$${\Bbb T^n/A}\cong{\hatg/G}\longrightarrow {\rm B}G\longrightarrow {\rm
B}\hatg,$$
but the fundamental group of ${\rm B}\hatg$ (which is of course isomorphic to
$Q$) will in general act non-trivially on the cohomology of the coset space
$\hatg/G$. In fact the action of $\pi_1({\rm B}G)$ on $H^1
(\hatg/G;\Bbb Z)$ is exactly the action of $Q$ on $\Bbb Z^n$ given by
$\bar\phi$. If we write $l(\_)$ for the soluble length, and $Q$ is soluble,
then clearly the following inequalities hold.
$$l(Q)\leq l(G)\leq l(\hatg)\leq l(Q)+1$$
There is however no analogue of lemma 1\cddot 2, and $\hatg$ is not in general
nilpotent. $\hatg$ depends on the choice of lifting $\bar\phi$, which will not
be unique. Even if we take $G$ to be a group of order 4 expressed as an
extension of $C_2$ by $C_2$, we may obtain $\hatg$ isomorphic to either
$\sone\times C_2$ or the orthogonal group ${\rm O}_2$ depending on our choice
of $\bar\phi$ from $C_2$ to $\{\pm 1\}$.
\par
There is an analogue of lemma 1\cddot 4, where we replace $C$ by $A$, and
$\Bbb Z$ by $\Bbb Z^n$ considered as a $\Bbb Z Q$ module via $\bar\phi$, with
exactly the same proof. Not all representations of $G$ will extend to $\hatg$,
although in each particular case it seems easy to describe the ones that will.
Lemma 1\cddot 6 breaks down because $G$ is not necessarily normal in $\hatg$.
In spite of all these shortcomings I believe that $\hatg$ will prove useful in
the study of various groups. For example, there is a family of 3-groups
of 3-rank two, which may be presented as split extensions of
$C_{3^t}\oplus C_{3^t}$ by $C_3$, where the action is induced from a
non-trivial action of $C_3$ on $\Bbb Z^2$. These groups are nilpotent of class
$2t$, so are unsuited to study via central extensions.
\par
\beginsubsection Classification of Circle groups on $p^n$ components.
\pra \noindent
Here we classify Lie groups $G$ with maximal connected subgroup is $\sone$
and group of components of order dividing $p^3$. $\aut(\sone)$ has order two,
so for $p$ odd $\sone$ is central in $G$. Central extensions of $\sone$ by $H$
are classified up to equivalence by $\co 2 H \sone$, which can be identified
with $\cohz 3 H$ via the Bockstein for the sequence $\Bbb Z \rightarrowtail
\Bbb R \twoheadrightarrow \sone$. In each of the cases we examine, either
$\cohz 3 H$ is trivial, or the action of $\aut(H)$ is transitive on the
non-zero elements, so there is at most one non-split extension. The
groups $H$ for which $\cohz 3 H$ is non-trivial are $(C_p)^2$, $(C_p)^3$,
$C_{p^2}\times C_p$, the dihedral group $D_8$ and $P_2$, the non-abelian group
of exponent $p$ (see appendix). For odd $p$, section 2 contains sufficient
information to check the third cohomology group of the groups $H$.
The case when $p = 2$ may be calculated similarly, or see [CE] for $Q_8$
and [Ev2] for $D_8$. The split extensions are just direct products $\sone
\times H$, and the non-split extensions are as described below. The comments
indicate why they are not isomorphic to the split extensions.
\par
\def\gr#1#2{\medskip\noindent$H\cong #1$.\ \ $H^3\cong #2$}
\gr{(C_p)^2}{\Bbb F_p}. The non-split extension is not abelian, and may be
obtained from either of the non-abelian groups $G$ of order $p^3$ as $\tilg$.
\par
\gr{C_{p^2}\times C_p}{\Bbb F_p}. As above the non-split extension is not
abelian, and may be obtained from any of the non-abelian extensions of $
C_p$ by $C_{p^2}\times C_p$ using the circle construction. \par
\gr{(C_p)^3}{\Bbb F_p^3}. The non-split extension is the direct product of
$C_p$ and the non-split extension of $\sone$ by $(C_p)^2$. This follows from
the fact that any element of $H^3(C_p^3)$ will restrict trivially to some $C_p$
subgroup. \par
\gr{P_2}{\Bbb F_p\oplus \Bbb F_p}. The non-split extension has nilpotency class
three, and may be obtained from any of the groups of order $p^4$ and nilpotency
class three via the circle construction.
\par
\gr{D_8}{\Bbb F_2}. The generator restricts non-trivially to either $C_2\times
C_2$ subgroup (this may be verified readily using the spectral sequence for
this extension), so the non-split extension has a non-abelian subgroup of index
two. \par
\medskip
Our approach to the classification of the non-central extensions is slightly
different. If $G$ is a non-central extension of $\sone$ by a group of order
$2^n$, let $K$ be the centraliser in $G$ of $\sone$. Then $K$ has index 2, and
the isomorphism type of $K$ as a group equipped with a map from $C_2$ to
$\out(K)$ is an invariant of the group. For fixed $K$ such a structure is
determined by an involution in $\out(K)$, and an isomorphism between
structures specified by involutions $t_1$ and $t_2$ is determined by an element
$s$ of $\out(K)$ such that $st_1=t_2s$. Thus the first stage of our
classification of non-central extensions of $\sone$ by a group of order $2^n$
is to take each $K$, a central extension of $\sone$ by a group of order
$2^{n-1}$, and find all conjugacy classes of involutions in $\out(K)$ which
restrict non-trivially to $\sone$. We shall classify only groups with 2, 4, and
8 components, so we need consider only one case in which $K$ is non-abelian,
that when $K$ is the non-split central extension of $\sone$ by $(C_2)^2$. First
we shall consider the case when $K$ is abelian. In this case $\out(K)=\aut(K)$,
and if we pick a representative for each of our conjugacy classes of
involution, these define distinct $\Bbb ZC_2$-module structures on $K$
containing the submodule $\widehat \sone$, that is $\sone$ with the non-trivial
$C_2$-action. For each such $M$ we calculate $H^2(\bee C_2;M)$ and the action
of $\aut_{\Bbb ZC_2}(M)$ (which is isomorphic to the centraliser in $\aut(K)$
of the involution used to define $M$) upon it. Finally we decide whether or not
elements in distinct $\aut(M)$ orbits really do give non-isomorphic groups. If
$K$ is not abelian, we pick a representative for each conjugacy class of
involution in $\out(K)$ restricting non-trivially to $\sone$, and use each of
these to define a $\Bbb ZC_2$-module structure $M'$ on the centre $Z(K)$ of
$K$. Either there are no extensions of $K$ by $C_2$ corresponding to this
involution, or such extensions are classified up to equivalence by $H^2(\bee
C_2;M')$ (see [Br]), depending on the vanishing of an obstruction in $H^3(\bee
C_2;M')$. (In the only case we need the obstruction group is trivial.) Once
again we are left to determine whether or not inequivalent extensions really
give non-isomorphic groups. There are sufficiently few cases to examine when
$K$ is non-abelian that {\it ad hoc} methods suffice. In the abelian case the
following lemma, suggested to the author by R.~E.~Borcherds, is useful.
\par
\proclaim \Lemma 1\cddot 10. If $H$ is a finite group with a fixed map from
$C_p$ to ${\rm Aut}(H)$, then a non-split extension of $H$ by $C_p$ and
the split extension cannot be isomorphic.\par
\proof If $G$ is an extension of $H$ by
$C_p$, it splits if and only if there is an element of order $p$ in
$G\setminus H$, so a non-split extension has the same number of elements of
order $p$ as $H$, whereas a split extension has more. \qed
\par \noindent
This lemma also applies to the case when $H$ is a compact abelian Lie group of
dimension one, because these contain only finitely many elements of order $p$.
\par
The non-central extensions of $\sone$ by $H$ are described below. We consider
$\sone$ as the complex numbers of unit modulus, with automorphism $z\mapsto
\bar z$.
\par
$H$ of order 2. In this case $M$ must be $\sone$ with the non-trivial $C_2$
action. $H^2(\bee C_2;M)\cong\Bbb F_2$, and we obtain two groups.
$$\eqalign{\widehat D &=\langle \sone, T\mid T^2=1, z^T=\bar z\rangle \cr
           \widehat Q &=\langle \sone, T\mid T^2=-1, z^T=\bar z\rangle \cr}$$
\def\hatd{{\widehat D}}
\def\hatq{{\widehat Q}}
$\hatd$ is isomorphic to ${\rm O}_2$, and $\hatq$ is isomorphic to the subgroup
of ${\rm SU}_2$ of elements of the forms
$$\hbox{$\pmatrix{z&0\cr 0&\bar z\cr}\qquad\hbox{and}\qquad
\pmatrix{0&z\cr -\bar z& 0 \cr}$}.$$
\par
$H$ of order 4. $\aut(\sone\times C_2)\cong C_2\times C_2$, so there are two
$C_2$ structures on $\sone \times C_2$ restricting non-trivially to $\sone$.
Considering $C_2$ as a subgroup of $\Bbb C^\times$, the two actions are
$$\vcenter{\halign{#\hfil&\qquad#\hfil&\enspace#\hfil\cr
$M_1$&$(z,a)$&$\mapsto (\bar z,a)$\cr
$M_2$&$(z,a)$&$\mapsto (\bar za,a)$\cr}}$$
$H^2(\bee C_2;M_1)\cong (\Bbb F_2)^2$, but two elements are exchanged by
$\aut(M_1)$, and we obtain three groups, with presentations
$$\langle \sone, A,T\mid A^2=T^2=1\quad z^A=z\quad A^T=A\quad z^T=\bar
z\rangle$$
$$\langle \sone, T\mid T^4=1\quad z^T=\bar z\rangle$$
$$\langle \sone, A,T\mid A^2=1\quad T^2=-1\quad z^A=z\quad A^T=A\quad z^T=\bar
z\rangle$$
\par
$H^2(\bee C_2;M_2)\cong \Bbb F_2$, and so we obtain two groups.
$$\langle \sone, A,T\mid A^2=T^2=1\quad z^A=z\quad A^T=-A\quad z^T=\bar
z\rangle$$
$$\langle \sone, T\mid T^4=-1\quad z^T=\bar z\rangle$$
Of these five groups, three have quotient isomorphic to $C_2\times C_2$, and
two $C_4$.
\par
We do not give presentations for the twenty-one non-central extensions of
$\sone $ by groups of order eight, but we list the possible `modules'
(we shall refer to a non-abelian $K$ equipped with a choice of involution in
$\out(K)$ as a $C_2$-`module'), and then
tabulate the groups that arise for each of them. In the abelian cases lemma
1\cddot10 together with the calculation of the isomorphism types of the finite
quotients suffices to show that extensions whose classes are in distinct
$\aut(M)$ orbits in $H^2(\bee C_2;M)$ do give rise to distinct groups.
\par \medskip
$\aut(\sone\times C_4)$ is a group of order 16 containing seven conjugacy
classes of involution, four of which restrict non-trivially to $\sone$.
Considering $C_4$ as contained in $\Bbb C^\times$, the following four modules
are representatives for the distinct isomorphism types.
$$\vcenter{\halign{#\hfil&\qquad#\hfil&\enspace#\hfil\cr
$M_1$&$(z,a)$&$\mapsto (\bar z,a)$\cr
$M_2$&$(z,a)$&$\mapsto (\bar za,a)$\cr
$M_3$&$(z,a)$&$\mapsto (\bar z,a^{-1})$\cr
$M_4$&$(z,a)$&$\mapsto (\bar za^2,a^{-1})$\cr
}}$$
\par
$\aut(\sone \times C_2^2)$ is isomorphic to $C_2\times S_4$, so has five
conjugacy classes of involutions, three of which restrict non-trivially to
$\sone$. The following are representatives of each class.
$$\vcenter{\halign{#\hfil&\qquad#\hfil&\enspace#\hfil\cr
$M_5$&$(z,a,b)$&$\mapsto (\bar z,a,b)$\cr
$M_6$&$(z,a,b)$&$\mapsto (\bar z,b,a)$\cr
$M_7$&$(z,a,b)$&$\mapsto (\bar zab,a,b)$\cr
}}$$
\par\medskip
There is also the non-split central extension of $\sone$ by $C_2\times C_2$,
which we shall refer to as $\tilq$, since it could be obtained from the
quaternion group by applying the circle construction. It may be presented as
follows.
$$\tilq =\langle\sone,A,B|A^2=B^2=1\quad z^A=z^B=z\quad [A,B]=-1\rangle$$
$\inn(\tilq)$ is isomorphic to $(C_2)^2$, $\out(\tilq)$ is isomorphic to
$C_2\times S_3$, and $\aut(\tilq)$ is the split extension of $(C_2)^2$ by
$C_2\times S_3$ where $C_2$ acts trivially and $S_3$ acts as $GL_2(\Bbb F_2)$.
It follows that any involution in $\out(\tilq)$ will give rise to a split
extension of $\tilq$ by $C_2$. There are two conjugacy classes of involution in
$\out(\tilq)$ which restrict non-trivially to $\aut(\sone)$, and we may choose
as lifts to $\aut(\tilq)$ of representatives of them the following elements.
$$\hbox{$\eqalign{t_1&:z\mapsto \bar z\cr t_2&:z\mapsto\bar z}$}\qquad
\hbox{$\eqalign{A&\mapsto A\cr A&\mapsto B}$}\qquad
\hbox{$\eqalign{B&\mapsto B\cr B&\mapsto A}$}$$
For each of these the $\Bbb ZC_2$-module structure of $Z(\tilq)$ is isomorphic
to $\widehat\sone$, the non-trivial $\Bbb ZC_2$-module structure on $\sone$.
Since $H^2(\bee C_2;\widehat\sone)$ has order 2 there are in each case two
equivalence classes of extensions of $\tilq$ by $C_2$. In each case these two
groups can be shown to be distinct by comparing the orders of elements in the
components of the group not contained in $\tilq$.
\par
The results of the above calculations are contained in figure~1\cddot1. Note
that in each case there is a split extension, and the isomorphism type of the
finite quotient of the split extension is listed first.
\par
\midinsert
$$\vcenter{\halign{\hfil#\hfil&\qquad\hfil#\hfil&\qquad\hfil#\hfil&\qquad
\hfil#\hfil\cr
Module&$|H^2(\bee C_2;M)|$&Extensions&Quotients\cr
$M_1$&4&3&$C_2\times C_4$, $C_2\times C_4$, $C_8$\cr
$M_2$&2&2&$C_2\times C_4$, $C_8$\cr
$M_3$&4&3&$D_8$, $D_8$, $Q_8$\cr
$M_4$&2&2&$D_8$, $Q_8$\cr
$M_5$&8&3&$(C_2)^3$, $(C_2)^3$, $C_2\times C_4$\cr
$M_6$&2&2&$D_8$, $D_8$\cr
$M_7$&2&2&$(C_2)^3$, $C_2\times C_4$\cr
$\widehat\sone=Z(\tilq)\hbox{ and }t_1$&2&2&$(C_2)^3$, $(C_2)^3$\cr
$\widehat\sone=Z(\tilq)\hbox{ and }t_2$&2&2&$D_8$, $D_8$\cr
}}$$
\figure 1/1/The non-central extensions of $\sone$ by groups of order eight/
\endinsert
We summarise our calculations in the following theorem.
\proclaim \Theorem 1\cddot11. For $p$ an odd prime there are 8 (resp. 3, 1)
groups consisting of $p^3$ (resp. $p^2$, $p$) circles. There are 29 (resp. 8,
3) groups consisting of 8 (resp. 4, 2) circles, but in only 8 (resp. 3, 1) of
these is the $\sone$ subgroup central. \qed \par
\beginsubsection Remarks. It is desirable to have a stronger result than lemma
1\cddot10, but informed opinion seems to be that for a fixed map from $Q$ to
$\aut(N)$ it may be possible for the split extension and a non-split extension
to be isomorphic groups. No examples seem to be known for $Q$ and $N$ both
finite, but J.~C.~Rickard suggested the following.
\proclaim \prop 1\cddot12. Let $N$ be $\prod_{n=1}^\infty C_4\times
\prod_{n=1}^\infty C_2$, and let $Q$ be $C_2$ acting trivially on $N$. Then
$$\co 2 Q N \cong \prod_{n=1}^\infty \co 2 {C_2}{C_4} \times
\prod_{n=1}^\infty \co 2 {C_2}{C_2},$$
and any extension $G$ whose class restricts to zero in
$\prod_{n=1}^\infty \co 2 {C_2}{C_4}$ is isomorphic to the split extension.
\par
\proof The split extension and all such $G$ are isomorphic to $N$. \qed\par

\beginsection 2. The Cohomology Rings of Various $p$-Groups.\par
In this section we shall calculate the integral and mod-$p$ cohomology rings of
an infinite family of $p$-groups for odd $p$, which we shall call $P(n)$. The
group $P(n)$ is defined for each $n\geq 3$, has order $p^n$, and may be
presented as below.
$$\langle A,B,C\mid A^p=B^p=C^{p^{n-2}}=[A,C]=[B,C]=1\quad [A,B]=C^{p^{n-3}}
\rangle$$
$P(3)$ is the non-abelian group of order $p^3$ and exponent $p$, which we also
refer to as $P_2$, and $P(4)$ is the second group on Burnside's list of groups
of order $p^4$ (see appendix). The centre of $P(n)$ is the subgroup generated
by $C$ of index $p^2$. Applying the circle construction described in section~1
to the whole centre we obtain the same group for all $n$, the unique
non-abelian group consisting of $p^2$ circles, which we shall call $\tilp$. The
metacyclic groups with a cyclic subgroup of index $p$ which we shall refer to
as $M(n)$, where $p^n$ is the order, also occur as subgroups of $\tilp$, but
their cohomology rings have been determined by other means [Th3], so we
consider them only briefly. \par \mark{\ }
\goodbreak \noindent
{\bf The Calculation of $\cohz * \tilp$}\par\noindent
We now begin our calculation of $\cohz * \tilp$ by examining the
spectral sequence with integer coefficients for $\tilp$ considered as an
extension of \sone\ by $C_p\oplus C_p$.
The $E_2$ page is readily seen to be generated by elements $\alpha, \beta\in
E_2^{2,0}$, $\gamma\in E_2^{3,0}$ and $\tau \in E_2^{0,2}$ subject only to the
relations $p\alpha=p\beta=0$, $p\gamma=0$ and $\gamma^2=0$. Note that $\tau $
has infinite order.
Since $E_2^{i,j}$ is trivial for $j$ odd, we see that all the even
differentials must vanish. The behaviour of the differentials is summarised in
the following lemma.
\proclaim \Lemma 2\cddot1. In the above spectral sequence there are exactly two
non-zero differentials, $d_3$ and $d_{2p-1}$. $d_3(\tau)$ is a non-zero
multiple of $\gamma$, and $E_4$ is generated by the classes of the elements
$\alpha, \beta, p\tau,\ldots, p\tau^{p-1}, \tau^p$ and $\tau^{p-1}\gamma$
(see figure 2\cddot 1). All of
these generators are universal cycles except for $\tau^{p-1}\gamma$, which is
mapped by $d_{2p-1}$ to a non-zero multiple of $\alpha^p\beta-\beta^p\alpha$.
The $E_\infty$ page is generated by the elements $\alpha, \beta, p\tau,\ldots,
p\tau^{p-1}, \tau^p$ subject only to the relations they satisfy as elements of
$E_2$, and  the relation $\alpha^p\beta =\beta^p\alpha$.
\par
\midinsert
$$\spec{0.65}{\cr \zeta&-&\zeta\alpha,\zeta\beta&-&
\sentry{\zeta\alpha^2,\zeta\alpha\beta,\cr
\zeta\beta^2}&\cr -&-&-&-&-\cr
\chi_{p-1}&-&-&\eta&-&\eta\alpha,\eta\beta&\cr &\cr &\cr
\chi_2&\cr -\cr
\chi_1\cr - \cr 1&-&\alpha,\beta&-&\alpha^2,\alpha\beta,\beta^2&-&
\sentry{\alpha^3,\alpha^2\beta,\cr
\alpha\beta^2,\beta^3\cr}&\cr} $$
\figure 2/1/The $E_4$ page of the spectral sequence of lemma 2\cddot 1/
\endinsert\par
\proof The derived subgroup of $\tilp$ consists of the subgroup of
its central $\sone $ of
\mark{The Cohomology Rings of Various $p$-Groups}
order $p$, so there can be no homomorphism from $\tilp$ to $\sone$
that restricts to an isomorphism from the centre to $\sone$. It follows by
considering the natural isomorphism $\cohz 2 G \cong \hom(G,\sone)$ that
the element $\tau$ cannot survive to $E_\infty$, so we must have $d_3(\tau)$ a
non-zero multiple of $\gamma$. This determines $d_3$ completely.
It may be checked that $E_4$ is isomorphic to
the subring of $E_2$ generated by
$\alpha, \beta, p\tau,\ldots, p\tau^{p-1}, \tau^p$ and $\tau^{p-1}\gamma$. All
these elements must be universal cycles, with the possible exception of
$\tau^{p-1}\gamma$, because the groups in which their images under $d_n$ lie
are already trivial. (The $E_4$ page of the spectral sequence is depicted in
figure~2\cddot1, where $\chi_i=p\tau^i$, $\zeta=\tau^p$, and
$\eta=\tau^{p-1}\gamma$.)
The only remaining potentially non-zero differential is
$d_{2p-1}(\tau^{p-1}\gamma)$. To complete this proof it suffices to show that
in the $E_\infty $ page the relation $\alpha^p\beta=\beta^p\alpha$ must
hold.\par
Let $Q$ be the quotient of $\tilp$ by its \sone\ subgroup, and take generators
$\alpha',\beta'$ for $\cohz 2 Q $ and $\gamma'$ for $\cohz 3 Q $. The
statement that $\gamma $ does not survive to $E_\infty $ in the spectral
sequence is equivalent to the statement that $\gamma'$ is mapped to zero by the
inflation map from $Q$ to $\tilp$. Now we calculate $\phi(\gamma')$, where
$\phi$ is the integral cohomology operation $\delta_p\pone\pi_*$, where $\pi_*$
is the map induced by the change of coefficients from $\Bbb Z$ to $\Bbb F_p$,
$\pone$ is a reduced power, and $\delta_p$ is the Bockstein for the sequence
$\Bbb Z \rightarrowtail \Bbb Z \twoheadrightarrow \Bbb F_p$. Taking $y,y' \in
\cohp 1 Q$ such that $\delta_p(y)=\alpha'$, and $\delta_p(y')=\beta'$, we see
that $$\phi(\gamma')=\delta_p\pone\pi_*(\gamma')=
\delta_p\pone(\beta_p(y)y'-\beta_p(y')y)=
\delta_p(\beta_p(y)^py'-\beta_p(y')^py)=\alpha'^p\beta'-\beta'^p\alpha'.$$
It follows that
$$\alpha^p\beta-\beta^p\alpha=\inf(\alpha'^p\beta'-\beta'^p\alpha')=
\inf(\phi(\gamma'))=\phi\inf(\gamma')=0. \qed$$\par
We are now ready to state our theorem on $\cod * \tilp$.
\proclaim \Theorem 2\cddot2. Let $p$ be an odd prime, and let $\tilp$ be the
group defined above. Then $H^*({\rm B}\tilp;\Bbb Z)$ is generated by elements
$\alpha$,\allowbreak$\beta$,\allowbreak$\chi_1,\ldots,\chi_{p-1}$,
\allowbreak$\zeta$, with
$$\deg(\alpha)=\deg(\beta)=2 \quad \deg(\chi_i)=2i \quad \deg(\zeta)=2p,$$
subject to the following relations:
$$p\alpha=p\beta=0$$
$$\alpha^p\beta=\beta^p\alpha$$
\hfil\hbox{$\alpha\chi_i=\cases{0 \cr -\alpha^p}\qquad
\beta\chi_i=\cases{0 &for $i<p-1$ \cr -\beta^p &for $i=p-1$}$}\hfil
$$\chi_i\chi_j=\cases{p\chi_{i+j} &$i+j<p$\cr p^2\zeta &$i+j=p$ \cr
p\zeta\chi_{i+j-p} &$p<i+j<2p-2$ \cr
p\zeta\chi_{p-2}+\alpha^{2p-2}+\beta^{2p-2}-\alpha^{p-1}\beta^{p-1}
&$i=j=p-1$}$$ \pra
Chern classes of representations of $\tilp $ generate the whole ring.
An automorphism of $\tilp$ sends $\chi_i$ to $\chi_i$
(resp.\ $(-1)^i\chi_i$) and $\zeta$ to $\zeta$ (resp.\ $-\zeta$) if it
fixes (resp. reverses) \sone. The effect of an automorphism on $\alpha$,
$\beta$ may be determined from their definition. Considered as elements of
$\hom(\tilp,\sone)$, $\alpha$ has kernel $\langle \sone, B\rangle$ and sends
$A$ to $e^{2\pi i/p}$, and $\beta$ has kernel $\langle \sone, A\rangle$ and
sends $B$ to $e^{2\pi i/p}$.
If we let $H$ be the subgroup generated by $B$ and elements of \sone\
we may define $$\chi_i=\cases{\cor_H^{\tilp}(\tau'^i) &for $i<p-1$ \cr
\cor_H^{\tilp}(\tau'^{p-1})-\alpha^{p-1} &for $i=p-1$}$$
where $\tau'$ is any element of $\cohz 2 H $ restricting to \sone\ as the
generator $\tau$.
Similarly, $\zeta =c_p(\rho)$, where $\rho$ is an irreducible representation
of $\tilp$ restricting to \sone\ as $p$ copies of the representation
$\xi$ with $c_1(\xi)=\tau $. \par
\proof First we note that in the $E_\infty$ page of the above spectral sequence
all the group extensions that we need to examine are extensions of finite
groups by the infinite cyclic group, so are split. The elements $\alpha$ and
$\beta$ defined in the statement above clearly yield generators for
$E_\infty^{2,0}$, and the relations between them are exactly the relations that
hold between the corresponding elements in the spectral sequence.
Let $\beta' $ in $\cod 2 H$ be the restriction to $H$ of
$\beta$, and take any choice of $\tau'$ as in the statement. We may show by
considering $\beta'$ and $\tau'$ as homomorphisms from $H$ to \sone\ that
conjugation by $A^i$ induces the map on $\cod 2 H$ that fixes $\beta'$
and sends $\tau'$ to $\tau'-i\beta'$. Now applying the formula for
$\res^G_K\cor^G_H$ (see for example [Br]) it follows that $\chi_i$
restricts to \sone\ as $p\tau^i$, so yields a generator for $E_\infty^{0,2i}$.
\par
Any irreducible representation of $\tilp$ has degree 1 or $p$, because $\tilp$
has an abelian subgroup of index $p$.
Let $\rho$ be the representation of $\tilp$ induced from a 1-dimensional
representation of $H$ with first Chern class $\tau'$. $\rho$ restricts to
\sone\ as $p$ copies of the representation with first Chern class $\tau$, so
its total Chern class restricts to \sone\ as $(1+\tau)^p$, and so $c_p(\rho)$
yields a generator for $E_\infty^{0,2p}$, and
$c_i(\rho)={1/p}{p\choose i}\chi_i+P_i(\alpha,\beta)$ for some polynomial
$P_i$. We shall show later that $P_i =0$.
\par
The restriction to $H$ of $\alpha$ is trivial, so by Fr\"obenius reciprocity
$$\alpha\cor^\tilp_H(\tau'^i)=\cor^\tilp_H(\res^\tilp_H(\alpha)\tau'^i)=0,$$
and the expressions given for $\alpha\chi_i$ follow. By calculating
$\alpha(\beta\chi_i)=\beta(\alpha\chi_i)$, we may deduce that $\beta\chi_i=0$
for $i < p-1$ and $\beta\chi_{p-1}=\lambda(\alpha^{p-1}\beta-\beta^p)
-\alpha^{p-1}\beta$ for some scalar $\lambda$. To show that $\lambda =1$ we
use the restriction map to $H$, and the formula for corestriction followed by
restriction.
$$\eqalign{\res_H^\tilp(\beta \chi_{p-1})&=\beta^\prime \sum_{i=0}^{p-1}
(\tau^\prime+i\beta^\prime)^{p-1}\cr
&= \beta^\prime \sum_{j=0}^{p-1} \tau^{\prime p-1-j} \beta^{\prime j}
\sum _{i=0}^{p-1} i^j}$$
Newton's formula tells us that
$$\sum_{i=1}^{p-1} i^j \equiv \cases{0 \enspace (p) &for $j \not\equiv 0
\enspace (p-1)$ \cr 1\enspace (p) &for $j \equiv 0 \enspace (p-1)$ }$$
so $\res^\tilp_H(\beta\chi_{p-1})=-\beta'^p$, and the required relation
follows. \par
We now know $\res^\tilp_\sone(\chi_i\chi_j)$, $\alpha\chi_i\chi_j$, and
$\beta\chi_i\chi_j$, which together imply the relations given for
$\chi_i\chi_j$. To complete the proof of the theorem we must determine the
effect of automorphisms of $\tilp$ on the $\chi_i$. We know that an
automorphism sends $c_i(\rho)$ to itself or $(-1)^i$ times itself depending
whether or not it reverses the sense of \sone, so it will suffice to show that
$\chi_i={1/p}{p\choose i}c_i(\rho)$. The character of $\rho$ is zero except on
$\sone$, so if $\theta $ is a 1-dimensional representation of $\tilp$
restricting trivially to \sone, then $\rho\otimes\theta$ is isomorphic to
$\rho$. If we apply the formula expressing $c.(\rho\otimes\theta)$ in terms of
$c.(\rho)$ and $c.(\theta)$ (see [At]) we obtain
$$c_i(\rho)=c_i(\rho \otimes \theta)=
\sum_{j=0}^i {p-i+j \choose j}c_1(\theta)^j c_{i-j}(\rho).$$
and hence inductively
$$c_i(\rho)c_1(\theta)=\cases {0 &for $i<p-1$ \cr -c_1(\theta)^p &for $i=p$}.$$
Since $\alpha$ and $\beta$ are possible values for $c_1(\theta)$ the required
result follows. We may show inductively that $\chi_i$ is in the subring
generated by Chern classes because $\chi_1$ is, and $\chi_1\chi_{i-1}$,
${1/p}{p\choose i}\chi_i$ are coprime multiples of $\chi_i$.
\qed
 
We are now ready to state our theorem on the integral cohomology of $P(n)$.
\proclaim \Theorem 2\cddot3. $\!$Let p be an odd prime and let $P(n)$ be as
defined above. Then $H^*({\rm B}P(n);\Bbb Z)$ is generated by elements
$\alpha,\allowbreak\beta,\allowbreak\mu,\allowbreak\nu,\allowbreak
\chi_1,\ldots,\chi_{p-1},\allowbreak\zeta$, with
$$\deg(\alpha)=\deg(\beta)=2\quad \deg(\mu)=\deg(\nu)=3\quad \deg(\chi_i)=2i
\quad\deg(\zeta)=2p$$
subject to the following relations:
$$p\alpha=p\beta=0\quad p\mu=p\nu=0\quad p^{n-3}\chi_1=0\quad p^{n-2}\chi_i=0
\quad p^{n-1}\zeta=0$$
$$\alpha\mu=\beta\nu$$
$$\alpha^p\beta=\beta^p\alpha\quad \alpha^p\mu=\beta^p\nu$$
$$\hbox{$\alpha\chi_i=\cases{0 & \cr -\alpha^p & }$\qquad
$\beta\chi_i=\cases{0 &for $i<p-1$ \cr -\beta^p &for $i=p-1$}$}$$
$$\hbox{$\mu\chi_i=\cases{0 & \cr -\beta^{p-1}\mu & }$\qquad
$\nu\chi_i=\cases{0 &for $i<p-1$ \cr -\alpha^{p-1}\nu &for $i=p-1$}$}$$
$$\chi_i\chi_j=\cases{p\chi_{i+j} &$i+j<p$ \cr
p^2\zeta &$i+j=p$ \cr p\zeta\chi_{i+j-p} &$p<i+j<2p-2$ \cr
p\zeta\chi_{p-2}+\alpha^{2p-2}+\beta^{2p-2}-\alpha^{p-1}\beta^{p-1}
&$i=j=p-1$}$$
$$\mu\nu=\cases{0 &for $n>3$ \cr \lambda\chi_3 &for $n=3,\ p>3, \lambda \in
\Bbb F_p^\times$ \cr
3\lambda\zeta &for $n=3,\ p=3,\  \lambda=\pm 1$}$$
\pra
Chern classes of representations of $P(n)$ generate $H^{\rm ev}({\rm B}P(n);
\Bbb Z)$.
Under an automorphism of $P(n)$ which restricts to the centre as
$C \mapsto C^j$, $\chi_i$ is mapped to $j^i\chi_i$, and $\zeta$ is mapped to
$j^p\zeta$. The effect of automorphisms on $\alpha$ and $\beta$ is
determined by the natural isomorphism $H^2({\rm B}G;\Bbb Z) \cong
\hom(G,\Bbb R/\Bbb Z)$, under which \pra
\hfil\hbox{$\eqalign{\alpha :A &\mapsto 1/p \cr
B &\mapsto 0 \cr C &\mapsto 0 }\quad \eqalign{\beta : A &\mapsto 0 \cr
B &\mapsto 1/p \cr C &\mapsto 0} \quad \eqalign{\chi_1 : A &\mapsto 0 \cr
B &\mapsto 0\cr C &\mapsto 1/{p^{n-3}}.}$}\hfil
\pra
An automorphism of $P(n)$ which sends $\alpha$ to $n_1\alpha+n_2\beta$,
$\beta$ to $n_3\alpha+n_4\beta$ and restricts to the centre as $C\mapsto C^j$
sends $\mu$ to $j(n_4\mu+n_3\nu)$ and $\nu$ to $j(n_2\mu+n_1\nu)$.
If $\gamma'$ in $H^2({\rm B}\langle B,C \rangle;\Bbb Z)$ is such that it
maps to the
following element of $\hom(\langle B,C \rangle, \Bbb R/\Bbb Z)$
$$\eqalign{\gamma ' : B &\mapsto 0 \cr C &\mapsto 1/{p^{n-2}}},$$
then $\chi_i$ is defined as follows:
$$\chi_i=\cases {\cor _{\langle B,C \rangle}^{P(n)}(\gamma'^i) &for $i<p-1$
\cr \cor_{\langle B,C \rangle}^{P(n)}(\gamma'^{p-1})-\alpha^{p-1}
&for $i=p-1$.}$$
\pra
These are, up to scalar multiples, equal to $c_i(\rho)$, where $\rho$ is a
$p$-dimensional irreducible representation of $P(n)$, whose restriction
to $\langle C \rangle$ is a sum of $p$ copies of the representation $\theta$,
with $c_1(\theta)=\res_{\langle C\rangle}^{\langle B,C \rangle}(\gamma')$.
In fact,
$ c_i(\rho)=\textstyle{{1 / p} {p \choose i}}\chi_i.$
Also, we may define $\zeta= c_p(\rho)$.
\par
\proof We examine the spectral sequence for ${\rm B}P(n)$ as an \sone-bundle
over ${\rm B}\tilp$. $E_2^{*,0}$ is isomorphic to $\cohz * {{P(n)}}$ and
$E_2^{*,*}$ is freely generated by $E_2^{*,0}$ and an element $\xi$ of infinite
order in $E_2^{0,1}$. We know that $\cod 2 {{P(n)}}\cong\hom(P(n),\sone)
\cong C_{p^{n-3}}\oplus C_p\oplus C_p$, so $d_2(\xi)$ must be
$\pm p^{n-3}\chi_1$. If we wanted to calculate the cohomology of the metacyclic
groups $M(n)$ described above, the differential
in this spectral sequence would send $\xi$ to
$\pm p^{n-3}\chi_1+\gamma$ for some non-zero $\gamma$ in $\langle \alpha,\beta
\rangle$. It is now easy to see that $E_\infty$ is generated by the elements
$\alpha,\allowbreak\beta,\allowbreak
\mu=\beta\xi,\allowbreak\nu=\alpha\xi,\allowbreak
\chi_1,\ldots,\chi_{p-1}$ and $\zeta$ subject to the relations they satisfy as
elements of $E_2^{*,*}$ together with $p^{n-3}\chi_1=0$, $p^{n-2}\chi_i=0$, and
$p^{n-1}\zeta=0$. For each~$m$, the filtration of $\cod m {P(n)}$ given by
the $E_\infty$ page is trivial, so we may use the same symbols to denote
elements of $\cod m {P(n)}$, and the relations that hold in $E_\infty$
determine all the relations that hold in $\cod m {P(n)}$ except for the
product of the two odd dimensional generators.
\par
We know that $p\mu\nu=0$, and
the relation $\alpha\mu=\beta\nu$ implies that $\alpha\mu\nu=\beta\mu\nu=0$,
and so $\mu\nu$ must be a multiple of $p^{n-3}\chi_3$ for $p\geq 5$ (resp.\
$3^{n-2}\zeta$ for $p=3$). Note that these elements restrict to zero
on all proper subgroups of $P(n)$.
In the case of $P(3)$, Lewis [Lew] shows that $\mu\nu$ is not zero by
considering the spectral sequence for $P(n)$ considered as an extension of a
maximal subgroup by $C_p$. A similar method will work in general, but
we offer an alternative proof that involves expressing $\mu$ and
$\nu$ as Bocksteins of elements of $\cohp 2 {P(n)}$. This proof is
contained in lemma~2\cddot4 and corollary~2\cddot5.
\par
The effect of automorphisms on $\chi_i$ and $\zeta$ is easily seen to be as
claimed from their alternative definitions as Chern classes. To determine the
effect of automorphisms on $\mu$ and $\nu$, we note that an automorphism of
$P(n)$ restricting to the centre as $C\mapsto C^j$ extends to an endomorphism
of $\tilp$ which wraps the central circle $j$ times around itself, so induces a
map of the above spectral sequence to itself sending $\xi$ to $j\xi$. This
completes the proof of theorem~2\cddot3 modulo lemma~2\cddot4 and its corollary.\qed\par
 
We now examine the spectral sequence with $\Bbb F_p$ coefficients for the
central extension $C_{p^{n-2}}\rightarrowtail P(n)\twoheadrightarrow C_p\oplus
C_p$. Take generators so that $\cohp * {C_p\oplus C_p}\cong
\Bbb F_p[x,x']\otimes \Lambda[y,y']$, where $\beta_p(y)=x$, $\beta_p(y')=x'$,
and $\cohp * {C_{p^{n-2}}}\cong \Bbb F_p[t]\otimes\Lambda[u]$, where
$\beta_p(u)=t$ for $n=3$ (resp. $\beta_p(u)=0$ for $n\geq 4$). Then the $E_2$
page is isomorphic to $\Bbb F_p[x,x',t]\otimes \Lambda[y,y',u]$, and the first
two differentials are as described in the following lemma.
\proclaim \Lemma 2\cddot4. With notation as above, identify the
elements $x,x',y,y'$ in the spectral sequence
with their images in $\cohp * {P(n)}$ under the inflation map.
\item{1)} Let $n\geq 4$. Then $d_2$ is trivial, and $d_3(t)$ is a non-zero
multiple of $xy'-x'y$. The elements $x,x',\allowbreak yy',u'y,\allowbreak
u'y'$, form a basis for $\cod 2
{P(n)}$, where $u'$ is any element of $\cod 1 {P(n)}$ restricting to
$C_{p^{n-2}}$ as $u$. \pra
\item{2)} Let $n=3$. Then $d_2(u)$ is a non-zero multiple of $yy'$, $d_2(t)=0$,
and $E_3$ is generated by $y,y',x,x',[uy],[uy']$ and $t$ subject to the
relation $yy'=0$ and those implied by the relations in $E_2$. In particular
$[uy]y'=-[uy']y$ but this element is non-zero. As in the case $n\geq 4$,
$d_3(t)$ is a non-zero multiple of $xy'-x'y$. Let $Y,Y'$ be elements of $\cod 2
{P(3)}$ such that $x,x',Y,Y'$ form a basis for $\cod 2 {P(3)}$, and let
$X=\beta_p(Y)$, $X'=\beta_p(Y')$. Then $yY',xy,\allowbreak xy',x'y',
\allowbreak X,X'$, form a basis
for $\cod 3 {P(3)}$ and $xX,xX',\allowbreak x'X',x^2y,\allowbreak x^2y',xx'y',
\allowbreak x'^2y',YX'$, form a basis
for $\cod 5 {P(3)}$. The $E_4$ page of this spectral sequence is depicted in
figure 2\cddot 2, where boldface numerals indicate numbers of generators
required in the case $p>3$, and plain numerals the case $p=3$.
\par
\midinsert
$$\spec{1.06}{-&-&{\bf 1}\quad{1}&{\bf 1}\quad{1}&
{\bf 0}\quad{2}&{\bf 2}\quad{2}\cr
-&[t^2y],[t^2y']&-&{\bf 3}\quad{4}&-&{\bf 4}\quad{6}\cr
-&-&[uy][ty']&{\bf 1}\quad{1}&-&{\bf 2}\quad{2}\cr
-&[ty],[ty']&-&\sentry{x[ty],x[ty'],\cr x'[ty']}&-&{\bf 4}\quad{4}\cr
-&[uy],[uy']&[uy]y'&\sentry{x[uy],x[uy'],\cr x'[uy],x'[uy']}&-&{\bf 6}
\quad{6} \cr 1&y,y'&x,x'&\sentry{xy,xy',\cr x'y'}&
\sentry{x^2,xx',\cr x'^2}&\sentry{x^2y,x^2y',\cr xx'y',x'^2y'}}$$
\figure 2/2/The $E_4$ page of the spectral sequence of lemma 2\cddot4/
\endinsert
\par
\proof 1). In this case $H^1$ has order $p^3$, so $u$ must survive. The element
$xy'-x'y$ is the image under $\pi_*$ of a generator for $\cohz 3
{(C_p\oplus C_p)}$, so must be killed by some differential. We have already
shown that it cannot be killed by $d_2$, so the only possibility is that $t$
survives until $E_3$ and kills it. The rest of the statement follows easily.
\par
2). In this case $H^1$ has order $p^2$, so $d_2(u)$ must be non-zero. It is
true in general that if $G$ is a central extension of $C_p$ by $Q$, then in the
corresponding spectral sequence with $\Bbb F_p$ coefficients
$d_2:E_2^{0,1}\rightarrow E_2^{2,0}$ must kill the extension class. This
follows by naturality, since one may regard the extension class as defining a
homotopy class of maps from ${\rm B}Q$ to $K(C_p,2)$ such that ${\rm B}G$ is the
 ${\rm B}C_p$-bundle induced by the path-loop fibration over $K(C_p,2)$. Since
all subgroups of $P(3)$ of order $p^2$ are copies of $C_p\oplus C_p$, the
extension class of $P(3)$ must restrict to zero on all cyclic subgroups, so
must be a multiple of $yy'$. The transgression commutes with the Bockstein so
$d_2(t)=0$ and $d_3(t)=\beta_pd_2(u)$.
\par
Given the values of these differentials it is routine to compute the $E_4$ page
of the spectral sequence. If we write $E_r^n=\bigoplus_{i+j=n}E_r^{i,j}$, then
$\{[uy],[uy'],x,x'\}$ forms a basis for $E_4^2=E_\infty^2$, and
$\{[ty],[ty'],[uy]y',xy,xy',x'y'\}$ forms a basis for $E_4^3=E_\infty^3$. The
spectral sequence operation $_F\beta$ introduced by Araki [Ar] and Vasquez
[Va] maps $[uy]$ to $[ty]$ and $[uy']$ to $[ty']$, so if $Y$ and $Y'$ are
chosen to yield the generators for $E_4^{1,1}$ their Bocksteins yield
generators for $E_4^{1,2}$. A basis for $E_4^5$ is given by the eight elements
of the statement, which we know to be universal cycles, and the elements
$[t^2y]$, $[t^2y']$. $E_4^4$ consists of universal cycles, and the universal
coefficient theorem tells us that $H^5$ has order $p^8$, so
$[t^2y]$ and $[t^2y']$ cannot be universal cycles. \qed \par
The author has determined all the differentials in the above spectral sequence,
using theorems 2\cddot14 and 2\cddot15. The non-trivial differentials are
$d_2,\allowbreak d_3,
\allowbreak d_4,\allowbreak d_{2p-2}$, and $d_{2p-1}$.
The $E_\infty$ page of this spectral sequence in the case $p=7$
is depicted in figure~2\cddot4 (on page~65). \par
\proclaim \Corollary 2\cddot5. In $\cohz * {P(n)} $ the product
$\mu\nu$ is non-zero if and only if $n=3$. \par
\proof In the notation of lemma~2\cddot4 it suffices to determine
$\delta_p(u'y)\delta_p(u'y')$ in the case $n\geq 4$, and
$\delta_p(Y)\delta_p(Y')$ in the case $n=3$. In the case when $n=3$,
$$\delta_p(Y)\delta_p(Y')=\delta_p(Y\beta_p(Y'))=\delta_p(YX').$$
The kernel of $\delta_p : \cohp 5 {P(3)} \rightarrow \cohz 6 {P(3)}$ is equal to
$\pi_*(\cohz 5 {P(3)})$, which is generated by $xX$, $xX'$ and $x'X'$, so by
lemma~2\cddot4 $\delta_p(YX')$ is non-zero.\par
In the case when $n=4$, $\cohz i {P(4)} $ has exponent $p$ for $i=2,3$, so
$\pi_*$ is injective from these groups,
and $\ker\beta_p:\cod 2 {P(4)} \rightarrow
\cod 3 {P(4)}$ is equal to $\beta_p(\cod 1 {P(4)})$. $\beta_p(yy')=xy'-x'y=0$,
so we may choose the element $u'$ in lemma~2\cddot4 so that $\beta_p(u')=\lambda yy'$
for some non-zero $\lambda$. Then we have
$$\delta_p(u'y)\delta_p(u'y')=\delta_p(u'y\beta_p(u'y'))=
\delta_p(u'y(\lambda yy'y'-u'x'))=0.$$
The case when $n\geq 5$ is similar but simpler, since $u'$ may be chosen so
that $\delta_p(u')=p^{n-4}\chi_1$, which implies that $\beta_p(u')=0$. \qed
\par
\noindent{\bf Remarks.} The author can not think of a method to determine the
non-zero constant $\lambda$ in the equation for $\mu\nu$ in the case
$\lambda=3$. If we replace the generators $\alpha$ and $\nu$ by
$\lambda^{-1}\alpha$ and $\lambda^{-1}\nu$ respectively and leave the other
generators unchanged, the new generating set will satisfy exactly the same
relations as the old ones, except that the constant $\lambda$ is replaced by 1.
\par
Theorem~2\cddot3 contains independent proofs of Thomas'
result that the even degree subring of
$\cohz * {P(n)}$ is generated by Chern classes [Th2],
and Lewis' calculation of $\cohz * {P(3)} $ [Lew].
Our notation differs slightly from
that of Lewis. We have renumbered the generators $\chi_i$ (note that
$\chi_1$ vanishes
for $n=3$). Also our $\chi_{p-1}$ and Lewis' $\chi_{p-2}$ are related
by the formula $$\chi_{p-2}^{\rm Lewis}=\chi_{p-1}+\alpha^{p-1}+\beta^{p-1}.$$
Our result disagrees with that of AlZubaidy [Al2]. We now use a similar method
to obtain the result of Lewis [Lew] and Thomas [Th1] concerning the metacyclic
groups $M(n)$, which may be presented as follows.
$$M(n)=\langle A,B \mid A^p=B^{p^{n-1}}=1\qquad [B,A]=B^{p^{n-2}}\rangle$$
\par
\proclaim \Theorem 2\cddot6. Let $M(n)$ be the metacyclic $p$-group defined
above. Then $\cohz * {M(n)}$ is generated by elements $\alpha,\allowbreak
\chi_1,\ldots,\chi_{p-1},\allowbreak \zeta,\allowbreak\eta$, with
$$\deg(\alpha)=2\quad \deg(\chi_i)=2i \quad\deg(\zeta)=2p\quad \deg(\eta)=2p+1$$
subject to the following relations:
$$p\alpha=0\quad p^{n-2}\chi_i=0 \quad p^{n-1}\zeta=0 \quad p\eta=0$$
$$\hbox{$\alpha\chi_i=\cases{0 & \cr -\alpha^p & }$\qquad
$\eta\chi_i=\cases{0 &for $i<p-1$ \cr -\alpha^{p-1}\eta &for $i=p-1$}$}$$
$$\chi_i\chi_j=\cases{p\chi_{i+j} &$i+j<p$ \cr
p^2\zeta &$i+j=p$ \cr p\zeta\chi_{i+j-p} &$p<i+j<2p-2$ \cr
p\zeta\chi_{p-2}+\alpha^{2p-2}
&$i=j=p-1$}$$
\pra
The even-dimensional generators are multiples of Chern classes of irreducible
representations of $M(n)$. Considered as elements of $\hom(G,{\Bbb R/\Bbb Z})$,
$\alpha $ sends $A$ to $1/p$ and $B$ to 0, and $\chi_1$ sends $A$ to 0 and $B$
to $1/{p^{n-2}}$. Under an automorphism of $M(n)$ which sends $\alpha$
to $n_1\alpha+n_2p^{n-3}\chi_1$ and restricts to the centre as $B^p\mapsto
B^{jp}$, $\chi_i$ is mapped to $j^i\chi_i$, $\zeta$ to $j^p\zeta$, and $\eta$
to $jn_1\eta$. Let $H$ be the subgroup of $M(n)$ generated by $B$. If
$\beta'$ is the element of $\cohz 2 H$
mapping $B$ to $1/{p^{n-1}}$, then we may define
$$\chi_i=\cases{\cor_H^{M(n)}(\beta'^i) &for $i < p-1$,\cr
\cor_H^{M(n)}(\beta'^{p-1})-\alpha^{p-1} &for $i=p-1$.}$$
\par
\proof As in the case of $P(n)$, we consider the spectral sequence for ${\rm B}
M(n)$ as an \sone-bundle over ${\rm B}\tilp$. The $E_2$ page is as in theorem
2\cddot3, and we may take $d_2(\xi)=\beta-p^{n-3}\chi_1$.  Since $A$ is in the
kernel of $\beta$ this is consistent with the relation $A^p=1$ in $M(n)$. The
$E_\infty$ page is seen to be generated by the elements
$\alpha,\allowbreak\chi_i,\allowbreak\zeta$, and
$\xi\alpha(\alpha^{p-1}-\beta^{p-1})$, subject to the relations they satisfy as
elements of $E_2$ together with $p^{n-2}\chi_i=0$, $p^{n-1}\zeta=0$. These
relations completely determine those that hold in $\cohz * {M(n)}$. The action
of ${\rm Aut}(M(n))$ is determined as in theorem 2\cddot3. \qed\par
\noindent{\bf Remarks.} Ths result confirms Thomas' calculation of $\cohz *
{M(n)}$ [Th1] which generalised Lewis' calculation for $\cohz * {M(3)}$ [Lew],
but note
that our generator in degree $2p-2$ is not the usual one. Their method is
simpler than the one used here (if one is not also interested in $P(n)$), but
does not seem to yield the action of
${\rm Aut}(M(n))$, which the author believes
to be a new result. The method involves considering $M(n)$ as an extension of
$C_{p^{n-1}}$ by $C_p$. The $C_{p^{n-1}}$ subgroups of $M(n)$ are not
characteristic however, so ${\rm Aut}(M(n))$ does not act on this
extension, or on the corresponding spectral sequence. \par

\beginsection Calculation of $\cohp * \ptwo$. \par
\beginsubsection The Spectral Sequence.
We proceed as for the corresponding calculations in integral cohomology.
We start by examining the spectral sequence with mod-$p$ coefficients for the
central extension:
$$\sone\rightarrowtail\ptwo\twoheadrightarrow C_p^{\widebar A}
\oplus C_p^{\widebar B}.$$
If we write
$$\cohp * {C_p^{\widebar A}\oplus C_p^{\widebar B}}=
\Bbb F_p[x,x']\otimes \Lambda[y,y'],$$
where $\beta(y)=x$ and $\beta(y')=x'$, and
$$\cohp * \sone =\Bbb F_p[t]$$
then we see that:
$$E_2^{i,j}=\Bbb F_p[x,x',t]\otimes\Lambda[y,y'].$$
\par We must have $d_3(t)=\lambda(xy'-x'y)$ for some non-zero $\lambda$,
and so $E_4^{*,*}$ is generated by elements $y,\allowbreak y',\allowbreak
x,\allowbreak x',\allowbreak [t^p],\allowbreak
[tyy'],\ldots,[t^{p-1}yy']$, and $ [t^{p-1}(xy'-x'y)]$, subject to the relations
that they satisfy in $E_2^{*,*}$, together with the relations
$$t^ixy'=t^ix'y \hbox{  for } i \not\equiv -1 \qquad (p).$$
If we write $z$ for $[t^p]$, $c_i$ for $[t^{i-1}yy']$, and $d$ for
$[t^{p-1}(xy'-x'y)]$, then $E_4^{*,*}$ is as shown in figure 2\cddot3.
\midinsert
$$\spec{0.9}{&\cr z&zy,zy'&\sentry{zx,zx',\cr zyy'}&
\sentry{zxy,zxy',\cr zx'y'}&\ldots\cr
-&\cr -&-&c_p&d&c_px,c_px'&dx,dx'&\ldots\cr
-&\cr -&-&c_{p-1}& \cr &\cr &\cr
-&-&c_2 &\cr -&\cr
1&y,y'&x,x'&\sentry{xy,xy',\cr x'y'}&\sentry{x^2,xx',\cr x'^2,yy'}\cr}$$
\figure 2/3/The $E_4$ page of the spectral sequence discussed above/
\endinsert
\par
Using the universal coefficient theorem, and our calculation of $\cohz * \ptwo$,
we find that
$$|\cohp {2i} \ptwo |=|\cohp {2i-1} \ptwo |=p^{i+2}
\hbox{  for } i=1,\ldots,p-1,$$
and we deduce that $c_i$ is a universal cycle for $i \le p-1$, and hence
 $E_{2p-1}^{*,*} \cong E_4^{*,*}$.
\par
The mod-$p$ reduction of $\zeta\in \cohz {2p} \ptwo $ restricts
non-trivially to \sone, so $z$ is a universal cycle. Applying $P^1$, and
$\beta P^1$ to $xy'-x'y$, we obtain:
$$\eqalign{0&=P^1(xy'-x'y) =x^py'-x'^py, \cr
0&=\beta(x^py'-x'^py) =x^px'-x'^px.}
$$
We know that $z,\allowbreak z(\lambda y+\mu y')$ are universal cycles, so
the only way the above relations can be introduced into $E_\infty^{*,*}$ is if
$$\eqalign{d_{2p-1}(c_p) &=\lambda(x^py'-x'^py) \cr
d_{2p-1}(d) &=\lambda'(x^px'-x'^px)}$$
for some non-zero scalars $\lambda,\lambda'$.
\par
Now it is easy to see that $E_{2p}^{*,*} \cong E_\infty^{*,*}$, and is
generated by elements $y,y',\allowbreak x,x',\allowbreak z,\allowbreak
c_2,\ldots,c_{p-1}$, with bidegrees as above, subject to the relations:
$$y^2=y'^2=0\qquad xy'=x'y$$
$$x^py'=x'^py\qquad x^px'=x'^px$$
$$c_iy=c_iy'=0\qquad c_ix=c_ix'=0 \qquad\hbox{for all $i$}$$
$$c_ic_j=0\qquad \hbox{for all $i,j$.}$$
\par
\beginsubsection Multiplicative Relations. We use the notation $\pi_*$ and
$\delta_p$ introduced during our integral calculations for the maps
$$\eqalign{\pi_*:\cohz n G &\longrightarrow \cohp n G  \cr
\delta_p:\cohp n G &\longrightarrow \cohz {n+1} G .}$$
It is easy to see that we may choose the elements
$\pi_*(\alpha),\allowbreak\pi_*(\beta),\allowbreak\pi_*(\chi_i)$, and
$\pi_*(\zeta)$ respectively to yield the generators $x,\allowbreak
x',\allowbreak c_i$, and $z$ of $E_\infty^{*,*}$. We define
generators of $\cohp * \ptwo $ by:
$$\eqalign{x &= \pi_*(\alpha) \cr x' &= \pi_*(\beta) \cr
c_i &= \pi_*(\chi_i) \qquad \hbox{for $i=2,\ldots,p-1$} \cr
z &= \pi_*(\zeta).}$$
Then elements $y,y'$ are uniquely defined by the equations:
$$\beta(y)=x\qquad \beta(y')=x'.$$
The relations between these generators follow from the spectral sequence and
the above definitions, and we obtain the following theorem, which has also been
obtained independently by Moselle [Mo].
\par
\proclaim \Theorem 2\cddot7.
$\cohp * \ptwo$ is generated by elements $y,y',\allowbreak x,x',\allowbreak
c_2,\ldots,c_{p-1},z$ with
$$\deg(y)=\deg(y')=1,\quad \deg(x)=\deg(x')=2, \quad\deg(c_i)=2i,\quad
\deg(z)=2p,$$
subject to the following relations:
$$xy' =x'y$$
$$\hbox{$\eqalign{ x^py' &=x'^py \cr \beta(y) &=x}$\qquad
$\eqalign{ x^px' &=x'^px \cr \beta(y')&=x'}$}$$
$$\hbox{{$c_iy=\cases{0 & \cr -x^{p-1}y &}$}\quad
{$c_iy'=\cases{0 &for $i<p-1$ \cr -x'^{p-1}y' &for $i=p-1$ }$}}$$
$$\hbox{{$c_ix=\cases{0 & \cr -x^p &}$}\quad
{$c_ix'=\cases{0 &for $i<p-1$ \cr -x'^p &for $i=p-1$ }$}}$$
$$\hbox{$c_ic_j=\cases{0 &for $i+j<2p-2$ \cr x^{2p-2}+x'^{2p-2}-x^{p-1}x'^{p-1}
&for $i=j=p-1$}$.}$$
\pra
An automorphism of $\ptwo$ sends $c_i$ to $c_i$ (resp. $(-1)^ic_i$) and $z$
to $z$ (resp. $-z$) if it fixes (resp. reverses) \sone. The effect of an
automorphism on $y,y',\allowbreak x,x'$ may be determined by their definitions.
Using the natural isomorphism $\cohp 1 \ptwo \cong \hom(\ptwo,\Bbb F_p)$, we
define:
$$\hbox{$\eqalign{y:A&\mapsto 1\quad (p) \cr B&\mapsto 0 \quad (p)}$\quad
$\eqalign{y':A&\mapsto 0\quad (p) \cr B&\mapsto 1\quad (p),}$}$$
and then define $x$ and $x'$ by the equations:
$$x=\beta (y)  \qquad x'=\beta (y').$$
We may define $c_i$ (resp. $z$) to be the image under $\pi_*$ of $\chi_i$
(resp. $\zeta$), as defined during the statement of theorem 3. $\pi_*(\chi_1)$
is $\lambda yy'$ for some non-zero $\lambda$.
\par
\proof We have done most of this already. We should check that the two
definitions of $x$, as $\pi_*(\alpha)$ and $\beta(y)$, agree. It suffices for
this to check that $\delta_p(y)=\alpha$, because $\beta=\pi_*\delta_p$.
This is true by naturality of the Bockstein for maps of short exact sequences
of coefficient modules, in particular for the map:
$$\matrix{\relax\Bbb Z &\mapright{\times p}
&\relax\Bbb Z &\mapright{} &\relax\Bbb F_p \cr
\mapdown{} &&\mapdown{\times {1/p}}&&\mapdown{}\cr
\relax\Bbb Z&\mapright{}&\relax\Bbb Q&\mapright{}&{\Bbb Q/\Bbb Z}.}$$
Also, we note that $\delta_p:\cohp {2i+1} \ptwo \rightarrow \cohz {2i+2}
\ptwo$ is
injective for all $i$, and $\delta_p(c_iy)=\chi_ix$, so the relations
involving $c_i$ and $y,y'$ follow.
\qed\par
Much of the action of the Steenrod algebra $\steen$ on $\cohp * \tilp$ is
determined implicitly by theorem~2\cddot7. For example $c_i$ is expressible as
$\cor^\tilp_H(t^i)$ ($-x^{p-1}$ when $i=p-1$), so $\pone(c_i)$ is easily
calculated. An example of this calculation occurs in the proof of
theorem~2\cddot14.
The only further piece of information required to specify the action of
$\steen$ is $\pone(z)$, which does not follow immediately from theorem
~2\cddot7, but may be detemined as follows.
\proclaim \prop 2\cddot8. With notation as in theorem~2\cddot7,
$$\pone(z)=zc_{p-1}.$$
\par
\proof Throughout this proof $K$ shall stand for any of the $p+1$ subgroups of
$\tilp$ of index $p$ (each of which is isomorphic to $\sone\oplus C_p$), and
$\res=\res^\tilp_K$. Let $\cohp * K = \Bbb F_p[t,\bar x]\otimes\Lambda[\bar
y]$, where $t$ restricts to \sone\ as the mod-$p$ reduction of the standard
generator of $\cohp 2 \sone$ (called $\tau$ in the proof of theorem~2\cddot2),
and $\beta(\bar y)=\bar x$. There are $p-1$ possible choices for $\bar y$ (and
hence for $\bar x$), and $p$ choices for $t$. If $\rho$ is the representation
of $\tilp$  mentioned in the statement of theorem~2\cddot2, then
$$\res(\pi_*c.(\rho))=\prod^{p-1}_{i=0}(1+t+i\bar x)=1-\bar x^{p-1}+t^p-
\bar x^{p-1}t,$$
and this expression is independent of the choice of $t$ and $\bar x$. It
follows that $$\hbox{$\res(c_i)=\cases{0&for $i<p-1$\cr -\bar x^{p-1}
&for $i=p-1$}$\qquad$\res(z)=t^p-\bar x^{p-1}t$.}$$
It follows that $\res(\pone(z))=\bar x^{2p-2}t-\bar x^{p-1}t^p$. A typical
element of $\cohp {4p-2} \tilp$ is of the form $\lambda zc_{p-1} +zP(x,x')
+Q(x,x')$. We shall show that $zc_{p-1}$ is the only such element restricting
correctly to each $K$. \par
First, note that $\res(Q(x,x'))=\lambda'\bar x^{2p-1}$, and that $\lambda'=0$
for $K=\langle \sone,AB^i\rangle$ (resp.\ $K=\langle\sone,B\rangle$) iff
$x'-ix$ (resp.\ $x$) divides $Q$. Hence $x^px'-x'^px$ must divide $Q$, so $Q$
is zero in $\cod * \tilp$. To complete the proof it suffices to show that for
each $\lambda''\in\Bbb F_p$, no polynomial $P(x,x')$ homogeneous of degree
$p-1$ can restrict to all $K$ as $\lambda''\bar x^{p-1}$. For every choice of
$K$ except $\langle \sone,A\rangle$ we may choose $\bar x$ to be $\res(x')$. We
have $p$ such $K$, and $x$ will restrict to each of them as a distinct multiple
of $\bar x$. Hence $P$ as above would have to satisfy $P(X,1)-\lambda''=0$ for
all $X\in\Bbb F_p$, but $P(X,1)$ cannot have $p$ roots. \qed\par

\beginsection Calculation of $\cohp * {P(n)}$ for $n \geq 4$.
\par
\beginsubsection The Spectral Sequence. We follow the method used in the
corresponding calculations for integral cohomology, and consider the spectral
sequence of the fibration:
$${\sone/Z(P(n))}\longrightarrow {\rm B}P(n)\longrightarrow {\rm B}
\widetilde P(n).$$
This has
$$E_2^{i,j}\cong H^i({\rm B}\ptwo; H^j(\sone;\Bbb F_p)),$$
so if we set
$$ H^*(\sone;\Bbb F_p)=\Lambda[u],$$
then $E_2^{*,*}$ is the anticommutative ring generated by $\cohp * \ptwo$ and
$u$.
We may apply lemma 1\cddot 1, and we deduce
that the spectral sequence collapses iff
$n>3$. For $n>3$ the relations in $\cohp * {P(n)}$ follow easily from those
in $\cohp * \ptwo$. In the case $n=3$, $d_2(u)=\lambda(yy')$ for some
non-zero $\lambda$, and so $E_\infty^{*,*}$ requires many generators which
appear in $E_\infty^{*,1}$. There is a corresponding increase in the number
of elements required to generate $\cohp * {P(3)}$ and consequently in the
number of relations we must calculate. We are forced to use new techniques,
which shall be introduced in the next chapter.
\par
\proclaim \Theorem 2\cddot9. Let $n$ be at least 4.
Then $\cohp * {P(n)}$ is generated
by elements  $u,y,y',\allowbreak x,x',\allowbreak c_2,\ldots,c_{p-1},z$ with
$$\deg(u)=\deg(y)=\deg(y')=1,\quad \deg(x)=\deg(x')=2, \quad\deg(c_i)=2i,\quad
\deg(z)=2p,$$
subject to the following relations:
$$ xy'=x'y$$
$$\hbox{$\eqalign{ x^py' &=x'^py \cr \beta(y) &=x}$\qquad
$\eqalign{ x^px' &=x'^px \cr \beta(y')&=x'}$}$$
$$\beta(u)={\textstyle\cases{y'y &for $n=4$\cr 0 &for $n>4$}}$$
$$\hbox{{$c_iy=\cases{0 & \cr -x^{p-1}y &}$}\quad
{$c_iy'=\cases{0 &for $i<p-1$ \cr -x'^{p-1}y' &for $i=p-1$ }$}}$$
$$\hbox{{$c_ix=\cases{0 & \cr -x^p &}$}\quad
{$c_ix'=\cases{0 &for $i<p-1$ \cr -x'^p &for $i=p-1$ }$}}$$
$$\hbox{$c_ic_j=\cases{0 &for $i+j<2p-2$ \cr x^{2p-2}+x'^{2p-2}-x^{p-1}x'^{p-1}
&for $i=j=p-1$}$.}$$
\pra
Under an automorphism of $P(n)$ which restricts to the centre as $C \mapsto
C^j$, $c_i$ is mapped to $j^ic_i$ and $z$ is mapped to $jz$. We define
$u,y,y'$ by regarding them as elements of $\hom(P(n),\Bbb F_p)$ then
$$\hbox{$\eqalign{u:A&\mapsto 0 \cr B&\mapsto 0 \cr C&\mapsto 1}\qquad
\eqalign{y:A&\mapsto 1 \cr B&\mapsto 0 \cr C&\mapsto 0}\qquad
\eqalign{y':A&\mapsto 0 \cr B&\mapsto 1 \cr C&\mapsto 0}$.}$$
This determines the effect of automorphisms on $u,y,y'$ and on $x=\beta(y)$,
$x'=\beta(y')$. We may define
$$c_i=\pi_*(\chi_i)\quad\hbox{and}\quad z=\pi_*(\zeta).$$
\par
\proof Almost all of this follows from our work above.
We may, of course choose $u$ to be any element of $\cohp 1 {P(n)}$ not in the
span of $y,y'$. With our definition it follows that
$\delta_p(u)=p^{n-4}\chi_1$, so we obtain
$$\beta(u)={\textstyle\cases{\lambda yy' &for some $\lambda \neq 0$ for $n=4$
\cr 0 &for $n>4$.}}$$
It can be checked by explicit calculation with 1- and 2-cochains that in the
case $n=4$
$$\beta(u) =y'y=-yy',$$
and hence also that for $n=4$
$$\pi_*(\chi_1)=y'y.$$
We leave the details to the interested reader.
\qed\par
Note that the action of $\steen$ on $\cohp * {P(n)}$ for $n\geq 4$ is
determined completely by information contained in theorem~2\cddot9 and
proposition~2\cddot8.
\par

\beginsection The Massey Product.
\par
This section is not original, but it introduces techniques we shall use to
construct various explicit elements in $\cohp * {P_2}$, and to find relations
between them. Throughout this section $R$ will be a
commutative ring on which $\pi_1(X)$ acts trivially, and $C\chain (X;R)$
the singular cochain complex of $X$ with coefficients in $R$.
Any cochains or
cohomology classes represented by a single letter will be homogeneous, and we
shall write, for example $(-1)^u$ for $(-1)^{\deg(u)}$. In particular,
$(-1)^{uv}$ shall mean $(-1)^{\deg(u)\deg (v)}$, not $(-1)^{\deg (uv)}$.
\beginsubsection Definition.
Let $u,v,w$ be elements of $C\chain (X;R)$, with
$$[uv]=0\quad [vw]=0\qquad \hbox{in $\cosr * X$.}$$
Then choose $a\in C^{u+v-1}(X;R)$ and $b\in C^{v+w-1}(X;R)$ with
$$\delta a=uv\qquad \delta b=vw,$$
and define the Massey product $\langle [u],[v],[w] \rangle$ by
$$\langle [u],[v],[w] \rangle=[(-1)^uub-aw]\in {\cosr {u+v+w-1} X /
(u\cosr {v+w-1} X + w\cosr {u+v-1} X)}.$$
The Massey product is only well-defined modulo $u\cosr {v+w-1} X +
w\cosr {u+v-1} X$ because of the freedom of choice of $a$ and $b$. It is
easily seen to be linear in each of its arguments.
Before stating some more properties enjoyed by the Massey product we state a
relation due to Hirsch [Hi] between the cup-0 and cup-1 products. We recall that
a cup-0 product is a product on cochains inducing the cup product on
cohomology, for example the standard product defined by the Alexander-Whitney
formula. A cup-1 product is a natural cochain transformation of degree $-1$
$$\smile_1:(C\chain(\_)^{\otimes 2})\longrightarrow C\chain(\_),$$
satisfying the following `coboundary formula' for all (homogeneous) cochains
$a$ and $b$:
$$\delta (a\smile_1 b)=-\delta a\smile_1 b -(-1)^aa\smile_1 \delta b
+ab+(-1)^{ab}ba.$$
\proclaim \Theorem 2\cddot 10. (Hirsch).  With the standard
choice of cup-0 product, there
is a choice of cup-1 such that for all cochains $a,b,c$
$$(ab)\smile_1 c=(-1)^aa(b\smile_1 c) +(-1)^{bc}(a\smile_1 c)b.$$
\par
Note that Hirsch uses Steenrod's original definition [St] of the cup-1 product,
which satisfies a slightly different coboundary formula to the one stated
above. This explains the difference between our statement of theorem~2\cddot10
and Hirsch's original statement. If we write $\smile_{\rm S}$ for Steenrod's
definition of the cup-1 product, then we may define a choice of cup-1 in the
modern sense by
$$a\smile_1b=(-1)^{a+b-1}(a\smile_{\rm S}b),$$
and this is a choice for cup-1 which satisfies the identity of theorem
2\cddot10. \par
\proclaim \Lemma 2\cddot11. The Massey product
satisfies the following identities, which
are valid whenever all the terms are defined, for any $u,v,w,x,y\in\cosr * X$:
$$\langle u,v,w \rangle x+(-1)^uu\langle v,w,x\rangle\equiv 0
\quad \hbox{mod $uH^*x$}\eqno{(1)}$$
$$\eqalignno{
(-1)^u\langle\langle u,v,w \rangle,x,y \rangle +\langle u,\langle v,w,x
\rangle,y\rangle+&(-1)^v\langle u,v,\langle w,x,y\rangle\rangle\equiv 0 \cr
&\hbox{mod $uH^*+H^{u+v-1}wH^{x+y-1}+yH^*$}&(2)}$$
$$(-1)^{wu}\langle u,v,w \rangle +(-1)^{uv}\langle v,w,u\rangle
+(-1)^{vw}\langle w,u,v \rangle \equiv 0
\quad \hbox{mod $uH^*+vH^*+wH^*$}\eqno{(3)}$$
$$\langle u,v,w\rangle +(-1)^{uv+vw+wu}\langle w,v,u \rangle \equiv 0
\quad\hbox{mod $uH^*+wH^*$}\eqno{(4)}$$
\par
\proof The verification of these is fairly straightforward. The proofs of (3)
and (4) use theorem 2\cddot10. For example, to prove (3), if we let
$u,v,w$ also stand for representative cocycles, then pick $a,b$ and $c$
such that $$\delta a=uv\quad \delta b=vw \quad\delta c=wu,$$
then
$$(-1)^{wu}\langle u,v,w \rangle +(-1)^{uv}\langle v,w,u\rangle
+(-1)^{vw}\langle w,u,v \rangle $$
contains the element
$$\eqalignno{&(-1)^{wu+u}ub-(-1)^{wu}aw+(-1)^{uv+v}vc-(-1)^{uv}bu+
(-1)^{vw+w}wa-(-1)^{vw}cv \cr =&-(-1)^{wu}(uv)\smile_1 w -(-1)^{uv}(vw)\smile_1u
-(-1)^{vw}(wu)\smile_1v \cr =&\delta ((u\smile_1v)\smile_1w).&
\scriptstyle\blacksquare}$$
\par
\beginsubsection Matrix Massey Products. This generalisation of the Massey
product is due to May [Ma].
We consider now homogeneous `vectors' of cocycles $(u_i),(v_{ij}),(w_j)$ of
degrees $u,v,w$ respectively, where $1\leq i\leq l$, and $1\leq j\leq m$,
with the property that
$$\eqalign{\sum_i[u_iv_{ij}]&=0 \hbox{  for all $j$} \cr
\sum_j[v_{ij}w_j]&=0\hbox{  for all $j$.}}$$
We choose $(a_j)$ and $(b_i)$ such that
$$\delta a_j=\sum_iu_iv_{ij}\qquad\delta b_i =\sum_jv_{ij}w_j,$$
then define the matrix Massey product $\langle [u_i],[v_{ij}],[w_j] \rangle$ by
$$\eqalign{
\langle [u_i],[v_{ij}],[w_j] \rangle&=(-1)^u\sum_iu_ib_i-\sum_ja_jw_j\cr
&\in {\cosr {u+v+w-1} X /(\sum_iu_i\cosr {v+w-1} X +\sum_jw_j\cosr {u+v-1} X
)}.}$$
The obvious generalisations of the properties claimed in lemma~2\cddot11
are valid for the matrix Massey product, with very similar proofs.
For example suppose that we are
given $(u_i),\allowbreak (v_{ij}),\allowbreak (w_{jk}),\allowbreak (x_k)$ for
$1\leq i\leq l$, $1\leq j\leq m$, $1\leq k\leq n$ with
$$\eqalign{\sum_i[u_iv_{ij}]&=0\hbox{  for all $j$}\cr
\sum_j[v_{ij}w_{jk}]&=0\hbox{  for all $i,k$}\cr
\sum_k[w_{jk}x_k]&=0\hbox{  for all $j$}.\cr}$$
We wish to show that
$$\sum_k\langle [u_i],[v_{ij}],[w_{jk}]\rangle[x_k]+
(-1)^u\sum_i[u_i]\langle [v_{ij}],[w_{jk}],[x_k]\rangle\equiv 0
\hbox{    mod }\sum_{i,k}u_iH^*x_k.$$
We choose $(a_j),(b_{ik}),(c_j)$ such that
$$\delta a_j=\sum_iu_iv_{ij}\quad \delta b_{ik}=\sum_jv_{ij}w_{jk}
\quad\delta c_j= \sum_kw_{jk}x_k,$$
then
$$\sum_k\langle [u_i],[v_{ij}],[w_{jk}]\rangle[x_k]-
(-1)^u\sum_i[u_i]\langle [v_{ij}],[w_{jk}],[x_k]\rangle\hbox{  contains}$$
$$(-1)^u\sum_{i,k}u_ib_{i,k}x_k-\sum_{j,k}a_jw_{jk}x_k
+(-1)^{u+v}\sum_{i,j}u_iv_{ij}c_j-(-1)^u\sum_{i,k}u_ib_{i,k}x_k=
-\delta (\sum_ja_jc_j).$$
\par We now give an example of a calculation using Massey products, which we
shall refer to in our later work.
\proclaim \Lemma 2\cddot 12. Let $p$ be a prime
not equal to 2, and let $y$ generate
$\cohp 1 {C_p}$. Then $\langle y,y,y\rangle$ is a unique element of
$\cohp 2 {C_p}$, and
$$\langle y,y,y\rangle=\cases{0&for $p>3$ \cr \beta (y) &for $p=3$.}$$
\par
\proof
We apply the `Jacobi identity' (lemma 2\cddot 11, equation (3)) and obtain
$$\langle y,y,y\rangle+\langle y,y,y\rangle+\langle y,y,y\rangle\equiv 0
\hbox{  mod $0$}.$$
Hence for $p>3$ we obtain $\langle y,y,y\rangle=0$.
$\beta (y)$ generates $\cohp 2 {C_p}$, and $\langle y,y,y\rangle$ is
defined modulo 0, so for $p=3$
it only remains to find the constant $\lambda$ in the
equation $$\langle y,y,y\rangle=\lambda \beta (y).$$
To find $\lambda$ we
resort to the definition in terms of cochains in the bar resolution.
\pra
For any $G$, and any cocycle $y\in C^1({\rm B}G;\Bbb F_p)$, the 1-cochain $a$
 defined by
$$a([g])=-{1\over2}y([g])^2 \quad\hbox{satisfies}\quad
\delta a([g|h])=y([g])y([h]),$$
so for any such $y$, a choice for $\langle y,y,y\rangle$ is
$$\langle y,y,y\rangle([g|h])=-y([g])y([h])(y([g])+y([h])).$$
Now let $G=\langle A \rangle \cong C_3$, and define $y$ by
$$y:[A^r]\mapsto r \quad\hbox{mod (3)}$$
Now define $\bar y\in C^1({\rm B}G;\Bbb Z)$ by
$$y:[A^r]\mapsto r \quad\hbox{mod (3), where $0\leq r\leq 2$}.$$
Then $b={1\over3}\delta \bar y$ is a cocycle representing $\beta(y)$, and satisfies
$$b([A^r|A^s])=\cases {0 &for $0\leq r+s\leq 2$ \cr
1 &for $3\leq r+s \leq 5$,}$$
where we require $0\leq r,s \leq 2$.
Now it remains to solve for $\lambda$ in the equation over $\Bbb F_3$:
$$\lambda b([A^r|A^s])+rs(r+s)=f(r)+f(s)-f(r+s),$$
where $f$ is some function from $\Bbb F_3$ to itself.
We obtain $\lambda =1$, as required.
\qed\par

\beginsection{Calculation of $\cohp * {P_2}$.}
\par
In this section we use the information given by the spectral sequence for the
fibration
$$\sone\longrightarrow {\rm B}P_2\longrightarrow{\rm B}\ptwo
\eqno{(1)}$$
together with the properties of Massey products stated in the previous section
to determine $\cohp * {P_2}$. To show that a certain product in
$\cohp 5 {P_2}$ is non-zero we briefly consider the spectral sequence
(discussed in lemma~2\cddot4) for the central extension
$$\langle C\rangle \longrightarrow P_2 \longrightarrow C_p^{\widebar A}
\oplus C_p^{\widebar B}.$$
\par
We now examine the spectral sequence for the fibration $(1)$. As before
$E_2^{*,*}$ is the anticommutative ring generated by $\cohp * \ptwo$
and $u$, a generator of $H^1(\sone ; \Bbb F_p)$. By naturality of
the spectral sequence for changes of coefficients, we see that
$$d_2(u)=\lambda yy'$$
for some non-zero $\lambda$.
Hence we see that
$$\eqalign{E_\infty^{*,0}&\cong {\cohp * \ptwo / yy'\cohp * \ptwo }\cr
E_\infty^{*,1}&\cong \ker \times yy' : \cohp * \ptwo
 \rightarrow \cohp * \ptwo ,}
$$
and it follows that
$E_\infty^{*,*}$ is generated by the elements $[uy],[uy'],
\allowbreak [ux],[ux'],\allowbreak [uc_2],\ldots,[uc_{p-1}],\allowbreak
y,y',\allowbreak x,x',\allowbreak c_2,\ldots,c_{p-1}, \allowbreak z$,
subject to the
relations implied by those that hold in $E_2^{*,*}$, together with the relation
$yy'=0$. $E_\infty^{*,0}$ is the subring of $\cohp * \ptwo$ given by the image
of the restriction from $\ptwo$. Throughout this section we shall use the same
symbol for an element of $\cohp * \ptwo$ and its image under the restriction.
We shall write $H^*$ for $\cohp * {P_2}$, unless there is danger of confusion.
\par
We see immediately from the spectral sequence that $$H^1H^1=\{0\},$$
that $H^2$ has order $p^4$, and that we need to introduce two new generators
in $H^2$ not in the image of the restriction. The Massey product of any three
elements of $H^1$ is defined, and is a unique element of $H^2$. We define
$$\eqalign{Y&=\langle y,y,y'\rangle\cr Y'&=\langle y',y',y\rangle,}$$
and claim that $x,x',Y$ and $Y'$ generate $H^2$. This is a consequence of the
following lemma
\proclaim \Lemma 2\cddot 13. Define
$\widebar y,d \in\cohp 1 {\langle A,C\rangle}$ and
$\widebar y',d' \in\cohp 1 {\langle B,C\rangle}$ by
$$\hbox{$\eqalign{\widebar y: A^rC^s &\mapsto r \cr d:A^rC^s &\mapsto s}\quad
\eqalign{\widebar y': B^rC^s &\mapsto r \cr d':B^rC^s &\mapsto s,}$}$$
so that
$$\widebar y=\res ^{P_2}_{\langle A,C\rangle}(y),\quad
\widebar y'=\res ^{P_2}_{\langle B,C\rangle}(y').$$
If $Y$ and $Y'$ are defined as above, then
$$\hbox{$\eqalign{\res^{P_2}_{\langle A,C\rangle}(Y)&=\widebar yd\cr
\res^{P_2}_{\langle B,C\rangle}(Y)&=0}\quad
\eqalign{\res^{P_2}_{\langle A,C\rangle}(Y')&=0\cr
\res^{P_2}_{\langle B,C\rangle}(Y')&= -\widebar y'd'.}$}$$
\par
\proof Recall from the proof of lemma 2\cddot12 that
the 1-cochain $a$ defined by
$$a([g])=-{1\over 2}y([g])^2$$
satisfies
$$\delta a([g|h])=y([g])y([h]).$$
Define a 1-cochain $b$ by
$$b([B^rA^sC^t])=-t,$$
so that
$$\delta b([B^rA^sC^t|B^{r'}A^{s'}C^{t'}])=r's=
y([B^rA^sC^t])y'([B^{r'}A^{s'}C^{t'}]).$$
Now we may choose as cocycle representing $Y$ the following:
$$Y([g|h])=y([g])({1\over 2}y([g])y'([h])-b([h])).$$
$y$ restricts to $\langle B,C\rangle $ as the zero cocycle, so
$$\res^{P_2}_{\langle B,C\rangle}(Y)=0.$$
Similarly, $y'$ restricts to $\langle A,C\rangle$ as the zero cocycle, and $b$
restricts to a cocycle representing $-d$.
Hence
$$\res^{P_2}_{\langle A,C\rangle}(Y)=\widebar y d.$$
The results for $Y'$ are proven similarly.
\qed\par
We know that
$$\hbox{$\eqalign{\res^{P_2}_{\langle A,C\rangle}(x)&=\beta(\widebar y)\cr
\res^{P_2}_{\langle B,C\rangle}(x)&=0}\quad
\eqalign{\res^{P_2}_{\langle A,C\rangle}(x')&=0\cr
\res^{P_2}_{\langle B,C\rangle}(x')&=\beta(\widebar y'),}$}$$
so it follows that $x,x',Y,Y'$ are linearly independent.
\par We are now ready to state our theorems determining the ring structure of
$\cohp * {P_2}$. We consider separately the cases $p>3$ in theorem 2\cddot14,
and $p=3$ in theorem 2\cddot15.
\proclaim \Theorem 2\cddot14. Let $p$ be greater than 3. $\cohp * {P_2}$
is generated
by elements $y,y',\allowbreak x,x',\allowbreak Y,Y',\allowbreak
X,X',\allowbreak d_4,\ldots,d_p,c_4,\ldots,c_{p-1},z$, with
$$\deg(y)=\deg(y')=1\quad\deg(x)=\deg(x')=\deg(Y)=\deg(Y')=2$$
$$\deg(X)=\deg(X')=3\quad\deg(d_i)=2i-1\quad\deg(c_i)=2i\quad\deg(z)=2p$$
$$\hbox{$\eqalign{\beta(y)&=x\cr \beta(Y)&=X}\quad
\eqalign{\beta(y')&=x'\cr \beta(Y')&=X'}$}$$
$$\beta(d_i)={\textstyle\cases{c_i &for $i<p$\cr 0 &for $i=p$}}$$
subject to the following relations:
$$yy'=0\qquad xy'=x'y$$
$$yY=y'Y'=0\qquad yY'=y'Y$$
$$Y^2=Y'^2=YY'=0$$
$$\hbox{$\eqalign{yX&=xY \cr Xy'&=2xY'+x'Y}
\qquad\eqalign{ y'X'&=x'Y'\cr X'y&=2x'Y+xY'}$}$$
$$XY=X'Y'=0\qquad XY'=-X'Y\qquad xX'=-x'X$$
$$x(xY'+x'Y)=x'(xY'+x'Y)=0$$
$$\hbox{$\eqalign{x^py'-x'^py&=0 \cr x^pY'+x'^pY&=0}\qquad
\eqalign{x^px'-x'^px=0\cr x^pX'+x'^pX=0}$}$$
$$\hbox{$\eqalign{c_iy&=\cases{0 &\phantom{f}\phantom{p}\cr -x^{p-1}y
&\phantom{f}\phantom{p}}\cr
c_ix&=\cases{0 &\phantom{f}\phantom{p}\cr -x^p &\phantom{f}\phantom{p}}\cr
c_iY&=\cases{0 &\phantom{f}\phantom{p}\cr -x^{p-1}Y
&\phantom{f}\phantom{p}}\cr
c_iX&=\cases{0 &\phantom{f}\phantom{p}\cr -x^{p-1}X
&\phantom{f}\phantom{p}}}\quad
\eqalign{c_iy'&=\cases{0 &for $i<p-1$\cr -x'^{p-1}y' &for $i=p-1$}\cr
c_ix'&=\cases{0 &for $i<p-1$\cr -x'^p &for $i=p-1$}\cr
c_iY'&=\cases{0 &for $i<p-1$\cr -x'^{p-1}Y' &for $i=p-1$}\cr
c_iX'&=\cases{0 &for $i<p-1$\cr -x'^{p-1}X' &for $i=p-1$}}$}$$
$$c_ic_j=\textstyle{\cases{0 &for $i+j<2p-2$ \cr
x^{2p-2}+x'^{2p-2}-x^{p-1}x'^{p-1} &for $i=j=p-1$}}$$
$$\hbox{$\eqalign{d_iy&=\cases{0 &\phantom{f}\phantom{p} \cr
-x^{p-1}Y &\phantom{f}\phantom{p} }\cr
d_ix&=\cases{0 &\phantom{f}\phantom{p}\cr -x^{p-1}y &\phantom{f}\phantom{p} \cr
x^{p-1}X &\phantom{f}\phantom{p}}\cr
d_iY&=0 \cr
d_iX&=\cases{0 &\phantom{f}\phantom{p}\cr -x^{p-1}Y &\phantom{f}\phantom{p}}}
\quad
\eqalign{d_iy'&=\cases{0 &for $i<p$ \cr
x'^{p-1}Y' &for $i=p$ }\cr
d_ix'&=\cases{0 &for $i<p-1$\cr -x'^{p-1}y' &for $i=p-1$ \cr
-x'^{p-1}X' &for $i=p$}\cr
d_iY'&=0 \cr
d_iX'&=\cases{0 &for $i\neq p-1$\cr -x'^{p-1}Y' &for $i=p-1$}}$}$$
$$d_id_j={\textstyle\cases{0 &for $i<p-1$ or $j<p-1$ \cr
x^{2p-3}Y-x'^{2p-3}Y'+x^{p-1}x'^{p-2}Y' &for $i=p$ and $j=p-1$}}$$
$$d_ic_j={\textstyle\cases{0 &for $i<p-1$ or $j<p-1$ \cr
x^{2p-3}y+x'^{2p-3}y'-x^{p-1}x'^{p-2}y' &for $i=j=p-1$\cr
-x^{2p-3}X+x'^{2p-3}X'-x^{p-1}x'^{p-2}X' &for $i=p$, $j=p-1$}}$$
\pra
We define $y,y'\in \cohp 1 {P_2} $ by
$$\eqalign{y:A^rB^sC^t&\mapsto r \cr y':A^rB^sC^t&\mapsto s,}$$
then define $x,x'$ by the equations
$$x=\beta(y)\qquad x'=\beta(y').$$
The equation $yy'=0$ implies that we may define unique elements $Y,Y'$ by
$$Y=\langle y,y,y'\rangle \qquad Y'=\langle y',y',y\rangle$$
and we also define
$$X=\beta(Y)\qquad X'=\beta(Y').$$
The effect of automorphisms of $P_2$ on the generators $y,y',\allowbreak
x,x',\allowbreak Y,Y',\allowbreak X,X'$ is determined by the above definitions.
An automorphism of $P_2$ which restricts to the centre as $C\mapsto C^j$ sends
$$\eqalign{d_i&\mapsto j^id_i\cr c_i &\mapsto j^ic_i \cr z&\mapsto jz.}$$
\pra
If we define $d'\in \cohp 1 {\langle B,C\rangle}$ by
$$d':B^rC^s\mapsto s,$$
then define $c'$ by
$$c'=\beta(d'),$$
then we may define
$$d_i={\textstyle\cases{\cor^{P_2}_{\langle B,C\rangle}(c'^{i-1}d')
&for $i<p-1$\cr \cor^{P_2}_{\langle B,C\rangle}(c'^{p-2}d')-x^{p-2}y
&for $i=p-1$\cr \cor^{P_2}_{\langle B,C\rangle}(c'^{p-1}d')+x^{p-2}X
&for $i=p$.}}$$
We may define $c_i$ in terms of either $d_i$ or $\chi_i\in \cohz {2i} {P_2}$,
using the equations
$$c_i=\beta(d_i)=\pi_*(\chi_i),$$
and we define $z$ by
$$z=\pi_*(\zeta).$$
We note that
$$\eqalign{\pi_*(\chi_2)&=\lambda(xY'+x'Y) \cr
\pi_*(\chi_3)&=\lambda'XX'}$$
for some non-zero $\lambda,\lambda'$.
\par
\proof We shall prove these assertions in the order in which they are stated.
Throughout this proof we adopt the convention that $\cor(\_)$ should stand for
$\cor^{P_2}_{\langle B,C\rangle}(\_)$, and similarly $\res(\_)$ should stand
for $\res^{P_2}_{\langle B,C\rangle}(\_)$.
\par
$yy'=0$ follows from the spectral sequence, as does $xy'=x'y$.
For the other relations in $H^3$, we note that
$$\eqalign{yY&=y\langle y,y,y'\rangle \cr
&\equiv\langle y,y,y\rangle y'\quad\hbox{modulo $yH^1y'=\{0\}$}\cr
&=0,}$$
and similarly for $y'Y'$.
Also
$$\eqalign{Yy'&=\langle y,y,y'\rangle y' \cr
&\equiv y\langle y,y',y'\rangle\quad\hbox{modulo $yH^1y'=\{0\}$}\cr
&\equiv y\langle y',y',y\rangle\quad\hbox{modulo $yyH^1+yy'H^1=\{0\}$}\cr
&=yY'.}$$
The spectral sequence implies that $yY'$ is non-zero.
\par
Before moving on to $H^4$ we note that
$$\langle y,y,Y\rangle \equiv 0 \quad\hbox{modulo $H^2y+H^2y'$},$$
which follows from
$$-\langle y,y,Y\rangle+\langle y,\langle y,y,y\rangle,y'\rangle
 -\langle\langle y,y,y\rangle,y,y'\rangle \equiv 0
\quad\hbox{modulo $H^2y+H^2y'$}.$$
Now
$$\eqalign{Y^2&=Y\langle y,y,y'\rangle\cr
&\equiv -\langle Y,y,y\rangle y'\quad\hbox{modulo $YH^1y'=\{0\}$}\cr
&\equiv -\langle y,y,Y\rangle y'\quad\hbox{modulo $YH^1y'+yH^2y'=\{0\}$}\cr
&\equiv 0 \qquad\hbox{modulo $H^2yy'+H^2y'^2=\{0\}$}}$$
Similarly, $Y'^2=0$ and $(Y+Y')^2=0$, because an automorphism of $P_2$
sending $y$ to $y-y'$ and fixing $y'$ sends $Y$ to $Y+Y'$, and we deduce that
$YY'=0$.
\par
The relations involving $yX$ and $y'X'$ follow by applying the Bockstein to the
relations $yY=0$ and $y'Y'=0$ respectively. For the remaining relations in
$H^4$ we introduce some matrix Massey products.
We consider
$$\langle (x,x'),\pmatrix{y' \cr -y},y\rangle,$$
which is defined modulo $H^2H^1$.
$$\eqalign{\langle (x,x'),\pmatrix{y' \cr -y},y\rangle y&\equiv
-(x,x')\pmatrix{\langle y',y,y\rangle \cr -\langle y,y,y\rangle}
\quad\hbox{modulo $xH^1y+x'H^1y=\{0\}$}\cr
&=-xY,}$$
and
$$\langle (x,x'),\pmatrix{y' \cr -y},y\rangle y' \equiv -(x,x')
\pmatrix{\langle y',y,y'\rangle \cr -\langle y,y,y'\rangle}
\quad\hbox{modulo $xH^1y'+x'H^1y'=\{0\}$}$$
but
$$\langle y',y,y'\rangle \equiv -2Y'\quad\hbox{modulo $\{0\}$}$$
so
$$\langle (x,x'),\pmatrix{y' \cr -y},y\rangle y'=2xY'+x'Y.$$
Similarly,
$$\eqalign{\langle (x',x),\pmatrix{y \cr -y'},y'\rangle y'=-x'Y' \cr
\langle (x',x),\pmatrix{y \cr -y'},y'\rangle y=2x'Y+xY',}$$
hence
$$X\in \langle (x,x'),\pmatrix{y' \cr -y},y\rangle
\quad X'\in \langle (x',x),\pmatrix{y \cr -y'},y'\rangle $$
and the remaining relations in $H^4$ follow. Our results allow us to deduce
that $X,X',\allowbreak yY',\allowbreak xy,\allowbreak xy',\allowbreak x'y'$
form a basis for $H^3$.
In $H^5$ the relations
$$XY=0\quad X'Y'=0\quad XY'=-X'Y \quad xX'=-x'X$$
follow by applying the Bockstein to the relations
$$Y^2=0\quad Y'^2=0\quad YY'=0\quad X'y=2x'Y+xY'$$
respectively. For the relations stated in $H^6$, we note that
$$xyY'=x'yY'=0,$$ and then apply the Bockstein to these relations,
noting also that
$$\beta(yY')=2xY'+2x'Y.$$
\par
We see now that in the $E_\infty^{*,*}$ page of the spectral sequence $Y$ yields
$\lambda[uy]$ and $Y'$ yields $-\lambda[uy']$ for some non-zero $\lambda$, and
deduce that
$$xY'+x'Y\equiv 0 \quad \hbox{modulo $\langle x^2,xx',x'^2,c_2\rangle$},$$
then the relations in $H^6$ imply that
$$xY'+x'Y=\lambda'c_2$$
for some $\lambda'$. $yY'$ is not in the kernel of the Bockstein, so $\lambda'$
must be non-zero, and we see that products of elements in $H^2$ generate $H^4$.
\par
It follows from lemma~2\cddot4 that
$H^5$ is generated by $x^2y,x^2y',xx'y',x'^2y',xX,xX',x'X'$ and $XY'$.
Now $XY'\notin {\rm Im}(\pi_*)$, and since $\cohz 6 {P_2}$ has
exponent $p$, we deduce that
$$XX'=-\beta(XY')\neq 0.$$
$$XX'x=0,$$
but it is apparent from the spectral sequence that in $H^6$,
$$\ker(\times x)=\langle c_3\rangle,$$
so
$$XX'=\lambda c_3\qquad \lambda\neq 0$$
and hence $H^6$ is generated by $XX'$ and products of elements of $H^2$.
\par
We already have the relations
$$x^py'=x'^py\qquad x^px'=x'^px,$$
and we obtain
$$x^pY'+x'^pY=0$$
by applying $\pone$ to the relation
$$xY'+x'Y=\lambda c_2,$$
but firstly we must find $\res(z)$. We recall that
$$z=\pi_*(\zeta)=\pi_*(c_p(\rho)),$$
and that $\rho$ restricts to $\langle B,C\rangle$ as the sum of $p$
representations, whose first Chern classes map under $\pi_*$ to
$c',c'+\widebar x',\ldots,c'+(p-1)\widebar x'$. Hence we see that
$$\eqalign{\res(z)&=\prod_{i=0}^{p-1}(c'+i\widebar x')\cr
&=c'^p-\widebar x'^{p-1}c'.}$$
Now we may calculate $\pone(c_2)$ as follows:
$$\eqalign{\pone(c_2)&=\pone(\cor(c'^2))\cr &=\cor(2c'^{p+1})\cr
&=2\cor(\res(z)c'+\res(x'^{p-1})c'^2)\cr &=2zyy'+2x'^{p-1}c_2 \cr &=0}$$
The last of these equations follows because
$$x'c_2=\pi_*(\beta\chi_2)=0.$$
Now
$$\eqalign{0&= \pone (xY'+x'Y)\cr &=x^pY'+xY'^p+x'^pY+x'Y^p \cr
&=x^pY'+x'^pY.}$$
Now we apply the Bockstein to obtain the relation
$$x^pX'+x'^pX=0.$$
\par
We may verify that in degrees greater than 6 all products of the generators
$y,y',\allowbreak x,x',\allowbreak Y,Y',\allowbreak X,X'$ may be expressed in
the form
$$\eqalign{f_1+f_2Y+f_3Y' \qquad &\hbox{for even total degree} \cr
f_1y+f_2y'+f_3X+f_4X' \qquad &\hbox{for odd total degree}}$$
where $f_i$ is a polynomial in $x$ and $x'$. With the exception of $xY'+x'Y$,
such expressions satisfy `the same' relations as elements of $H^*$ as they do
as elements of $E_\infty^{*,*}$. Elements that are expressible as above form a
subspace of $H^n$ of codimension 1 for $7\leq n\leq 2p$, so we introduce the
elements $c_i$ for $i>3$ and $z$ to our generating set.
$\delta_p$ sends elements of the above form and odd total degree to polynomials
in $\alpha $ and $\beta$, so any element mapping under $\delta_p$ to
$\chi_i$ (resp. $p\zeta$) will suffice to complete a basis for $H^{2i-1}$
(resp. $H^{2p-1}$). Hence we may add $d_i$ defined as above to our generating
set.
\par
Before verifying the relations involving $c_i$ and the low dimensional
generators, we consider the automorphism $\theta$ of $P_2$, given by
$$\eqalign{\theta : A &\mapsto B\cr B&\mapsto A,}$$
which has the effect of exchanging the `primed' and `unprimed' low dimensional
generators of $H^*$. We know the effect of this automorphism on $\chi_i$, and
we deduce that it sends
$$\theta^*:c_i\mapsto (-1)^ic_i.$$
\par
We already know the relations between $c_i$ and $y,y',x,x'$ because they are
exactly the relations that they satisfy in $\cohp * \ptwo$. For $c_iY$, note
that
$$\eqalign{\cor(c'^i)Y&=\cor(c'^i\res(Y)) \cr &=\cor(0)\cr&=0,}$$
so
$$c_iY={\textstyle\cases{0 &for $i<p-1$ \cr -x^{p-1}Y &for $i=p-1$}}$$
For the results concerning $c_iX$ we apply the Bockstein to the above results.
For the corresponding results concerning the `primed' generators, we apply the
automorphism $\theta$ described above. The results concerning $c_ic_j$ follow
from the relations in $\cohp * \ptwo$, or from those that hold in $\cohz *
{P_2}$.
\par
Now we examine the relations between $d_i$ and the low dimensional generators.
We do not yet know the effect of the automorphism $\theta$ on $d_i$, so we must
examine both the `primed' and `unprimed' relations.
$$\cor(c'^{i-1}d')y=\cor(c'^{i-1}d'\res(y))=0,$$
so
$$d_iy={\textstyle\cases{0 &for $i<p$\cr x^{p-2}Xy=-x^{p-1}Y &for $i=p$.}}$$
$$\cor(c'^{i-1}d')y'=\cor(c'^{i-1}d'\widebar y'),$$
where $\widebar y'$ is as defined during the statement of lemma 2\cddot13,
but
$$\eqalign{\cor(c'^{i-1}d'\widebar y')&=-\cor(c'^{i-1})Y'\cr
&={\textstyle\cases{0 &for $i<p$ \cr -(c_{p-1}+x^{p-1})Y'=(x'^{p-1}-x^{p-1})Y'
&for $i=p$.}}}$$
Hence
$$d_iy'={\textstyle\cases{0 &for $i<p$ \cr x'^{p-1}Y' &for $i=p$.}}$$
Similarly,
$$\eqalign{\cor(c'^{i-1}d')Y&=0  \cr
\cor(c'^{i-1}d')Y'&=\cor(c'^{i-1}d'\widebar y'd')=0,}$$
so $$\eqalign{d_iY&=0 \quad\hbox{for all $i$}\cr
d_iY'&=0 \quad\hbox{for all $i$.}}$$
For the relations involving $d_i$ and $x,x',X,X',$ we apply the Bockstein to
the above relations and substitute for terms involving $c_i$. For example,
$$0=\beta(d_{p-1}y)=c_{p-1}y-d_{p-1}x,$$
so
$$d_{p-1}x=c_{p-1}y=-x^{p-1}y.$$
\par
We now use the formula for $\res\cor$:
$$\res\cor=\sum_{j=0}^{p-1}c_{A^j}^*$$
where $c_{A^j}^*$ is the map of $\cohp * {\langle B,C\rangle}$ induced by
conjugation by $A^j$. It is easy to verify that
$$\eqalign{c_{A^j}^*: d'&\mapsto d'+j\widebar y'\cr
c'&\mapsto c'+j\widebar x'.}$$
Then
$$\res\cor(c'^i)=\sum_{j=0}^{p-1}(c'+j\widebar x')^i.$$
Similarly,
$$\res\cor(c'^{i-1}d')=\sum_{j=0}^{p-1}(c'+j\widebar x')^{i-1}(d'+j \widebar
y').$$
Hence
$$\eqalign{\cor(c'^{i-1}d')\cor(c'^{j-1}d')&=
\cor(c'^{i-1}d'\sum_{k=0}^{p-1}(c'+k\widebar x')^{j-1}(d'+k \widebar
y'))\cr &=\cor(c'^{i-1}d'\sum_{k=0}^{p-1}k\widebar y'(c'+k\widebar x')^{j-1})
\cr&=\cor(c'^{i-1}d'\sum_{l=0}^{j-1}{j-1 \choose l}\widebar y'
c'^{j-1-l}\widebar x'^l \sum_{k=0}^{p-1}k^{l+1}) }
$$
but
$$\sum_{k=1}^{p-1}k^{l+1}\equiv {\textstyle\cases{-1 &if $p-1$ divides $l+1$,
\cr 0 &otherwise}}$$
We may assume that $i>j$, so the only non-zero case we need consider is
$j=p-1$ and $l=p-2$.
$$\eqalign{\cor(c'^{p-1}d')\cor(c'^{p-2}d') &=
-\cor(c'^{p-1}d'\widebar y'\widebar x'^{p-2})\cr &=
\cor(c'^{p-1})Y'x'^{p-2} \cr&=
(c_{p-1}+x^{p-1})Y'x'^{p-2} \cr &=
x^{p-1}x'^{p-2}Y'-x'^{2p-3}Y',}$$
so
$$\eqalign{d_pd_{p-1}&=
(\cor(c'^{p-1}d')+x^{p-2}X)(\cor(c'^{p-2}d')-x^{p-2}y)\cr &=
x^{p-1}x'^{p-2}Y'+x^{2p-3}Y-x'^{2p-3}Y'}$$
and we obtain the relations between the $d_i$.
\par
Similarly
$$\cor(c'^{i-1}d')\cor(c'^j)=\cor(c'^{i-1}d'\sum_{l=0}^j{j\choose l}
c'^{j-l}\widebar x'^l\sum_{k=0}^{p-1}k^l),$$
so
$$\eqalign{\cor(c'^{i-1}d')\cor(c'^j)&={\textstyle\cases{
0 &for $j<p-1$\cr -\cor(c'^{i-1}d'\widebar x'^{p-1}) &for $j=p-1$}}
\cr&={\textstyle\cases{0 &for $i<p-1$ or $j<p-1$\cr
-(d_{p-1}+x^{p-2}y)x'^{p-1} &for $i=j=p-1$\cr
(d_p+x^{p-2}X)x'^{p-1} &for $i=p$ and $j=p-1$}}}$$
Hence we obtain the relations for $d_ic_j$ as claimed.
\par
All that remains to be calculated is the effect of automorphisms of $P_2$ on
$d_i$. The effect of automorphisms on $c_i$ and $z$ follow immediately from the
corresponding results for integral cohomology. Let $\phi$ be an automorphism of
$P_2$ that restricts to  the centre as
$$\phi:C\mapsto C^j.$$
We know that
$$\delta_p \phi^*(d_i)=\phi^*\delta_p(d_i)=j^i\delta_p(d_i),$$
so
$$\phi^*(d_i)\equiv j^id_i \qquad \hbox{modulo $\ker(\delta_p)$.}$$
Also $d_i$ was defined so that it corresponds, in the $E_\infty^{*,*}$ page of
the spectral sequence, to $\lambda[uc_{i-1}]$ for some non-zero $\lambda$.
If $j\neq \pm 1$ then $\phi$ does not extend to an automorphism of $\ptwo$, but
it can be extended to an endomorphism $\widebar \phi$ of $\ptwo$, defined on
$\sone$ by
$$\widebar \phi :v\mapsto jv \qquad\hbox{for all $v\in\sone\cong
{\Bbb R/\Bbb Z}$}.$$
It may be verified that $\widebar \phi^*(c_i)=j^ic_i$, by first verifying that
$\widebar \phi^*(\chi_i)=j^i\chi_i$. Note that $\widebar\phi^*$ is an
automorphism of $\cohp * \ptwo$ although for $j\neq \pm 1$ it is not an
automorphism of $\cohz * \ptwo$. The map induced by $\widebar\phi$;
$$\widebar\phi:{\ptwo/{P_2}}\mapsto {\ptwo/{P_2}}\cong \sone$$
has degree $j$, so $u$ is sent to $ju$ under the map of $H^*
(\sone;\Bbb F_p)$ induced by $\widebar\phi$.
\par
Hence the element $[uc_{i-1}]$ is sent to $j^i[uc_{i-1}]$ by the map of the
spectral sequence induced by $\widebar\phi$, so it follows that
$$\phi^*(d_i)\equiv j^id_i \qquad \hbox{modulo im$(\res^{\ptwo}_{P_2})$.}$$
$$\res^{\ptwo}_{P_2}(\cohp {2i-1} \ptwo )\cap\ker(\delta_p)=\{0\},$$
so we see that
$$\phi^*(d_i)= j^id_i.\qed$$
\par
\proclaim \Theorem 2\cddot15. Let $p=3$. Then
$\coht * {P_2}$ is generated by elements
$y,y',\allowbreak x,x',\allowbreak Y,Y',\allowbreak X,X',\allowbreak z$, with
$$\deg(y)=\deg(y')=1\quad\deg(x)=\deg(x')=\deg(Y)=\deg(Y')=2$$
$$\deg(X)=\deg(X')=3\quad\deg(z)=6$$
$$\hbox{$\eqalign{\beta(y)&=x\cr \beta(Y)&=X}\quad
\eqalign{\beta(y')&=x'\cr \beta(Y')&=X'}$}$$
subject to the following relations:
$$yy'=0\qquad xy'=x'y$$
$$yY=y'Y'=xy'\qquad yY'=y'Y$$
$$YY'=xx'\qquad Y^2=xY'\qquad Y'^2=x'Y$$
$$\hbox{$\eqalign{yX&=xY-xx' \cr Xy'&=x'Y-xY'
\cr XY&=x'X \cr XY'&= -X'Y }
\qquad\eqalign{ y'X'&=x'Y'-xx' \cr X'y&=xY'-x'Y
\cr X'Y'&=xX' \cr xX'&= -x'X }$}$$
$$XX'=0\qquad x(xY'+x'Y)=-xx'^2\qquad x'(xY'+x'Y)=-x'x^2$$
$$\hbox{$\eqalign{x^3y'-x'^3y&=0 \cr x^3Y'+x'^3Y&=-x^2x'^2}\qquad
\eqalign{x^3x'-x'^3x=0\cr x^3X'+x'^3X=0}$}$$
\pra
We define $y,y'\in \coht 1 {P_2} $ by
$$\eqalign{y:A^rB^sC^t&\mapsto r \cr y':A^rB^sC^t&\mapsto s,}$$
then define $x,x'$ by the equations
$$x=\beta(y)\qquad x'=\beta(y'),$$
The equation $yy'=0$ implies that we may define unique elements
in $\coht 2 {P_2}$ by forming the Massey product of any three elements of
$\coht 1 {P_2}$, and we define $Y,Y'$ by
$$Y=\langle y,y,y'\rangle \qquad Y'=\langle y',y',y\rangle$$
and we also define
$$X=\beta(Y)\qquad X'=\beta(Y'),$$
and note that $x,x'$ satisfy
$$x=\langle y,y,y\rangle \qquad x'=\langle y',y',y'\rangle.$$
The effect of automorphisms of $P_2$ on the generators $y,y',\allowbreak
x,x',\allowbreak Y,Y',\allowbreak X,X'$ is determined by the above definitions.
An automorphism of ${P_2}$ which restricts to the centre as $C\mapsto C^j$
sends
$$z\mapsto jz,$$
and we may define $z$ to be $\pi_*(\zeta)$. We also note that
$$\pi_*(\chi_2)=-xY'-x'Y-x^2-x'^2.$$
\par
\proof Many of the relations may be proven exactly as in the case $p>3$. The
relations involving $Y^2$,$Y'^2$, and $YY'$ must be proven differently, and we
exhibit a basis for $H^4$ before attempting to prove them.
\par
$yy'=0$ and $xy'=x'y$ follow from the spectral sequence. Now
$$\eqalign{yY&=y\langle y,y,y'\rangle \cr
&\equiv\langle y,y,y\rangle y'\quad\hbox{modulo $yH^1y'=\{0\}$}\cr
&=xy',}$$
and similarly for $y'Y'$. The relation $yY'=y'Y$ follows exactly as in the
proof of theorem 2\cddot14. For the relations involving $yX$ and
$y'X'$ we apply the
Bockstein to the relations $yY=xy'$ and $y'Y'=xy'$. As in the proof of
theorem 2\cddot14, any
$$Z\in \langle (x,x'),\pmatrix{y' \cr -y},y\rangle$$
satisfies
$$Zy=-xY+xx'\qquad Zy'=x'Y-xY'.$$
Similarly, any
$$Z'\in \langle (x',x),\pmatrix{y \cr -y'},y'\rangle$$
satisfies
$$Z'y'=-x'Y'+xx'\qquad Z'y=xY'-x'Y.$$
We deduce that $X,X'$ satisfy the relations claimed in $H^4$, and that
$X,X'$,$xy,xy',x'y'$,$Yy'$ form a basis for $H^3$. Using the relation for
$X'y$ we see that
$$\beta(Y'y)=-xY'-x'Y.$$
The relations
$$xY'y=x'^2y\qquad x'Y'y=x^2y'$$
follow easily from the relations we have already proven, now we apply the
Bockstein to them, and obtain
$$x(xY'+x'Y)=-xx'^2\qquad x'(xY'+x'Y)=-x'x^2.$$
It follows that the relation in $H^4$ yielding the relation in $E_\infty
^{3,1}$ of the spectral sequence,
$$x[uy']=x'[uy]$$
must be
$$0=xY'+x'Y+c_2+x^2+x'^2.$$
We deduce that $x^2,xx',x'^2,xY,xY',x'Y,$ and $x'Y'$ form a basis for $H^4$.
Now we shall return to the other relations we wish to prove in $H^4$.
\par
Consider the map
$$(\res^{P_2}_{\langle A,C\rangle},\res^{P_2}_{\langle B,C\rangle}):
\coht 4 {P_2} \longrightarrow \coht 4 {\langle A,C\rangle} \times
\coht 4 {\langle B,C\rangle}.$$
Lemma 2\cddot 13 tells us that $YY'$ is in the kernel of this map, and that
$xx'$,$xY'$ and $x'Y$ form a basis for this kernel. Hence
$$YY'=\lambda xx'+\lambda' xY'+\lambda''x'Y,\eqno{(*)}$$
for some $\lambda,\lambda',\lambda''$.
The automorphism $\theta$ of $P_2$ which exchanges the `primed' and `unprimed'
generators fixes $YY'$, so $$\lambda''=\lambda'.$$
Now we multiply equation $(*)$ by $y$, and obtain
$$xx'y=\lambda xx'y -\lambda'x'^2y.$$
Hence we see that
$YY'=xx'$.
For the remaining relations in $H^4$, we consider the effect of an automorphism
$\theta'$ of $P_2$ having the following effect on $H^1$:
$$\eqalign{\theta'^*: y&\mapsto y+y' \cr y' &\mapsto y'.}$$
We verify that
$$\eqalign{\theta'^*: x&\mapsto x+x' \cr x' &\mapsto x'\cr
Y&\mapsto Y+x'-Y'\cr Y'&\mapsto Y'+x',}$$
and apply $\theta'^*$ to the equation $YY'=xx'$, obtaining
$$YY'+x'^2-Y'^2+x'Y=xx'+x'^2,$$
and hence
$$Y'^2=x'Y.$$
We obtain the relation for $Y^2$ similarly.
\par
We apply the Bockstein to the relations
$$\eqalign{Y^2&=xY'\cr Y'^2&=x'Y \cr YY'&= xx' \cr 0&=xY'+x'Y+c_2+x^2+x'^2}$$
respectively to obtain the relations
$$\eqalign{XY&=-xX' \cr X'Y'&=-x'X\cr XY'&=-X'Y \cr xX'&=-x'X.}$$
Lemma 2\cddot 4 implies that $XY'$ is linearly
independent of elements of the form
$$f_1y+f_2y'+f_3X+f_4X',$$
where $f_i$ is a polynomial in $x$ and $x'$, and hence as in theorem 2\cddot14
we see
that no new generator is needed in $H^5$.
$XX'=\beta(YX')$ is in the image of $\pi_*$, which, in $H^6$, is spanned by
 $z,x^3$,$x^2x',xx'^2$, and $x'^3$. It is easily checked that
$$XX'x=XX'x'=0,$$
but
$$\ker(\times x)\cap\ker(\times x')\cap\langle z,x^3,x^2x',xx'^2,x'^3\rangle
=\{0\},$$
and so
$$XX'=0.$$
We already have the relations
$$x^3y'=x'^3y\qquad x^3x'=x'^3x,$$
and as in the case $p>3$, we prove the relation
$$x^3Y'+x'^3Y=-x^2x'^2$$
by applying $\pone$ to the relation
$$0=x'Y+xY'+c_2+x^2+x'^2.$$
Noting that $\res^{P_2}_{\langle B,C\rangle}(z)=c'^3-\widebar x'^2c'$,
we have
$$\eqalign{\pone(c_2+x^2+x'^2)&=\pone(\cor^{P_2}_{\langle B,C\rangle}(c'^2)+
x'^2) \cr
&= -\cor^{P_2}_{\langle B,C\rangle}(c'^4)-x'^4 \cr
&= -\cor^{P_2}_{\langle B,C\rangle}(\res^{P_2}_{\langle B,C\rangle}(z)c'+
\widebar x'^2c'^2)-x'^4 \cr
&=-zyy'-x'^2c_2+x'^2x^2-x'^4
\cr &=x'^2x^2}$$
and
$$\eqalign{\pone(xY'+x'Y)&=x^3Y'+x'^3Y+xY'^3+x'Y^3 \cr
&=x^3Y'+x'^3Y-x^2x'^2,}$$
so
$$x^3Y'+x'^3Y=-x^2x'^2.$$
We apply the Bockstein to this relation to obtain
$$x^3X'+x'^3X=0.$$
\par
It may be checked that all products of degree at least 6 of elements of degree
at most three may be expressed in the form
$$\eqalign{f_1+f_2Y+f_3Y' \qquad &\hbox{for even total degree} \cr
f_1y+f_2y'+f_3X+f_4X' \qquad &\hbox{for odd total degree}}$$
where $f_i$ is a polynomial in $x$ and $x'$. The relations we have given
between such elements are sufficient to imply the relations that hold between
the corresponding elements of the spectral sequence, hence our presentation of
the ring $\coht * {P_2}$ is complete. The effect of automorphisms on $z$
follows from its definition as $\pi_*(\zeta)$, and we already have the
required expression for $\pi_*(\chi_2)=c_2$.
\qed\par\noindent
{\bf Remarks.} Using the results of theorems~2\cddot14 and~2\cddot15, the
author has determined the differentials in the spectral sequence with
$\Bbb F_p$ coefficients for $P_2$ expressed as a central extension of $C_p$ by
$C_p$. The $E_\infty$-page in the case $p=7$ is depicted in figure~2\cddot4.
\pageinsert
 
\def\spect#1#2{\hbox{\dimen0=#1truein\vrule
\vbox{\tabskip= \dimen0\baselineskip= 31pt
\halign{\kern -.9\dimen0\rlap{$\scriptstyle##$}\hfil&&\rlap
{$\scriptstyle##$}\hfil\cr#2\crcr}\smallskip\hrule}}}

$$\spect{.9}{\cr
z&zy,zy'&zx,zx'&\sentry{zxy,zxy',\cr zx'y'}&\sentry{zx^2,zxx',\cr zx'^2}
&\sentry{zx^2y,zx^2y',\cr zxx'y',zx'^2y'}&\sentry{zx^3,zx^2x',\cr zxx'^2,zx'^3}
\cr
-\cr-\cr-&-&d_7\cr
-\cr-&-&d_6&c_6\cr
-\cr-&-&d_5&c_5\cr
-\cr-&-&d_4&c_4\cr
-&-&-\cr-&-&XY'&XX'&-\cr
-&X,X'&-&\sentry{xX,xX',\cr x'X'}&-&\sentry{x^2X,x^2X',\cr xx'X',x'^2X'}&-\cr
-&Y,Y'&yY' &\sentry{xY,xY',\cr x'Y,x'Y'}&-&\sentry{x^2Y,x^2Y',\cr xx'Y',x'^2Y'}
&-\cr
1&y,y'&x,x'&\sentry{xy,xy',\cr x'y'}&\sentry{x^2,xx',\cr x'^2}
&\sentry{x^2y,x^2y',\cr xx'y',x'^2y'}&\sentry{x^3,x^2x',\cr xx'^2,x'^3}}$$
\figure 2/4/The $E_\infty$-page for the extension
$C_7\rightarrowtail P_2\twoheadrightarrow
C_7\oplus C_7$ with $\Bbb F_7$ coefficients/
 
\endinsert
It is unreasonable to expect that by rechoosing generators we could make the
statement of theorem 2\cddot15 look more like the statement of
theorem 2\cddot 14. For
example, for $p>3$, $\cohp 2 {P_2}$ may be expressed as the direct sum of two
Aut($P_2$)-invariant subspaces:
$$\cohp 2 {P_2}=\langle x,x'\rangle\oplus\langle Y,Y'\rangle.$$
For $p=3$ the subspace $\langle x,x'\rangle$ has no Aut($P_2$)-invariant
complement.
\par
We recall that central
extensions of $C_p$ by $G$ are classified up to equivalence of extensions by
$\cohp 2 G$ (see [Br] or [Th3] for details), where $E$ is equivalent to $E'$ if
there is a homomorphism $\theta$ making the following diagram commute:
$$\matrix{C_p&\rightarrowtail&E&\twoheadrightarrow&G\cr
\mapdown{{\rm Id}}&&\mapdown{\theta}&&\mapdown{{\rm Id}}\cr
C_p&\rightarrowtail&E'&\twoheadrightarrow&G}$$
It is easily seen that a two generator group $G$ of exponent $p$ and order $p^4$
must be a central extension of $C_p$ by $P_2$. Conversely, such an extension
has exponent $p$ iff its extension class restricts to zero on all cyclic
subgroups, and in this case $G$ will be a two generator group iff its extension
class is non-zero. For $p>3$, the elements $\lambda Y+\mu Y'$ restrict to zero
on all cyclic subgroups, whereas for $p=3$ there is no such non-zero element.
We thus have verified a result due to Burnside [Bu], that there are two
generator groups of exponent $p$ and order $p^4$ only when $p>3$.
\par
The action of the Steenrod algebra $\steen$ on the 1- and 2-dimensional
generators of $\cohp * {P_2}$, and on the $c_i$ and $d_i$ (which are expressed
in terms of corestrictions from an abelian subgroup) is apparent. The generator
$z$ is the restriction of the generator of $\cohp {2p} \tilp$ of the same name,
so $\pone(z)$ is determined by lemma~2\cddot8. The following proposition
completes the description of the action of $\steen$ on $\cohp * {P_2}$.
\par
\proclaim \prop 2\cddot16. With notation as in theorems 2\cddot14 and
2\cddot15,
$$\eqalign{\pone(X)&=x^{p-1}X+zy\cr \pone(X')&=x'^{p-1}X'-zy'.}$$
\par
\proof The spectral sequence operation $_B\pone$ defined by Araki and Vasquez
([Ar], [Va]) on the $E_\infty$ page of the spectral sequence for ${\rm B}P_2$
as an \sone-bundle over ${\rm B}\tilp$ sends $ux$ to $ux^p$, so we deduce that
$\pone(X)\equiv x^{p-1}X$ modulo the image of the restriction from $\tilp$. Let
$K$ be a subgroup of $P_2$ of index $p$, and let $\cohp * K=\Bbb F_p[\bar x,c]
\otimes \Lambda[\bar y, d]$, where $\beta(d)=c$, $\beta(\bar y)=\bar x$, and
$d$, considered as a morphism from $K$ to $\Bbb F_p$ sends $C$ to 1. Then if we
let $\res=\res^{P_2}_K$ we have $\res(X)=\lambda(\bar xd-c\bar y)$, since $X$
is in the image of the Bockstein and the image of the Bockstein in $\cod 3 K$
is generated by $\bar xd-c\bar y$. We obtain
$$\res(\pone(X)-x^{p-1}X)=-\lambda(c^p-\bar x^{p-1}c)\bar y.$$
\par
If $P$ is an expression of degree $2p+1$ involving only $y,y',\allowbreak x$,
and $x'$, then $\res(P)$ is a multiple of $\bar x^p\bar y$, and if $Q$ is in
the span of $zy$ and $zy'$, then $\res(Q)$ is a multiple of $(c^p-\bar
x^{p-1}c)\bar y$. We know that $\pone(X)-x^{p-1}X=P+Q$ for some such choices of
$P$ and $Q$, and we deduce that for all $K$, $\res(P)=0$. Thus $\res(\beta
P)=0$, and hence $\beta P$ is a multiple of $x^px'-x'^px$, so is zero in
$\cohp * {P_2}$. The Bockstein is injective on the subspace of $H^{2p+1}$
generated by $x,x',\allowbreak y$, and $y'$, so we deduce that $P=0$. From
lemma~2\cddot13 we can determine $\lambda$ in the case $K=\langle A,C\rangle$
and $\bar y =\res(y)$ (resp.\ $K=\langle B,C\rangle$ and $\bar y=\res(y')$),
and conclude that $Q=-zy$. The result for $\pone(X')$ follows from this result
or may be deduced similarly. \qed\par
The above information tells us what happens in the Atiyah-Hirzebruch spectral
sequence for $P_2$.
\proclaim \Corollary 2\cddot17. In the Atiyah-Hirzebruch spectral sequence for
$P_2$ ($E_2^i\cong \cohz i {P_2}$ converging to a filtration of the
representation ring of $P_2$) the only non-zero differential is $d_{2p-1}$,
which sends $\mu$ to a multiple of $\zeta\beta$, and $\nu $ to a multiple of
$\zeta\alpha$.
\par
\proof In this spectral sequence the first potentially non-zero differential is
$d_{2p-1}$, which is $\delta_p\pone\pi_*$.\qed\par

\def\ind{{\rm Ind}}
\def\lemma{\Lemma}
\def\corollary{\Corollary}

\def\Sum{\sum}
\def\product{\prod}
\beginsection {3.} The Size of the Chern Subring and its Closure. \par
\mark{\ }
We recall that a group $G$ is said to have $p$-rank $n$ if $n$ is maximal such
that $G$ has a subgroup isomorphic to $(C_p)^n$.
In [Bl], Blackburn classified, for odd primes $p$, the $p$-groups having no
normal subgroup isomorphic to $C_p\times C_p\times C_p$. It can be shown,
independently of Blackburn's result, that all such groups satisfy the apparently
stronger condition of having $p$-rank at most two, see for example [Go].
The classification is as follows:
\proclaim{\Theorem 3\cddot 1}. (Blackburn) Let $p$ be an odd prime. Then the
$p$-groups of $p$-rank two are the following groups:
\pra\smallskip\noindent{\bf 1.}
The (non-cyclic) metacyclic $p$-groups.
\pra\smallskip\noindent{\bf 2.} $P(n)$, where $n\geq 3$ as defined in the
introduction to section 2.
\pra\smallskip\noindent{\bf 3.}
$B(n,\epsilon)$, where $n\geq 4$, there are two groups for each $n$,
depending whether $\epsilon$ is 1 or a quadratic
non-residue modulo $p$.
$B(n,\epsilon)$ has order $p^n$, and may be presented as:
$$\langle A,B,C\mid A^p=B^p=C^{p^{n-2}}=[B,C]=1\quad [A,C^{-1}]=B\quad
[B,A]=C^{\epsilon p^{n-3}}\rangle$$
For $p>3$ the groups $B(4,\epsilon)$ are the seventh and eighth groups on
Burnside's list of groups of order $p^4$, and for $p=3$, $B(4,-1)$ is the
tenth group on Burnside's list ([Bu] or the appendix to this dissertation),
and $B(4,1)$ is the eighth group.
\pra\smallskip\noindent{\bf 4.}
In the case $p=3$, every $3$-group of maximal \nilc except $C_3$, $C_9$ and the
wreath product of $C_3$ with itself (the sixth group on Burnside's list of
groups of order 81 [Bu]).
\qed
\par
Let us say that a group $G$ has property $C$ if $\ch G =\cohev G$ (recall that
$\ch G $ is the subring generated by Chern classes of representations of $G$).
Various of the $p$-groups of $p$-rank two have been shown to have property $C$,
for example in section 2 we verified Thomas' result that $P(n)$ has property
$C$ [Th2]. Thomas also verified that the split metacyclics have property $C$
and conjectured that all $p$-groups of $p$-rank two would have property $C$
[Th1], [Th2]. He also gave the example $A_4$ to show that the conjecture
could not be extended to arbitrary groups of $p$-rank at most two for all
primes. AlZubaidy claimed to have verified the conjecture for $p\geq
5$, see [Al1], [Al2], but some of his proofs are flawed.
Recently Tezuka-Yagita have shown that all metacyclic $p$-groups have property
$C$ [TY2], and Huebschmann has shown independently [Hu1] that all finite
metacyclic groups have property $C$. Moselle has
suggested that property $\widebar C$, that of satisfying $\chbar G=\cohev G$
might be more natural, and points out that $G$ has property $\widebar C$ if and
only if each of its Sylow subgroups does [Mo].
\par
We shall show that for $p\geq 5$ the groups $B(n,\epsilon)$ do not have
property $C$, but that for
\def\bthree{{B(n,\epsilon)}}
$p=3$ the groups $\bthree$ do have property $C$. We also show
that for all odd $p$ the groups $\bthree$ have property $\widebar
C$. Yagita has shown independently that for $p \geq 5$ the groups $\bthree $ do
not have property $C$ [Ya2].
\par
\mark{The Size of the Chern Subring and its Closure}
In [At] Atiyah showed that for any finite group $G$, $K^0({\rm B}G)$
is the completion of the representation ring of $G$ with respect to a certain
topology. The filtration of $K^0({\rm B}G)$ given by the $E_\infty$ page of
the Atiyah-Hirzebruch spectral sequence ($\co i G {K^j(*)} \Rightarrow
K^{i+j}({\rm B}G)$) gives rise to a filtration on the representation ring of
$G$. He conjectured that this filtration coincided with another filtration
defined algebraically, and remarked that this conjecture is equivalent to the
conjecture that $\ch G$ maps onto the $E_\infty $ page of the AHSS. Many
counterexamples have been found of composite order
[Wei]. We show that for $p\geq 5$
the groups $\bthree $ are counterexamples. These seem to be the first
counterexamples of prime power order.
\par
First we state a result (presumably well known) concerning abelian $p$-groups
and property $C$, and show that $p$-groups with an elementary abelian maximal
subgroup of rank at least 3 do not have property $\widebar C$.
\par
\proclaim \prop 3\cddot2. Let $A$ be an abelian $p$-group. Then the following
are equivalent:
\pra
1) $A$ has $p$-rank at most two.
\pra
2) $\ch A =\cohev A $.
\pra
3) $\chbar A=\cohev A$.
\par
\proof An abelian group has only 1-dimensional irreducible representations, so
$\ch A$ and $\chbar A$ are generated in degree two. For any finite group $G$,
$$\cher 2 G=\cohz 2 G \cong \hom(G,U(1)),$$
so 2) and 3) are equivalent to each other and to the statement that $\cohev A$
is generated in degree two. If $A$ is cyclic then $\cohz * A$ is a polynomial
algebra on a generator of degree two, and if $A$ is a product of two cyclic
groups then the K\"unneth theorem implies that $\cohev A$ is a polynomial
algebra on two generators of degree two. For the general case, let
$$A\cong C_{p^{n_1}}\oplus C_{p^{n_2}}\oplus\ldots\oplus C_{p^{n_m}},\qquad
\hbox{where}\quad n_1\leq n_2\leq\ldots\leq n_m.$$
$$\hbox{Then}\qquad
\cohp * A\cong \Bbb F_p[x_1,\ldots,x_m]\otimes\Lambda[y_1,\ldots,y_m],$$
where the $n_i$-th higher Bockstein maps $y_i$ to $x_i$. The image of
$\chbar A$ under reduction mod-$p$ is the subring generated by the $x_i$, and
the image of $\cohz * A$ is the universal cycles in the Bockstein spectral
sequence. (Recall that the Bockstein spectral sequence has $E_2^i\cong \cohp i
A$ and converges to $\Bbb F_p$ concentrated in degree zero, with differentials
the higher Bocksteins and $B_i^*$ the image under reduction mod-$p$ of
elements of order dividing $p^{i-1}$.) The first non-zero differential in this
spectral sequence sends $y_1y_2y_3$ to $x_1y_2y_3+\epsilon x_2y_1y_3+
\epsilon' x_3y_1y_2$, where $\epsilon,\epsilon'$ are either 0 or 1, so this
element is in the image of $\cohev A$, but not in the image of $\chbar A$.
\qed\par
\prop 3\cddot3. If $G$ is a $p$-group with a $(C_p)^n$ subgroup of index $p$
for some $n\geq 3$, then $\chbar G\neq \cohev G$.
\par
\proof The subgroup is maximal, so is normal. Call it $N$, and let $Q$ be the
quotient $G/N\cong C_p$. By proposition 3\cddot2 we may assume that $Q$ acts
non-trivially on $N$. We may identify $H^1=\cohp 1 N$ with $\hom(N,\Bbb F_p)$
and then we see that $$\cohp * N\cong S(\beta H^1)\otimes E(H^1),$$
where $S(\beta H^1)$ and $E(H^1)$ are respectively
the symmetric and exterior algebras on
the vector space $H^1$. Now $H^2$ splits as a direct sum of $\beta H^1$ and
$E^2(H^1)$, $H^3$ splits as a direct sum of $\beta H^1\otimes H^1$ and
$E^3(H^1)$, and the action of $Q$ on $H^i$ respects these splittings. Let $y$
be an element of $H^1$ with an orbit of length $p$ so that its images
under $Q$ are $y+iy'$, for some $y'$ fixed by $Q$,
and let $\alpha=\delta_p(y)$ and $\alpha'=\delta_p(y')$. Also let $w$ be an
element of $E^3(H^1)$ fixed by $Q$, and let $\xi=\delta_p(w)$.
Now writing ${\cal N}^G_N$ for Evens multiplicative
transfer map from $N$ to $G$ [Ev1], we have
$$\res^G_N{\cal N}^G_N(\alpha+\xi)=
\prod^{p-1}_{i=0}(\alpha+i\alpha'+\xi)=\alpha^p-\alpha'^{p-1}\alpha
-\alpha'^{p-1}\xi.$$
Since $\pi_*(\alpha'^{p-1}\xi)$ does not lie in $S(\beta H^1)$ it follows that
$\alpha'^{p-1}\xi$ is not in $\chbar N$, and hence the component of
${\cal N}^G_N(\alpha+\xi)$ in degree $2p+2$ cannot be in $\chbar G$.\qed\par

\beginsubsection The Groups $B(n,\epsilon)$.
Throughout this section we shall consider $\bthree$ as being presented,
as in the introduction, by elements $A$, $B$, and $C$ subject to the above
relations. We note that the centre of $\bthree$ is cyclic generated by $C^p$,
and that if we apply the circle construction to the whole centre we obtain a
\def\tilb{{\widetilde B}}
group $\tilb$ which is independent of both $n$ and $\epsilon$. This is because
there is only one circle group on $p^3$ components of nilpotency class three.
The subgroup generated by $A$, $B$, and $C^p$ is isomorphic to $P(n-1)$, with
centre also generated by $C^p$, and so we have a commutative diagram of
extensions.
$$\matrix{P(n-1) & \mapright{} & \bthree & \mapright{} & C_p \cr \mapdown{}
&& \mapdown{} && \mapdown{id} \cr
\ptwo & \mapright{} & \tilb & \mapright{} & C_p} $$
We note also that the subgroup generated by $B$ and $C$ is abelian and
isomorphic to $C_{p^{n-2}}\oplus C_p$. The corresponding subgroup of $\tilb$
(that is the one generated by $B$, $C$ and $\sone$) is also abelian and is
isomorphic to $\sone\oplus C_p\oplus C_p$.
The action induced on $\cohz * {P(n-1)}$ by conjugation by $C$ fixes $\chi_i$,
$\zeta$, $\alpha$ and $\nu$ and sends $\beta$ to $\beta+\alpha$, and $\mu$ to
$\mu+\nu$. The action induced on $\cohz * \ptwo$ is similar.
\proclaim \lemma 3\cddot4. The fixed point subring of $\cohz * \ptwo$ (resp.
$\cohev {G;\Bbb Z}$) is generated by the elements $\alpha$, $\chi_i$, $\zeta$,
and $\beta^n(\beta^p-\alpha^{p-1}\beta)$ for $n\geq 0$. \par
\proof
The element $\beta^p-\alpha^{p-1}\beta$ is the product of all conjugates of
$\beta$ so must be fixed. Also $\alpha(\beta^p-\alpha^{p-1}\beta)=0$ and
$\beta^n$ is sent to itself modulo multiples of $\alpha$, so
$\beta^n(\beta^p-\alpha^{p-1}\beta)$ is fixed. $\zeta^i P(\alpha,\beta)$ is
fixed if and only if $P(\alpha,\beta)$ is, so we only need to check that there
are no more fixed points in the subring generated by $\alpha$ and $\beta$. We
may take as basis for the degree $n$ piece of this the elements $\alpha^n,
\enspace\alpha^{n-1}\beta,\ldots\beta^n$ if $n\leq p$ and $\alpha^n,
\enspace\alpha^{n-1}\beta,\ldots\alpha^{n-p+1}\beta^{p-1},\enspace\beta^n$ if
$n\geq p$. In either case we are left to show that there are no fixed points of
the form $\Sum_{i=1}^m \lambda_i\alpha^{n-i}\beta^i$, where $m\leq p-1$.
Equating coefficients of $\alpha^{n-m+1}\beta^{m-1}$ in
$$\Sum_{i=1}^m \lambda_i\alpha^{n-i}\beta^i=\Sum_{i=1}^m \lambda_i\alpha^{n-i}
(\beta+\alpha)^i$$
we obtain
$$\lambda_{m-1}=\lambda_{m-1}+m\lambda_m$$
so $\lambda_m=0$, and inductively each $\lambda_i$ is zero. \qed
\proclaim \lemma 3\cddot5. The image of $\ch \tilb$ under restriction to
$\ptwo$ (resp. the image of $\ch \bthree $ under restriction to $P(n-1)$) is
generated by $\chi_i$, $\zeta$, $\alpha$, and $\beta^p-\alpha^{p-1}\beta$.
\par
\proof
By lemma 1\cddot 7 it suffices to prove the assertion for the Lie groups. A
direct proof for the finite groups would be very similar. $\tilb$ has an
abelian subgroup of index $p$, so has only 1- and $p$-dimensional irreducible
representations. By considering the natural isomorphisms
$$\cher 2 G = \cohz 2 G \cong \hom(G,\sone)$$
we see that $\chi_1$ and $\alpha$ span the image of $\cher 2 \tilb$.
For any $p$-dimensional representation $\rho$ of $\ptwo$ the representation
$\ind^\tilb_\ptwo(\rho)$ will split into $p$~ $p$-dimensional representations of
$\ptwo$, any of which will restrict to $\rho$. Hence $\chi_i $ and $\zeta$ are
in the image. If $\theta $ is the 1-dimensional representation of $\ptwo$ with
$c_1(\theta)=\beta$, then
$$c_p(\ind^\tilb_\ptwo(\theta) ) = \product^{p-1}_{i=0}(\beta+i\alpha)=
\beta^p-\alpha^{p-1}\beta.$$
We do not need any more generators because the image must be generated in
degrees at most $2p$. In these degrees we already have all of the fixed points
by lemma 3\cddot4. \qed\par
Let $H$ be the subgroup of $\tilb$ generated by $\sone$ and $B$, and define
elements $\tau'$, $\beta'$ of $\cohz 2 H \cong \hom(H,\sone)$ by the following
formulae.
$$ \hbox{$\eqalign{\tau'|_\sone&=1_\sone \cr
\tau'(B)&=0}\qquad \eqalign{\beta'|_\sone &=0\cr
\beta'(B)&=\exp({2\pi i /p})}$}$$
Let $M$ be the abelian maximal subgroup of $\tilb$, that is the
subgroup generated by $\sone$, $B$, and $C$. The restriction map from $\cohz *
M $ to $\cohz * H$ is clearly onto.
\proclaim \lemma 3\cddot6. If $\phi \in \cohz {2(p+n)} M$ restricts to $H$ as
$-\tau'^{p-1}\beta^{n+1}$, then $\cor^\tilb_M(\phi)$ restricts to $\ptwo$ as
$\beta^n(\beta^p-\alpha^{p-1}\beta)$. \par
\proof $M\ptwo=\tilb$, and $M\cap\ptwo= H$, so the restriction-corestriction
formula shows us that
$$\eqalignno{\res^\tilb_\ptwo\cor^\tilb_M(\phi)&= \cor^\ptwo_H \res^M_H(\phi)\cr
&= -\cor^\ptwo_H(\tau'^{p-1}\beta'^{n+1})\cr
&=-\cor^\ptwo_H(\tau'^{p-1})\beta^{n+1}\cr
&=-(\chi_{p-1}+\alpha^{p-1})\beta^{n+1}\cr
&=\beta^n(\beta^p-\alpha^{p-1}\beta).&\scriptstyle\blacksquare}$$
\par
\proclaim \corollary 3\cddot7. For $p\geq 5$, $\cohev \bthree $ is not
generated by Chern classes. \par
\proof
$\beta^{p+1}-\alpha^{p-1}\beta^2$ is in the image of the restriction from
$\bthree$
to $P(n-1)$ by lemma 3\cddot6. For $p\geq 5$ it may be verified that it is
not expressible in terms of the generators for the image of $\ch \bthree$ given
in lemma 3\cddot5. For $p=3$,
$$\beta^4-\alpha^2\beta^2=\chi_2^2-3\zeta\chi_i-\alpha^4.\eqno{{\scriptstyle
\blacksquare}}$$
\par \noindent
Yagita has also proved corollary 3\cddot7 [Ya2] using a method involving his
calculation of the Brown-Peterson cohomology of ${\rm B}P_2$ [TY].
\proclaim \corollary 3\cddot8. For $p\geq 5$ and the groups $\bthree$, the
Chern subring does not map surjectively to the $E_\infty$ page of the Atiyah
Hirzebruch spectral sequence.\par
\proof
Let $K$ be the abelian maximal subgroup of $\bthree$. Then by lemma 3\cddot5
and lemma 1\cddot 7 there is a $\phi \in \cher {2p+2} K $ satisfying
\def\btwo{{P(n-1)}}
$$\res^\bthree_\btwo\cor^\bthree_K (\phi)=\beta^{p+1}-\alpha^{p-1}\beta^2.$$
From the naturality properties of the AHSS it follows that $\cor(\phi)$ is a
universal cycle. We show that it is not the sum of a Chern class and a
universal boundary by considering the restriction to the subgroup generated by
$B$.
$$\cohz * {\langle B \rangle} = \Bbb Z[{\widebar \beta}]/(p{\widebar\beta}),$$
where $\widebar \beta$ is the restriction from $\btwo $ of $\beta$. The AHSS
for $\langle B\rangle$ collapses because $E_2^{i,j}$ is trivial if $i$ or $j$
is odd. Hence the universal boundaries $B_\infty(\bthree)$ restrict trivially
to $\langle B\rangle$. By lemma 3\cddot5 the image of $\ch \bthree$ under this
restriction is generated by ${\widebar\beta}^p$ and ${\widebar\beta}^{p-1}$.
So for $p\geq 5$,
$$\res^\bthree_{\langle B\rangle}\cor^\bthree_K (\phi)={\widebar\beta}^{p+1}
\notin  \res^\bthree_{\langle B\rangle}(B_\infty(\bthree)+\ch \bthree ).
\eqno{{\scriptstyle \blacksquare}}$$
\par
\proclaim \corollary 3\cddot9.
$\coh \tilb$ restricts onto the fixed point subring of $\coh \ptwo$.
\par
\proof
Immediate from lemma 3\cddot4 and lemma 3\cddot6.\qed\par
\proclaim \corollary 3\cddot10. $\cohev \tilb$ is generated by corestrictions
of Chern classes, and by Chern classes for $p=3$. \par
\proof
Consider the spectral sequence with integer coefficients for the extension
$$\tilp \mapright{} \tilb\mapright{}C_p.$$
For any extension with quotient $C_p$ the $E_2$ page of the corresponding
spectral sequence is generated by $E_2^{0,*}$, $E_2^{1,*}$ and $E_2^{2,0}\cong
\Bbb F_p$. The inflation of a generator for $\cod 2 {C_p}$ is a Chern class,
and yields a generator for $E_2^{2,0}$. Also corollary~3\cddot9 implies that
$E_2^{0,*}$ consists of universal cycles, and lemma~3\cddot6 implies that
corestriction of Chern classes (resp.\ Chern classes in the case when $p=3$)
yield generators for the even degree subring of $\eee 2, *, *,$.
\qed\par
We shall now work towards theorem~3\cddot13, which will tell us that
$\cohev {B(n,\epsilon)}$ is generated by corestrictions of Chern classes, and
will involve the application of theorem~1\cddot9. Corollary~3\cddot10 has
verified one of the hypotheses required for this, so it suffices now to check
that multiplication by the extension class $c(B(n,\epsilon))$ in $\cod 2
\tilb$ is injective on $\cohod \tilb$. This is proven for all $n$ as a
corollary of lemma~3\cddot12, but we first present a simpler proof in the case
when $n\geq 5$. The reader who objects to redundancy should turn directly to
the end of the proof of lemma~3\cddot11.
\par
Since $\tilp$ has index $p$ in $\tilb$ and trivial odd degree cohomology it
follows by the usual corestriction-restriction argument that $\cohod \tilb$ has
exponent $p$. It is easily checked that $c(B(n,\epsilon))$ is of the form
$p^{n-4}$ times an element of infinite order plus a non-zero element of the
image of the inflation from ${\tilb/\tilp}\cong C_p$. Thus it will follow that
$\cohev {B(n,\epsilon)}$ is generated by corestrictions of Chern classes (resp.
Chern classes for $p=3$) for $n\geq 5$ provided that multiplication by a
generator for $\inf(\cohz 2 {\tilb/\tilp})$ is injective on $\cohod \tilb$.
This is a case of the following lemma.
\proclaim \Lemma 3\cddot11. Consider a fibration
$$F\mapright{}E\mapright{\pi}{\rm B}C_p,$$
such that $H^*(F)$ is concentrated in even degrees and $H^*(E)$ maps onto
$H^*(F)^{C_p}$, and let $\gamma$ be a generator for $\cod 2 {C_p}$. Then
multiplication by $\pi^*(\gamma)$ is injective on $H^{\rm od}(E)$.
\par
\proof Consider the spectral sequence for the fibration. The element
$\pi^*(\gamma)$ yields a generator for $\eee 2, 0, 2, $ which we shall call
$\gamma$, and all elements in even total degree are universal cycles. It
suffices to prove that multiplication by $\gamma $ is an isomorphism from
$\eee \infty, 2i+1, j, $ to $\eee \infty, 2i+3, j,$ for all $i$ and $j$.
In fact we shall show inductively
that multiplication by the corresponding element of $\eee 2n,2,0,$ (which we
shall also call $\gamma $) is an isomorphism from $\eee 2n,i,j,$ to $\eee
2n,i+2,j,$ provided that either $i$ is odd or $i \geq 2n$. This is true for
$n=1$. Now $d_{2n}$ is trivial because it maps between even and odd values of
$j$. It is possible for $d_{2n+1}$ to be non-trivial, but only from odd total
degree to even total degree. Since by induction the odd degree part of $\eee
2n+1,*,*,$ is generated by $\eee 2n+1,1,*,$, $d_{2n+1}$ is completely
determined by its values on $\eee 2n+1,1,*,$. We obtain the following
commutative diagram, where the horizontal maps are isomorphisms by induction.
$$\matrix{\eee 2n+1,1,j,&\mapright{\times \gamma}&\eee 2n+1,3,j,
&\mapright{\times \gamma}&\eee 2n+1,5,j,&\mapright{}&\cdots\cr
\mapdown{d_{2n+1}}&&\mapdown{d_{2n+1}}&&\mapdown{d_{2n+1}}\cr
\eee 2n+1,2n+2,j-2n-1,&\mapright{\times\gamma}&
\eee 2n+1,2n+4,j-2n-1,&\mapright{\times\gamma}&
\eee 2n+1,2n+6,j-2n-1,&\mapright{}&\cdots\cr}$$
Hence multiplication by $\gamma$ induces isomorphisms on the kernels
and cokernels, and the inductive step is proven.\qed
\par
We shall now deduce that $\cohev G=\chbar G$ (resp. $\ch G$ for $p=3$) for any
$p$-group~$G$ of $p$-rank two and nilpotency class three using the following
lemma.
\proclaim \Lemma 3\cddot12. As above, let $M$ be the abelian subgroup of
$\tilb$ of index $p$. Then the restriction from $\cohod \tilb$ to $\cohod M$ is
injective. \par
\proof
First let us fix a presentation for $\tilb$.
$$\eqalign{
\tilb=\langle \sone,A,B,C\mid \hbox{$\sone$ central,}\qquad
A^p&=B^p=C^p=1\cr
[A,B]&=C\qquad
[B,C]=1\qquad [A,C]=e^{2\pi i/p}\rangle}$$
Note that $\tilb/\sone$ is generated by the images of $A$ and $B$, which we
shall call $A$ and $B$ too, so that our notation coincides with that of
section~2. Note that the subgroup $M$ is expressible as an extension of $\sone$
by $\langle B,C\rangle$. We shall show that the map of spectral sequences
induced by the following map of extensions is injective on the $E_\infty$ page.
$$\matrix{\relax\sone&\mapright{}&M&\mapright{}&\langle B,C\rangle\cr
\mapdown{\rm Id}&&\mapdown{i}&&\mapdown{\bar i}\cr
\relax\sone&\mapright{}&\tilb&\mapright{}&{P_2}}$$
Let $\eee *,*,*,$ stand for the spectral sequence for $\tilb$, and let
$\widebar{E}_*^{*,*}$ for the spectral sequence for $M$. Notice that
$\widebar{E}_*^{*,*}$ collapses. If we use the notation of section~2 for
elements of $\cohz * {P_2}$ and let $\cohz * \sone \cong \Bbb Z[\tau]$, then we
see that
$$\eee 2,*,*,\cong \Bbb Z[\tau]\otimes \cohz * {P_2}.$$
If we let $\rho$ be a 1-dimensional representation of $M$ such that
$\res^M_\sone(c_1(\rho))=\tau$, then
$$\res^\tilb_\sone(c.(\ind^\tilb_M(\rho)))=(1+\tau)^p.$$
Hence $p\tau^i$ and $\tau^p$ are universal cycles, the $E_\infty$ page is
periodic vertically with period $2p$, and $\eee \infty,*,*,\cong \eee 2p,*,*,$.
Since $\tilb'\cap\sone$ has order $p$ there can be no homomorphism from $\tilb$
to $\sone$ that restricts to $\sone$ as an isomorphism, so $d_3(\tau)$ must be
non-zero. The kernel of the map from $\eee 2,3,0, $ to $\widebar{E}_2^{3,0}$
is generated by $\nu$, and $d_3(\tau)$ must lie in this kernel (because
$\widebar{E}_*^{*,*}$ collapses), so without loss of generality
$d_3(\tau)=\nu$.
\par
It may be verified that $\eee 5,*,*,$ is generated by $\alpha$, $\beta$,
$\chi_2$, $\chi_4,\ldots,\chi_{p-1}$, $\zeta$, $\mu$, $\tau^i\beta\mu$,
$\tau^i\chi_2$, $\tau^i\chi_4,\ldots,\tau^i\chi_{p-2}$,
$\tau^i(\chi_{p-1}+\alpha^{p-1})$, $p\tau^i$, $\tau^p$, and $\tau^{p-1}\nu$,
where $1\leq i\leq p-1$, subject to various relations. For example,
we have that $[\tau^i\beta\mu]\alpha=0$ unless $i$ is congruent to $-1$
modulo~$p$.
Figure~3\cddot1 depicts this $E_5$ page in the case when $p\geq 7$.
\midinsert
 
\def\speck#1#2{\hbox{\dimen0=#1truein\vrule
\vbox{\tabskip= \dimen0\baselineskip= 35pt
\halign{\kern -.9\dimen0\rlap{$\scriptstyle##$}\hfil&&\rlap
{$\scriptstyle##$}\hfil\cr#2\crcr}\smallskip\hrule}}}
 
$$\speck{.7}{\cr p\tau^3 &-&-&- & \tau^3\chi_2 &\tau^3\beta\mu&
-&\tau^3\beta^2\mu&\tau^3\chi_4\cr
-&-&-&-&-&-&-&-&-\cr
p\tau^2 &-&-&- & \tau^2\chi_2 &\tau^2\beta\mu&-&\tau^2\beta^2\mu&\tau^2\chi_4\cr
-&-&-&-&-&-&-&-&-\cr
p\tau &-&-&-& \tau\chi_2 &\tau\beta\mu&-&\tau\beta^2\mu&\tau\chi_4\cr
-&-&-&-&-&-&-&-&-\cr
1 &-& \alpha,\beta &\mu &\sentry{\alpha^2,\alpha\beta,\cr \beta^2,\chi_2}&
\beta\mu & \sentry{\alpha^3,\alpha^2\beta,\cr \alpha\beta^2,\beta^3}
&\beta^2\mu&\sentry{\alpha^4,\alpha^3\beta,\cr\alpha^2\beta^2,\alpha\beta^3,\cr
\beta^4,\chi_4}
}$$
\figure 3/1/The $E_5$-page of the spectral sequence of lemma 3\cddot12 in the
case $p\geq 7$/
\endinsert
We now claim that $\tau^{p-1}\nu$ survives until the $E_{2p-1}$ page, and that
$d_{2p-1}(\tau^{p-1}\nu)$ is a non-zero multiple of $\zeta\alpha$ (which is
equal to $\delta_p\pone\pi_*(\nu)$---see proposition~2\cddot16 and
corollary~2\cddot17). To show this, we first consider the spectral sequence
with $\Bbb Z_{(p)}$ coefficients for the path-loop fibration over
$K(\Bbb Z,3)$. Let $i_3$ be the (image in $\Bbb Z_{(p)}$ coefficients of the)
fundamental class in $H^3(K(\Bbb Z,3))$, and $i_2$ that in
$H^2(\Omega K(\Bbb Z,3))$, so that $d_3(i_2)=i_3$. It is readily verified
that $i_2^{p-1}i_3$ survives until the $E_{2p-1}$ page, and that
$d_{2p-1}(i_2^{p-1}i_3)$ is a unit multiple of $\delta_p\pone\pi_*(i_3)$. The
claim concerning $\tau^{p-1}\nu$ now follows by naturality, because $\nu$ may
be thought of as a map from ${\rm B}P_2$ to $K(\Bbb Z,3)$ classifying the
$K(\Bbb Z,2)$ bundle ${\rm B}\tilb$.
\par
Even on $E_5$ pages the map from $\eee 5,{\rm od},2i,$ to $\widebar{E}_5^{{\rm
od},2i}$ is injective except when $i$ is congruent to $-1$ modulo $p$, and its
kernel is the submodule generated by $\tau^{p-1}\nu$. Thus it will suffice to
show that no differential lower than $d_{2p-1}$ can hit anything in the ideal
of $\eee 5,*,0,$ generated by $\zeta\alpha$. This follows because
multiplication by $\alpha$ is injective on this ideal but trivial on
$\eee 5,{\rm od},2i,$ for $1\leq i<p-1$.\qed
\par
\proclaim \Theorem 3\cddot13. If $G$ is any finite normal subgroup of $\tilb$ of
$p$-rank two with quotient isomorphic to \sone,
then $\cohev G=\chbar G$, and if $p=3$ then $\cohev G=\ch G$.
Examples of such $G$ include the $p$-groups $B(n,\epsilon)$ and for $p=3$ the
seventh, eighth and tenth groups on Burnside's list of groups of order 81 ([Bu]
or the appendix). \par
\proof
Let $c(G)$ be the extension class of $G$ in $\cod 2 \tilb$,
that is the element of $\cod 2 \tilb$ corresponding to the homomorphism from
$\tilb$ to $\sone$ with kernel $G$. Since $H=G\cap M$
is abelian of $p$-rank at most two, $\ch H=\cohev H$, so by theorem~1\cddot9
multiplication by $c(H)$ is injective on $\cohod M$. But now
$c(H)=\res^\tilb_M(c(G))$, and so by lemma~3\cddot12 multiplication by $c(G)$
is injective on $\cohod \tilb$. The result now follows from theorem~1\cddot9,
since corollary~3\cddot10 tells us that $\cohev \tilb$ is generated by
corestrictions of Chern classes and by Chern classes for $p=3$.
The groups listed in the statement do have $p$-rank two, and they
are normal subgroups of $\tilb$, because they have nilpotency class three, and
a central cyclic subgroup with quotient group $P_2$, so yield $\tilb$ when the
circle construction is applied. \qed
\par
\proclaim \corollary 3\cddot14. For $p\geq 5$, if $G$ is a $p$-group of
$p$-rank two, then $\cohev G$ is generated by corestrictions of Chern classes.
\par
\proof This is immediate from a combination of Blackburn's classification [Bl]
(or see theorem~3\cddot1), Huebschmann's results on the metacyclic groups
[Hu1], Thomas' result for the groups $P(n)$ [Th2] (or theorem~2\cddot3), and
theorem~3\cddot13.\qed\par

\beginsection {4.} Further Calculations at the Prime Three. \par \mark{\ }
In this section we shall make a more detailed study of the integral cohomology
of the group $\tilb$ and the $p$-groups $B(n,\epsilon)$, particularly when
$p=3$. In this case the proof of lemma~3\cddot12 describes most of the
behaviour of a spectral sequence converging to a filtration of $\cohz * \tilb$.
Using this spectral sequence we determine the cohomology of~$\tilb$, and then
calculate the cohomology of $B(n,\epsilon)$ using the spectral sequence for
${\rm B}B(n,\epsilon)$ as an $\sone$-bundle over ${\rm B}\tilb$. The main
result of the section is that (for $p=3$) the integral cohomology ring of
$B(n,1)$ is isomorphic to that of $B(n,-1)$, so that there are two groups of
order $3^5$ with isomorphic integral cohomology rings. The author believes that
these are the first examples of such $p$-groups, although non-isomorphic groups
of composite order with isomorphic cohomology rings have been known for some
time, and Larson has even exhibited metacyclic examples [La]. For $p>3$ it
seems far more difficult to determine the cohomology of $\tilb$. Nevertheless,
and in spite of the title of this section, we present a proof that the
cohomology groups (ignoring the ring structure) of $B(n,\epsilon)$ are
independent of $\epsilon$ for all odd $p$, without actually determining these
groups. One response to the main result of this section has been ``Are you sure
that the groups are different?'', so before starting our calculations we show
that as $\epsilon$ varies in $\Bbb F_p^\times$ there really are two
non-isomorphic groups of the form $B(n,\epsilon)$.
\par
\proclaim \Lemma 4\cddot1. For any fixed odd prime $p$ and $n\geq 5$, define
for each $\epsilon \in \Bbb F_p$ a group $G(\epsilon)$ of order $p^n$ by the
presentation
$$\eqalign{
G(\epsilon)=\langle A,B,C,D\mid A^p=B^p=C^p=1\quad C\hbox{ central,}\quad
B&=[A,D]\cr\quad C=[A,B]\quad [B,D]&=1\quad D^{p^{n-3}}=C^\epsilon\rangle.}$$
Then there are three isomorphism classes of such $G$, depending as
$\epsilon $ is 0, a quadratic residue, or a quadratic non-residue mod-$p$.
\par
\proof Note that the relations given suffice to express any element of $G$
in the form $A^iB^jC^kD^l$. Note that the subgroup of $G$ generated by $A$ and
$B$ is normal, and isomorphic to $P_2$, with quotient cyclic of order
$p^{n-3}$. First we claim that if $g\in \langle A,B\rangle$, then the order of
$gD^j$ is at least the order of $D^j$. First note that
$$(gD^j)^p=(\prod_{k=0}^{p-1}g^{D^{jk}})D^{jp},$$
and that the expression on the right is not 1 unless $D^{jp}$ is. Now since
$D^p$ commutes with elements of $\langle A,B\rangle$, we see that
$$(gD^j)^{p^2}=D^{jp^2}.\eqno{(*)}$$
It follows that if $\epsilon\neq 0$, then $\langle A,B\rangle $ is the
subgroup of $G$ of elements of order $p$, and we deduce that $G(0)$ is not
isomorphic to $G(\epsilon)$ for $\epsilon\neq 0$ because it contains more
elements of order $p$. From now on, let us consider only the case
$\epsilon\neq 0$.
\par
\mark{Further Calculations at the Prime Three}
The elements $A$ and $D$ generate $G$, and any automorphism of $G$ must
the form
$$\hbox{$\eqalign{D\mapsto&A^iB^jC^kD^l\cr A\mapsto&A^rB^sC^t}\quad
\eqalign{\hbox{where }&l\neq 0\hbox{ mod $p$}\cr \hbox{where }&1\leq r\leq
p-1.}$}$$
If we define $D'$ to be the image under such an automorphism of $D$, $A'$ to be
the image of~$A$, $B'=[A',D']$ and $C'=[A',B']$, then the primed elements
define a presentation of the form $G(\epsilon')$. To complete the proof it
suffices to show that $\epsilon'r^2=\epsilon$. From equation $(*)$ it follows
that $D'^{p^{n-3}}=D^{lp^{n-3}}$. Now
$$\eqalign{B'&=[A^rB^sC^t,A^iB^jC^kD^l]\cr
&\equiv[A^r,D^l] \quad\hbox{mod $\langle C\rangle$}\cr
&\equiv B^{rl} \quad\hbox{mod $\langle C\rangle$},\cr
C'&=[A^rB^sC^t,B^{rl}C^{?}]=C^{r^2l},}$$
from which it is apparent that $D'^{p^{n-3}}=C'^{\epsilon/{r^2}}$.\qed\par
We shall now begin our calculation of $\cohz * \tilb$ in the case when $p=3$ by
considering the spectral sequence for $\tilb$ expressed as an extension of
$\sone$ by $P_2$. This spectral sequence was studied for arbitrary $p$ during
the proof of lemma~3\cddot12, where $d_3$ and part of $d_{2p-1}$ were
calculated. In the case when $p=3$ not much more work is required to complete
the spectral sequence. For the rest of this chapter we shall consider only the
prime three unless otherwise stated.
\proclaim \Lemma~4\cddot2. Consider the spectral sequence with integral
coefficients for the extension
$$\sone\longrightarrow \tilb \longrightarrow P_2.$$
Use notation as in the proof of lemma~3\cddot12, so that
$$E_2^{*,*}\cong \Bbb Z[\tau]\otimes\cohz * {P_2}.$$
Then $d_3(\tau)=\nu$, and $E_5^{*,*}$ is generated by $\alpha$, $\beta$,
$\chi_2$, $\zeta$, $\mu$, $3\tau$, $\tau(\chi_2+\alpha^2)$, $\tau\beta\mu$,
$3\tau^2$, $\tau^2\nu$, $\tau^2(\chi_2+\alpha^2)$, $\tau^2\beta\mu$, and
$\tau^3$. All of these are universal cycles except for $\tau^2\nu$ and
$\tau^2\beta\nu$, which satisfy
$$d_5(\tau^2\nu)=\zeta\alpha\qquad
d_5(\tau^2\beta\mu)=\zeta(\beta^2+\alpha^2+\chi_2).$$
There are no other non-zero differentials, so $E_\infty^{*,*}$ is generated by
$\alpha$, $\beta$,
$\chi_2$, $\zeta$, $\mu$, $3\tau$, $\tau(\chi_2+\alpha^2)$, $\tau\beta\mu$,
$3\tau^2$, $\tau^2(\chi_2+\alpha^2)$,
$\tau^3$ and $\tau^2(\beta^2-\alpha^2)\mu$.\par
\proof The claims concerning $d_3$ and $E_5^{*,*}$ were proved in
lemma~3\cddot12. Of the generators required for $E_5^{*,*}$, all are universal
cycles with the possible exception of $\tau^2\nu$, $\tau^2(\chi_2+\alpha^2)$
and $\tau^2\beta\mu$. It was shown in lemma~3\cddot12 that
$d_5(\tau^2\nu)=\zeta\alpha$. For the remaining cases we again (as in
lemma~3\cddot12) compare the spectral sequence $E_*^{*,*}$ with the spectral
sequence $\widebar E_*^{*,*}$ for $M$ (the abelian maximal subgroup of $\tilb$)
using the map induced by the following map of extensions.
$$\matrix{\relax\sone&\mapright{}&M&\mapright{}&\langle B,C\rangle\cr
\mapdown{\rm Id}&&\mapdown{i}&&\mapdown{\bar i}\cr
\relax\sone&\mapright{}&\tilb&\mapright{}&{P_2}}$$
Define generators $\beta'$, $\gamma'$ for $\cohz 2 {\langle B,C\rangle}
\cong \hom({\langle B,C\rangle},\sone)$ by
$$\hbox{$\eqalign{\beta':B&\mapsto \omega\cr C&\mapsto 1}\qquad
\eqalign{\gamma':B&\mapsto 1\cr C&\mapsto \omega,}$}$$
where $\omega=\exp({2\pi i/3})$, so that $\beta'=\res(\beta)$. Also define a
generator $\mu$ for $\cohz 3 {\langle B,C\rangle}$ by $\mu'=\res(\mu)$. Now we
know that $\widebar E_*^{*,*}$ collapses and that $\widebar E_2^{*,*}\cong
\Bbb Z[\tau',\beta',\gamma']\otimes \Lambda[\mu']$.
\par
To show that $d_5(\tau^2(\chi_2+\alpha^2))=0$, we note that it lies in
$\eee 5,9,0,$, which is generated by $\beta^3\mu$ and $\zeta\mu$, so maps
injectively to $\widebar E_5^{9,0}$ (recall that
$\res(\zeta)=\gamma'^3-\beta'^2\gamma'$). The result follows because $\widebar
E_*^{*,*}$ collapses. Similarly, we know that $d_5(\tau^2\beta\mu)$ must be in
the kernel of $\res^{*,*}_*$. The equation
$(\tau^2\nu)\beta^2=(\tau^2\beta\mu)\alpha$ implies that
$$\alpha d_5(\tau^2\beta\mu)=\zeta\alpha\beta^2,$$
so $d_5(\tau^2\beta\mu)=\zeta\beta^2+Q$ for some $Q$ in the kernel of
multiplication by $\alpha$ as a map from $\eee 5,10,0,$ to $\eee 5,12,0,$. This
kernel is spanned by $\zeta(\alpha^2+\chi_2)$ and $\beta^5-\alpha^2\beta^3$,
which map under $\res_*^{*,*}$ to $(\beta'^2\gamma'-\gamma'^3)\beta'^2$ and
$\beta'^5$ respectively. It follows that $d_5(\tau^2\beta\mu)$ is as claimed.
The $E_\infty$-page of this spectral sequence is depicted in figure~4\cddot1
(on page~87), but with elements named by the elements of $\cohz * \tilb$
that yield them, rather than as in this statement.
\qed\par
\proclaim \Corollary 4\cddot3. In the spectral sequence of lemma~4\cddot2, the
subring of $\eee \infty,*,*,$ generated by $\tau^3$ and $\zeta$ is of the form
$\Bbb Z[x,y]/(3y)$, where $x=\tau^3$ and $y=\zeta$. As a module for this
subring $\eee \infty,*,*,$ is generated by 1, $3\tau$, $3\tau^2$,
$\alpha^{i+1}$, $\beta^{i+1}$, $\alpha\beta^{i+1}$, $\alpha^2\beta^{i+1}$,
$\chi_2$, $\beta^i\mu$, $\tau\beta^{i+1}\mu$, $\tau(\alpha^2+\chi_2)\beta^i$,
$\tau^2(\beta^2-\alpha^2)\beta^i\mu$ and $\tau^2(\alpha^2+\chi_2)\beta^i$,
where $i\geq 0$. Each of these module generators has order three, except that
1, $3\tau$ and $3\tau^2$ have infinite order. Multiplication by $\tau^3$ is
injective on $\eee \infty,*,*,$, and the kernel of multiplication by $\zeta$ is
generated (as a $\Bbb Z[\tau^3,\zeta]$-module) by $3\tau$, $3\tau^2$,
$\alpha^{i+1}$, $\alpha\beta^{i+1}$, $\alpha^2\beta^{i+1}$ and
$\beta^2+\chi_2$.
\par
\proof Almost immediate from lemma~4\cddot2. \qed\par
An immediate consequence of lemma~4\cddot2 is that $\cohz * \tilb$ may be
generated by nine elements of even degree and three of odd degree, since $\eee
\infty,*,*,$ can be. It has, however, already been shown in corollary~3\cddot10
that six elements suffice to generate $\cohev \tilb$. We shall now define six
elements of even degree and one of odd degree which (as we shall demonstrate
in lemma~4\cddot7) generate $\cohz * \tilb$. To determine the relations
that they satisfy we shall have to calculate their restrictions to $M$,
which are listed as lemma~4\cddot6. Information given by restriction to proper
subgroups does not quite suffice to show that the elements generate $\cohz *
\tilb$, essentially because there are elements of degree four that restrict
trivially to every proper subgroup. We are forced, therefore, to prove
separately in lemma~4\cddot5 that the elements of even degree defined below
generate $\cohev \tilb$, using the same method as in corollary~3\cddot10.
\proclaim \Definition 4\cddot4. Define $\tau',\beta',\gamma'\in \cohz 2 M \cong
\hom(M,\sone)$ by
$$\eqalign{\tau':zB^iC^j&\mapsto z\cr
\beta':zB^iC^j&\mapsto \omega^i\cr
\gamma':zB^iC^j&\mapsto \omega^j,}$$
where $\omega=\exp({2\pi i/3})$. Now define elements $\alpha$, $\beta$,
$\gamma$, $\delta_1$, $\delta_2$, $\delta_3$, $\mu$ of $\cohz * \tilb$ as
follows. Let $\alpha$, $\beta$ and $\mu$ be the images under inflation from
$\tilb/\sone=P_2$ of the elements with the same names. Let $\rho$ be a
3-dimensional irreducible representation of $\tilb$, whose restriction to $M$
contains the 1-dimensional representation with first Chern class $\tau'$, and
define $$\delta_3=c_3(\rho).$$
Define also
$$\eqalignno{\delta_1&=\cor^\tilb_M(\tau')\cr
\delta_2&=\cor^\tilb_M(\tau'^2)\cr
\gamma&=\cor^\tilb_M(\tau'^2(\beta'-\gamma'))=
\delta_2\beta-\cor^\tilb_M(\tau'^2\gamma').&{\scriptstyle\blacksquare}}$$
\par
\proclaim \Lemma 4\cddot5. The elements $\alpha$, $\beta$, $\gamma$,
$\delta_1$, $\delta_2$ and $\delta_3$ defined in 4\cddot4 generate $\cohev
\tilb$. \par
\proof If we consider $\beta$ as a homomorphism from $\tilb $ to $\sone$, it
has image $C_3$ and kernel $\tilp$, with $C_3$ acting on $\tilp$ via outer
automorphisms. This is, modulo choice of generators, the extension considered
in corollary~3\cddot10. Recall that in corollary~3\cddot10 it was shown (for
$p=3$) that Chern classes generate $\cohev \tilb$. The method involved checking
that their restrictions to $\tilp$ generate $\cohz * \tilp^{C_3}$. By a similar
argument we may show that the elements of the statement generate $\cohev
\tilb$ by showing that the restrictions of $\alpha$, $\gamma$, $\delta_1$,
$\delta_2$ and $\delta_3$ to $\tilp$ (that is, to $\ker(\beta)$) generate $\cohz
* \tilp^{C_3}$, which was described in lemma~3\cddot4. Some caution is required
since the notation used in lemma~3\cddot4 for elements of $\cohz * \tilp$
clashes with our current notation (which was designed to harmonise with our
notation for elements of $\cohz * {{\tilb/\sone}}$ rather than with the
cohomology of any subgroup). We claim that the restrictions are as decribed
below, where the elements of $\cohz * \tilb$ are as defined in~4\cddot4, and their
images are named as in lemma~3\cddot4.
$$\res(\alpha)=\alpha\qquad\res(\delta_1)=\chi_1\qquad
\res(\delta_2)=\chi_2+\alpha^2$$
$$\res(\gamma)=\beta^3-\alpha^2\beta\qquad \res(\delta_3)=\zeta.$$
For the elements of degree two this may be checked by considering them as
homomorphisms from $\tilb$ to $\sone$. For the elements defined as
corestrictions from $M$ to $\tilb$ we use the formula $\res^\tilb_\tilp
\cor^\tilb_M=\cor^\tilp_{M\cap\tilp}\res^M_{M\cap\tilp}$, which holds because
$\tilp M=\tilb$. For example,
$$\eqalign{\res^\tilb_\tilp(\delta_2\beta-\cor^\tilb_M(\tau'^2\gamma'))&=
0-\cor^\tilp_{M\cap\tilp}(\tau'^2\gamma')\cr
&=-\cor^\tilp_{M\cap\tilp}(\tau'^2)\beta\cr
&=-(\chi_2+\alpha^2)\beta\cr
&=\beta^3-\alpha^2\beta.}$$
For $\delta_3$, we note that the representation $\rho$ restricts to $\tilp$ as
the (unique) irreducible representation of degree three that restricts to
$\sone $ as three copies of the identity, so $c_3(\rho)$ must restrict to
$\zeta$. It is clear now from lemma~3\cddot4 that the restrictions of $\alpha$,
$\gamma$, $\delta_1$, $\delta_2$ and $\delta_3$ generate $\cohz * \tilp^{C_3}$,
and the result follows. \qed\par
\proclaim \Lemma 4\cddot6. The restrictions to the subgroup $M$ of $\tilb$ of
the elements defined in~4\cddot4 are as follows, where the equation for
$\res(\mu)$ is used to define $\mu'$, a generator for $\cohz 3 M\cong \Bbb F_3$.
$$\res(\alpha)=0\quad\res(\beta)=\beta'\quad\res(\mu)=\mu'\quad
\res(\delta_1)=3\tau'+\beta'\quad\res(\gamma)=\gamma'^3-\beta'^2\gamma'$$
$$\res(\delta_2)=3\tau'^2-\tau'\beta'+\beta'^2-\gamma'\beta'-\gamma'^2\qquad
\res(\delta_3)=\tau'^3+\tau'^2\beta'-\tau'\gamma'\beta'-\tau'\gamma'^2$$
\par
\proof The claims for $\alpha$, $\beta$ and $\mu$ are trivial. For $\delta_1$,
$\delta_2$ and $\gamma$, we use the restriction-corestriction formula, which in
this case tells us that $\res^\tilb_M\cor^\tilb_M$ is equal to the sum of
conjugation by the distinct powers of $A$. By considering them as homomorphisms
from $\tilb$ to $\sone$, it is easily verified that conjugation by $A$ fixes
$\beta$ and sends $\tau'$ to $\tau'-\gamma'$ and $\gamma'$ to $\gamma'-\beta'$.
For example,
$$\res^\tilb_M(\delta_1)=\sum_{i=0}^2 c_{A^i}^*(\tau')=
\tau'+(\tau'-\gamma')+(\tau'+\gamma'+\beta')=3\tau'+\beta'.$$
For $\delta_3$, which we recall was defined as $c_3(\rho)$, we note that since
the restriction to $M$ of $\rho$ contains a summand with Chern class $\tau'$
and is invariant under conjugation by $A$ it must also contain summands with
Chern class $\tau'-\gamma'$ and $\tau'+\gamma'+\beta'$. Thus
$$\res(\delta_3)=\tau'(\tau'-\gamma')(\tau'+\gamma'+\beta').\qed$$
\par
\proclaim \Lemma 4\cddot7. The seven elements defined in~4\cddot4 generate
$\cohz * \tilb$. Their relation to the generators for the $E_\infty$-page of
the spectral sequence of lemma~4\cddot2 is as follows.
\item{1.} The generators $\alpha$, $\beta$ and $\mu$ yield the spectral
sequence elements of the same name.
\item{2.} The generator $\delta_1$ (resp. $\delta_2$, $\delta_3$) yields the
generator $3\tau$ (resp. $3\tau^2$, $\tau^3$) for $\eee \infty,0,2,$ (resp. for
$\eee \infty,0,4,$, $\eee \infty,0,6,$).
\item{3.} The subgroup of $\cohz 4 \tilb$ generated by $\alpha^2$,
$\alpha\beta$, $\beta^2$ and $\delta_1\beta$ is mapped bijectively to $\eee
\infty,4,0,$ (which is generated by $\alpha^2$, $\alpha\beta$, $\beta^2$ and
$\chi_2$).
\item{4.} The subgroup of $\cohz 6 \tilb$ generated by $\alpha^3$,
$\alpha^2\beta$, $\alpha\beta^2$, $\beta^3$ and $\gamma$ is mapped
bijectively to $\eee \infty,6,0,$ (which is generated by $\alpha^3$,
$\alpha^2\beta$, $\alpha\beta^2$, $\beta^3$ and $\zeta$).
\item{5.} $\delta_2\beta$ (resp. $\delta_2\mu$) yields the generator
$\tau(\alpha^2+\chi_2)$ for $\eee \infty,4,2,$ (resp. $-\tau\beta\mu$ for
$\eee \infty,5,2,$).
\item{6.} $\delta_2^2-\delta_3\delta_1$ (resp. $\delta_2^2\mu$)
yields the generator
$-\tau^2(\alpha^2+\chi_2)$ for $\eee \infty,4,4,$ (resp.
$\tau^2(\beta^2-\alpha^2)\mu$ for $\eee \infty,7,4,$). \pra\noindent
The $E_\infty$-page of this spectral sequence is depicted in figure~4\cddot1.
\par
\proof To show that the elements generate $\cohz * \tilb$ it will suffice to
show that they yield a set of generators for the $E_\infty$-page of the
spectral sequence of lemma~4\cddot2. Thus it suffices to prove the numbered
assertions made above. Many of these assertions follow easily from
lemma~4\cddot6. For example the restriction to $M$ of $\delta_2$ is $3\tau'^2$
modulo terms of lower degree in $\tau'$, the restriction from
$\eee \infty,0,4,$ to $\widebar\eee\infty,0,4,$ is injective, and hence
$\delta_2$ must yield the generator $3\tau^2$. Similarly $\gamma$ yields an
element of $\eee\infty,6,0,$ which restricts to $\widebar\eee\infty,6,0,$ as
$\gamma'^3-\beta'^2\gamma'$, while the image of the subgroup of
$\eee\infty,6,0,$ generated by $\alpha^3$, $\alpha^2\beta$, $\alpha\beta^2$ and
$\beta^3$ in $\widebar\eee\infty,6,0,$ is generated by~$\beta'^3$, so $\gamma$
is independent of this subgroup, and so the five elements of item~4 above
yield generators for $\eee \infty,6,0,$. All the assertions except that of
item~3 can be proved similarly, and we will consider them proved.
\par
\midinsert

$$\spec{.64}{\cr \delta_3&-&\delta_3\alpha,\delta_3\beta&
\delta_3\mu&{\scriptstyle\ldots}\cr
-&-&-&-&-&-&-&-&-&-\cr
\delta_2&-&-&-&\delta^2_2-3\delta_3\delta_1&-&\delta_2\gamma,\delta_2^2\beta&
\delta_2^2\mu&\delta_2^2\beta^2&\delta_2^2\beta\mu\cr
-&-&-&-&-&-&-&-&-&-\cr
\delta_1&-&-&-&\delta_2\beta&\delta_2\mu&\delta_2\beta^2&\delta_2\beta\mu&
\delta_2\beta^3&\delta_2\beta^2\mu\cr
-&-&-&-&-&-&-&-&-&-\cr
1&-&\alpha,\beta&\mu&\sentry{\alpha^2,\alpha\beta,\cr \beta^2,\delta_1\beta}&
\beta\mu&\sentry{\alpha^3,\alpha^2\beta,\cr \alpha\beta^2,\beta^3,\cr
\gamma}&\beta^2\mu&\sentry{\alpha^4,\alpha^2\beta^2,\cr \alpha\beta^3, \beta^4,
\cr \gamma\beta}&\beta^3\mu,\gamma\mu}$$
 
\figure 4/1/The $E_\infty$-page of the spectral sequence of lemmata 4\cddot2 and
4\cddot7/
\endinsert
 
It is clear that $\delta_1\beta$ yields an element of $\eee\infty,4,0,$. To
prove that it does not lie in the subgroup of $\eee\infty,4,0,$ generated by
$\alpha^2$, $\alpha\beta$ and $\beta^2$ is more difficult. From the spectral
sequence it can be seen that $\cohz 4 \tilb$ is isomorphic to $\Bbb Z\oplus
C_3^4$. The products of the elements of definition~4\cddot4 that have degree
four are $\delta_2$, $\delta_1^2$, $\delta_1\beta$, $\delta_1\alpha$,
$\alpha^2$, $\alpha\beta$ and $\beta^2$. In the group generated by these
elements the subgroup of elements of order at most three is generated by
$\delta_1^2-3\delta_2$, $\delta_1\beta$, $\delta_1\alpha$,
$\alpha^2$, $\alpha\beta$ and $\beta^2$. We now determine some relations
between these elements. By Fr\"obenius reciprocity
$$\delta_1\alpha=\cor(\tau')\alpha=\cor(\tau'\res(\alpha))=0.$$
Similarly,
$$\eqalign{\delta_1^2=\cor(\tau')\cor(\tau')&=\cor(\tau'\res\cor(\tau'))\cr
&=\cor(3\tau'^2+\tau'\beta')\cr
&=3\cor(\tau'^2)+\cor(\tau')\beta'\cr
&=3\delta_2+\delta_1\beta.}$$
It follows that $\delta_1\alpha$ and $\delta_1^2-3\delta_2$ are dependent on
$\delta_1\beta$, $\alpha^2$, $\alpha\beta$, $\beta^2$. Lemma~4\cddot5 tells us
that $\alpha$, $\beta$, $\gamma$, $\delta_1$, $\delta_2$ and $\delta_3$
generate $\cohev \tilb$, so it must be that $\delta_1\beta$ is not contained
in the subgroup generated by $\alpha^2$, $\alpha\beta$ and $\beta^2$. \qed\par
\proclaim \Theorem 4\cddot8. $\cohz * \tilb$ is generated by the elements of
definition~4\cddot4 subject to the following relations.
$$3\alpha=3\beta=0\qquad 3\mu=0\qquad 3\gamma=0$$
$$\hbox{$\eqalign{\delta_1^2&=\delta_1\beta+3\delta_2\cr
\delta_1\mu&=\beta\mu}\qquad\eqalign{\gamma\delta_1&=\gamma\beta\cr
\delta_1\alpha&=0}\qquad\eqalign{\alpha\mu&=0\cr \gamma\alpha&=0}$}$$
$$\hbox{$\eqalign{\delta_1\delta_2&=\delta_2\beta+9\delta_3\cr
\delta_2\alpha&=0}\qquad\eqalign{\delta_1\beta^2&=\beta^3-\alpha^2\beta \cr
\alpha^3\beta&=\alpha\beta^3}$}$$
$$\delta_2^3-27\delta_3^2=-\delta_3(\beta^3-\alpha^2\beta)+\delta_2^2\beta^2+
\delta_2\beta^4-\gamma^2-\beta^6+\alpha^2\beta^4$$
\par
\proof If we assume that these relations hold, it follows easily from
corollary~4\cddot3 and lemma~4\cddot7 that they are all the relations that we
require. (Note in particular that the relations given imply that
$3(\delta_2^2-3\delta_3\delta_1)=0$.) Thus it remains to demonstrate that these
relations hold.
\par
That the elements $\alpha$, $\beta$, $\gamma$ and $\mu$ have order three is
immediate from their definitions. The relations claimed in degree~4 were
checked during the proof of lemma~4\cddot7. The relations $\alpha\mu=0$
and $\delta_1\mu=\beta\mu$ hold because the restrictions of these equations to
$M$ are valid, and the restriction is injective from $\eee\infty,5,0,$ to
$\widebar\eee\infty,5,0,$. The relations $\delta_2\alpha=0$, $\gamma\alpha=0$,
and $\gamma\delta_1=\gamma\beta$ follow readily by a Fr\"obenius reciprocity
argument. The relation for $\delta_1\delta_2$ may also be deduced by using
Fr\"obenius reciprocity as follows.
$$\eqalign{
\delta_2\delta_1=\cor(\tau'^2)\delta_1&=\cor(\tau'^2(3\tau'+\beta'))\cr
&=\cor(3\tau'^3)+\cor(\tau'^2)\beta \cr
&=\cor(\res(3\delta_3))+\delta_2\beta\cr
&=9\delta_3+\delta_2\beta.}$$
The relation $\alpha^3\beta=\alpha\beta^3$ holds because the corresponding
relation holds in the spectral sequence. We know (because $\delta_1\beta$
yields an element of $\eee\infty,4,0,$) that $\delta_1\beta^2$ is in the span
of $\alpha^3$, $\alpha^2\beta$, $\alpha\beta^2$, $\beta^3$ and $\gamma$, that
its restriction to $M$ is $\beta'^3$, and that it is annihilated by $\alpha$.
Since the kernel of the restriction from $\eee\infty,6,0,$ is spanned by
$\alpha^3$, $\alpha^2\beta$ and $\alpha\beta^2$ whereas the kernel of
multiplication by $\alpha$ is spanned by $\gamma$ and $\beta^3-\alpha^2\beta$,
these facts suffice to determine $\delta_1\beta^2$.
\par
It is clear from corollary~4\cddot3 and lemma~4\cddot7 that all we are left to
do is obtain an expression for $\delta_2^3$. It may be checked that the
relation given is correct when restricted to $M$, and that when it is
multiplied by $\alpha$ the valid equation $0=0$ is obtained. The result
follows, since it may be checked that in degree~12 the kernel of restriction to
$M$ is generated by $\delta_3\alpha^3$, $\delta_3\alpha^2\beta$,
$\delta_3\alpha\beta^2$, $\alpha^6$, $\alpha^2\beta^4$ and $\alpha\beta^5$, and
that multiplication by $\alpha$ is injective on this subgroup. \qed\par
\proclaim \Lemma~4\cddot9. For any odd prime~$p$, define $\delta',\beta\in\cohz
2 \tilb \cong \hom(\tilb,\sone)$, where $\tilb$ is presented as in the proof of
lemma~3\cddot12, by
$$\eqalign{\delta':A^jB^kC^lz&\mapsto z^p\cr
\beta: A^jB^kC^lz&\mapsto \exp(2\pi ik/p).}$$
(In the case when $p=3$, $\beta$ is the element $\beta$ defined in~4\cddot4 and
$\delta'$ is $\delta_1-\beta$.) then for any~$\epsilon$, the subgroup of $\tilb$
with extension class $p^{n-4}\delta'-\epsilon\beta$ is isomorphic to
$B(n,\epsilon)$.\par
\proof The assertions relating $\delta'$ and $\beta$ to $\delta_1$ and $\beta$
in the case $p=3$ are easily proved. If we write
$\phi=p^{n-4}\delta'-\epsilon\beta$, then as a homomorphism from $\tilb$ to
$\sone$, $\phi$ has the following form.
$$\phi:A^jB^kC^lz\mapsto z^{p^{n-3}}\omega^{-\epsilon j},$$
where $\omega=\exp({2\pi i/p})$. The kernel of $\phi$ (which we recall is one
of the equivalent interpretations of the group with characteristic class
$\phi$) is the subgroup generated by $A$, $C$, and $B'=B\eta^\epsilon$, where
$\eta=\exp({2\pi i/p^{n-2}})$. It is easily checked that the map from this
group to $B(n,\epsilon)$ (as presented in the statement of theorem~3\cddot1)
sending $A$ to $A$, $C$ to $B^{-1}$ and $B'$ to $C$ is an isomorphism. \qed
\par
Before calculating $\cohz * {B(n,\epsilon)}$ we require one more lemma which
will determine a generator for $\cohz 5 {B(n,\epsilon)}$.
\proclaim \Lemma~4\cddot10. Let $A$ be the abelian maximal subgroup of
$B(n,\epsilon)$, so that $A\cong C_{3^{n-2}}\oplus C_3$ and
$A=B(n,\epsilon)\cap M$. Also take generators $\widebar\tau$,
$\widebar\gamma$ and $\widebar\mu$ for $\cohz * A$, where $\widebar\tau=
\res^M_A(\tau')$, $\widebar\gamma=\res^M_A(\gamma')$, and
$\widebar\mu=\res^M_A(\mu')$. Note that $\widebar\gamma$ and $\widebar\mu$
have order $3$ while $\widebar\tau$ has order $3^{n-2}$. then there is an
element of $\cohz 5 {B(n,\epsilon)}$ that restricts to $A$ as
$\widebar\gamma\bar\mu$, and this element is in the image of the
corestriction from a $P(n-1)$ subgroup to $B(n,\epsilon)$.
\par
\proof First we claim that if $H$ is the subgroup of $A$ isomorphic to
$C_{3^{n-3}}\oplus C_3$, and $\widetilde\mu$ is a generator for
$\cohz 3 H$, then $\cor^A_H(\widetilde\mu)=\widebar\mu$. To show this,
note that if $y$ is a generator for $\coht 1 {C_{3^{n-3}}}$ (which is
isomorphic to $\hom({C_{3^{n-3}}},C_3)$), and we consider
$C_{3^{n-3}}$ as a subgroup of $C_{3^{n-2}}$, then the classical description
of the transfer (as a map from $G/G'$ to $H/H'$ where $H\subset G$) allows
us to check that $\cor(y)$ is a generator for $\coht 1 {C_{3^{n-2}}}$.
Now if we define $y\in \coht 1 H$ to be the inflation of a generator for
$\coht 1 {H/(1\oplus C_3)} $, and $y'$ the inflation of a generator for
$\coht 1 {H/(C_{3^{n-2}}\oplus C_3)}$ then Fr\"obenius reciprocity implies
that $\cor^A_H(yy')=\cor^A_H(y)y'$ is an element of $\coht 2 A$ not in the
image of $\cohz 2 A$. The claim now follows since corestriction commutes
with the Bockstein.  \par
Now note that $H=A\cap P(n-1)$, and that from the above paragraph we see
that we may take generators $\widetilde\tau$, $\widetilde\gamma$ and
$\widetilde\mu$ for $\cohz * H$, where $\widetilde\tau=
\res(\widebar\tau)$ and $\widetilde\gamma=\res(\widebar\gamma)$ while
$\widebar\mu=\cor(\widetilde\mu)$. If we use the notation of
theorem~2\cddot3 and lemma~3\cddot4 for elements of the group
$\cohz * {P(n-1)}$,
then after a suitable choice of generators the action of
$B(n,\epsilon)/P(n-1)\cong C_3$ on $\cohz * {P(n-1)}$ may be taken to
fix $\alpha$, $\chi_1$ and $\nu$, and send $\beta$ to $\beta+\alpha$ and
$\mu$ to $\mu+\nu$. The kernel of the restriction fom $P(n-1)$ to $H$ in
degree two is generated by $\alpha$ (this may be checked by considering
$\cohz 2 {P(n-1)}\cong \hom(P(n-1),\sone)$), which implies that the kernel
of the restriction from $P(n-1)$ to $H$ in degree three is generated by
$\nu$. Thus we see that
$$\eqalignno{\res^{B(n,\epsilon)}_A\cor^{B(n,\epsilon)}_{P(n-1)}(\beta\mu)&=
\cor^A_H\res^{P(n-1)}_H(\beta\mu)\cr
&=\pm\cor^A_H(\widetilde\gamma\widetilde\mu)\cr
&=\pm\widetilde\gamma\widetilde\mu&{\scriptstyle\blacksquare}}$$
\proclaim \Theorem 4\cddot11. For $n\geq 5$ (and $p=3$), $\cohz *
{B(n,\epsilon)}$ is generated by elements $\alpha$, $\delta_1$, $\mu$,
$\delta_2$, $\nu$, $\delta_3$ and $\gamma$ of degrees 2, 2, 3, 4, 5, 6 and 6
respectively. The generators other than $\nu$ may be taken to be the
restrictions from $\tilb$ of the generators for $\cohz * \tilb$ of the same
name. They are subject to the following relations.
$$3\alpha=0\qquad 3\gamma=0\qquad 3\mu=0\qquad 3\nu=0$$
$$3^{n-3}\delta_1=0\qquad 3^{n-2}\delta_2=0\qquad 3^{n-1}\delta_3=0$$
$$\hbox{$\eqalign{3\delta_2&=\delta_1^2(1-3^{n-4}\epsilon)\cr
9\delta_3&=\delta_1\delta_2(1-3^{n-4}\epsilon)\cr
\delta_2^3&=27\delta_3^2-\gamma^2}\qquad\eqalign{
\delta_1\alpha&=0\cr \delta_1\mu&=0\cr \gamma\delta_1&=0}\qquad\eqalign{
\alpha\mu&=0\cr \alpha\delta_2&=0 \cr \gamma\alpha&=0}$}$$
$$\delta_1\nu=0\qquad\mu\nu=0\qquad\delta_2\nu=-\gamma\mu\qquad\gamma\nu=
\delta_2^2\mu$$
\par
\proof We consider the spectral sequence for ${\rm B}B(n,\epsilon)$ as an
\sone-bundle over ${\rm B}\tilb$. It follows from lemma~4\cddot9 that the
differential in this spectral sequence sends $\xi$ (a generator for $\eee
2,0,1,$) to $\epsilon\beta-3^{n-4}\delta_1$ (because $n\geq 5$ we have that
$3^{n-4}\delta_1=3^{n-4}\delta'$). The relations given between the generators
that are in the image of the restriction from $\tilb$ are the relations of
theorem~4\cddot8 after the substitution of $\beta$ by
$3^{n-4}\epsilon\delta_1$. In the spectral sequence $\eee\infty,*,1,$
corresponds to the kernel of multiplication by $\epsilon\beta-3^{n-4}\delta_1$,
which is the ideal of $\cohz * \tilb$ generated by
$\delta_1\beta-\beta^2+\alpha^2$. This ideal is a free module of rank one for
$\Bbb F_3[\delta_3,\alpha]$, and the product of
$\delta_1\beta-\beta^2+\alpha^2$ with any of $\delta_1$, $\delta_2$, $\gamma$
or $\mu$ is zero. It follows that we only need one new generator $\nu$ in
degree five, and that the only extra relations that we need are expressions for
$\delta_1\nu$, $\delta_2\nu$, $\gamma\nu$ and $\mu\nu$ in terms of the
generators in the image of restriction from $\tilb$. We may take $\nu$ as in
lemma~4\cddot10, so that its restriction to $A$ (the subgroup of
$B(n,\epsilon)$ of lemma~4\cddot10) is $\widebar\gamma\bar\mu$. The
restrictions to $A$ of the other generators may be determined easily from
lemma~4\cddot6 (which describes $\res^\tilb_M$), if we note that
$\res^M_A(\tau')=\widebar\tau$, $\res^M_A(\mu')=\widebar\mu$,
$\res^M_A(\gamma')=\widebar\gamma$ and
$\res^M_A(\beta')=3^{n-4}\epsilon\widebar\tau$. Since $\res^\tilb_M$ is
injective on $\cohod \tilb$ (see lemma~3\cddot12) it follows that
$\res^{B(n,\epsilon)}_A$ is injective on the image of $\cohod \tilb$. To check
the relations given for $\delta_1\nu$, $\delta_2\nu$ and $\gamma\nu$ it
therefore suffices to check that they yield valid relations when restricted to
$A$. The relation $\mu\nu=0$ follows by Fr\"obenius reciprocity, because $\nu$
may be taken to be a corestriction from $P(n-1)$, and the product of any
element of $\cohz 5 {P(n-1)}$ and any element of $\cohz 3 {P(n-1)}$ is zero.
\qed \par
\proclaim \Corollary~4\cddot12. For $n\geq 5$ and $p=3$, the $p$-groups
$B(n,1)$ and $B(n,-1)$ defined in theorem~3\cddot1 are non-isomorphic, but have
isomorphic integral cohomology rings.
\par
\proof It was shown in lemma~4\cddot1 that these groups are distinct. We may
define an isomorphism from $\cohz * {B(n,-1)}$ to $\cohz * {B(n,1)}$ by sending
each of the generators (as listed in theorem~4\cddot11) for $B(n,-1)$ to the
generator of the same name for $B(n,1)$, except that we must send $\delta_2$ to
$(1-3^{n-4})\delta_2$ and $\delta_3$ to $(1+3^{n-4})\delta_3$. To check that
this does describe an isomorphism, it helps to observe that for either group,
the equations
$$3\delta_2^2=\delta_2\delta_1^2(1-3^{n-4}\epsilon)=9\delta_3\delta_1$$
imply that $\delta_2^2$ has order $3^{n-4}$. \qed\par
For $p\geq 5$ it is much harder to  calculate $\cohz * \tilb$. Nevertheless we
are able to prove the following theorem comparing the additive structure of
$\cohz * {B(n,\epsilon)}$ for arbitrary~$p$.
\proclaim \Theorem 4\cddot13. For all $m\geq 0$, for all odd $p$, and all
$n\geq 6$, the group $\cohz m {B(n,\epsilon)}$ is independent of $\epsilon$.
\par
\proof First we consider the groups $\cohz m \tilb$. The $E_2$-page of the
spectral sequence for $\tilb$ as an extension of $\sone$ by $P_2$ has the
property that $\eee 2,0,2j,$ is isomorphic to $\Bbb Z$, whereas $\eee 2,i,j,$
is a finite $p$-group for $(i,j)$ not of the form $(0,2j')$. It follows that
$\eee \infty,*,*,$ also has this property, and hence that $\cohz {2m}\tilb$ is
an extension of a finite $p$-group by $\Bbb Z$ (which must of course be split)
and $\cohz {2m+1} \tilb$ is a finite $p$-group. We know that $\tilb$ has a
subgroup isomorphic to $\tilp$ of index $p$, and that the torsion in
$\cohz * \tilp$ in degree $2m$ (resp.\ in degree $2m+1$) has exponent $p$
(resp.\ is trivial). Since $\cor^\tilb_\tilp\res^\tilb_\tilp$ is multiplication
by $p$ it follows that in $\cohz {2m} \tilb$ the torsion has exponent at most
$p^2$, and that $\cohz {2m+1} \tilb$ has exponent at most $p$ (the odd degree
case could also be deduced from lemma~3\cddot12). Now fix for each $m$ an
isomorphism
$$\cohz {2m} \tilb\cong T^{2m}\oplus\Bbb Z,$$
where $T^{2m}$ is the torsion subgroup of $\cohz {2m} \tilb$.
\par
Now we consider the spectral sequence for ${\rm B}B(n,\epsilon)$ as an
\sone-bundle over ${\rm B}\tilb$. With notation as in lemma~4\cddot9, we see
that the differential $d_2$ from $\eee 2,*,1,$ to $\eee 2,*,0,$ corresponds to
multiplication by $p^{n-4}\delta'-\epsilon\beta$ from $\cohz * \tilb$ to
itself. During theorem~3\cddot13 it was shown that $\cohev \tilb$ maps on to
$\cohev {B(n,\epsilon)}$, or equivalently that $\eee \infty,{\rm od},1,$ is
trivial. If we let $\phi$ stand for the map
$$\phi:T^{2m}\oplus\Bbb Z\longrightarrow T^{2m+2}\oplus \Bbb Z$$
corresponding to multiplication by $p^{n-4}\delta'-\epsilon\beta$ from $\cohz
{2m} \tilb$ to $\cohz {2m+2} \tilb$, then
$$\eee\infty,2m+2,0,\cong {(T^{2m+2}\oplus\Bbb Z)/\im(\phi)}\qquad
\eee\infty,2m,1,\cong \ker(\phi)$$
$$\eee\infty,2m+3,0,\cong
{\cohz {2m+3} \tilb/(p^{n-4}\delta'-\epsilon\beta)\cohz {2m+1}\tilb}.$$
Since $\cohod \tilb$ has exponent $p$ and $n>5$, multiplication by
$p^{n-4}\delta'$ annihilates the group
$\cohz {2m+1} \tilb$. Also $-\epsilon$ is a unit
modulo $p$, so
$$\eee \infty,2m+3,0,\cong{\cohz {2m+3} \tilb/\beta\cohz {2m+1} \tilb}.$$
Writing elements of $\cohz {2m} \tilb$ as `row vectors' via the isomorphism
with $T^{2m}\oplus \Bbb Z$, it may be checked that the map $\phi$ has matrix
$$\phi=\pmatrix{-\epsilon\beta& 0\cr p^{n-4}\delta'& p^l}:
T^{2m}\oplus\Bbb Z\longrightarrow T^{2m+2}\oplus \Bbb Z,$$
where $l$ is either $n-4$, $n-3$, or $n-2$, but depends only on $m$, not on
$\epsilon$. In any case, if $p^{n-4}$ annihilates $T^{2m}$, which certainly
happens for $n\geq 6$, then $\ker(\phi)$ and $T^{2m+2}\oplus \Bbb Z/\im(\phi)$
are independent of $\epsilon$.
\par
To complete the proof, it will suffice to show that $\cohz {2m+3}
{B(n,\epsilon)}$ is isomorphic to $\eee\infty,2m+2,1,\oplus \eee\infty,2m+3,0,$
for any choice of $\epsilon$. The group $\eee \infty,2m+3,0,$ is a subgroup of
$\cohz {2m+3} {B(n,\epsilon)}$, and it follows from lemma~3\cddot12 that this
subgroup maps injectively to $\cohz {2m+3} A$, where $A$ is the abelian maximal
subgroup of $B(n,\epsilon)$. Since $A$ is isomorphic to
$C_{p^{n-2}}\oplus C_p$, the exponent of $\cohz {2m+3} A$ is $p$, and hence if
$\xi$ is any element of $\cohz {2m+3} {B(n,\epsilon)}$ such that $p\xi$ is an
element of $\eee\infty,2m+3,0,$, then $p\xi=0$. \qed\par

\def\split{\colon}
\def\nonsplit{{\hbox{\cddot}}}
\beginsection 5. The Integral Cohomology of the Held Group.\par \mark{\ }
In this section we shall examine integral cohomology of some finite groups by
considering the restriction maps to the Sylow subgroups, and applying the
results of section~2. Groups that we shall consider include extensions of
$(C_p)^2$ by various subgroups of $GL_2(p)$, and the Held group, a sporadic
simple group of order $2^{10}.3^3.5^2.7^3.17$ [He], [Co]. For the Held group we
obtain explicitly the $p$-torsion subring of the integral cohomology for
all odd primes $p$. C.~B.~Thomas has
suggested [Th4] that the Atiyah-Hirzebruch-Rector spectral sequence [Re]
might be a useful tool in studying the modular representation rings of groups.
Since the $E_2$ page of this spectral sequence for the $p$-modular ring
consists of the cohomology of the group with coefficients of exponent prime to
$p$ our results for the Held group may have consequences concerning its
2-modular representation ring. \par
\noindent
{\bf Methods.}\par \noindent
If $M$ is any $p$-local module for a finite group $G$ (that is a module for
which multiplication by any prime other than $p$ is an isomorphism) then the
restriction map from the cohomology of $G$ with coefficients in $M$
to that of its Sylow $p$-subgroup $G_p$ is split injective,
with splitting map a multiple of the corestriction. There is a characterisation
of the image of this restriction due to Cartan-Eilenberg [CE] as the `stable
elements' of ${\co * {G_p} M}^{N(G_p)}$, that is the elements $\xi$ satisfying
$$\res^{gHg^{-1}}_{H\cap gHg^{-1}}c^*_g(\xi)=\res^H_{H\cap gHg^{-1}}(\xi).$$
Recent work of P.~Webb [We] gives, for each of a wide variety of $G$-simplicial
complexes, an alternative characterisation of the image of $\co * {G} M$ in
$\co * {G_p} M$ as the kernel of a map from the direct sum of the cohomology
of the vertex stabilisers to that of the edge stabilisers. If $G$ never
reverses an edge of the complex the maps involved are composed solely of
restriction maps, so this characterisation can be simpler than the usual one.
All the proofs given in this section use the stable element characterisation,
but others were used in some of the author's earlier proofs of these results.
The author is indebted to J.~C.~Rickard for
many discussions concerning Webb's work. \par
If $G_p$ (the Sylow subgroup of $G$) is abelian and we consider only trivial
coefficients then a theorem of Swan [Sw] tells us that the stable elements are
exactly the fixed points of $\cod * {G_p}$ under the action of its normaliser.
(In [Sw] this theorem is stated for arbitrary coefficients, but this is
incorrect. The author is indebted to D.~J.~Green for pointing this out
to him and
providing a counterexample.) The only cases of this theorem that we shall need
are implied by the following lemma.
\mark{The Integral Cohomology of the Held Group}
\proclaim \Lemma 5\cddot1. The corestriction from $\cohz i {C_p}$ to
$\cohz i {C_p\oplus C_p}$ is trivial for $i >0$.
\par
\proof $\cohz i {C_p\oplus C_p}$ has exponent $p$ for $i>0$, and the
restriction from $C_p\oplus C_p$ to $C_p$ is onto. \qed\par
\proclaim \Corollary 5\cddot2. If a Sylow $p$-subgroup $G_p$ of $G$ is elementary
abelian of rank two then $\cohz * G \cong {\cohz * {G_p}}^{N(G_p)}$. If $G_p$
is isomorphic to $P_2$ then the only non-trivial conditions for an element of
${\cohz * {G_p}}^{N(G_p)}$ to be stable arise from intersections $G_p\cap
G_p^g$ of order $p^2$. \qed\par
We require the following lemma which is due to Dickson [Di].
\proclaim \Lemma 5\cddot3. Let $GL_2(p)$ act in the usual way on a vector space
spanned by $x$ and $x'$. Then $\Bbb F_p[x,x']^{SL_2(p)}$ is a polynomial
algebra on generators $a$ and $b$ where
$$a=x^px'-x'^px\qquad b=x^{p(p-1)}+x^{(p-1)^2}x'^{p-1}+\ldots+x'^{p(p-1)},$$
and $\Bbb F_p[x,x']^{GL_2(p)}$ is a polynomial algebra on $a^{p-1}$ and $b$.
\par
\proof Consider the polynomials as functions from $\Bbb F_p^2$ to $\Bbb F_p$.
The action of $SL_2(p)$ on $\Bbb F_p^2$ is transitive on non-zero elements, so
if $P$ is fixed it must represent either the trivial function, in which case it
must divide by each of $x'$, $x+ix'$, and hence by $a$, or we may assume that
it represents the function sending every non-zero point of $\Bbb F_p^2$ to 1.
We shall describe such a polynomial as being identically one. We claim now
that if $P$ is identically one then the degree of $P$ is divisible by $p(p-1)$.
First, setting $x=0$ and letting $x'$ vary we see that $P$ is of the form
$xQ+x'^n$, and that $p-1$ divides $n$. Now
$$\sum \lambda_ix^ix'^{n-i}=P(x,x')=P(x,x+x')=\sum\lambda_ix^i(x+x')^{n-i},$$
where $\lambda_0=1$. Balancing coefficients of $xx'^{n-1}$ gives
$\lambda_1=\lambda_1+n\lambda_0$, so $p$ must divide $n$.
$$b=(\prod_{i=0}^{p-1}(x+ix'))^{p-1}+x'^{p(p-1)}\eqno(*)$$
Equation ($*$) makes apparent the fact that $b$ is identically one, because the
first term of the right hand side
is one if $x'=0$ and zero otherwise. It is also apparent from
($*$) that $b$ is fixed by $x'\mapsto x'$, $x\mapsto x+x'$, and there is a
similar expression with $x$ and $x'$ interchanged that illustrates the fact
that $b$ is fixed by $x'\mapsto x'+x$, $x\mapsto x$. These two elements
generate $SL_2(p)$, so $b$ is invariant as claimed. Now given $Q$ an arbitrary
homogeneous polynomial fixed by $SL_2(p)$, either $Q$ is identically zero and
we may divide $Q$ by $a$, or $Q$ is identically one and of degree divisible by
$p(p-1)$, in which case for some $n$, $Q-b^n$ is homogeneous of
degree the same as $Q$ and identically zero and fixed by $SL_2(p)$. \par
The group $GL_2(p)$ is generated by $SL_2(p)$ and the map that fixes $x$ and
sends $x'$ to $\lambda x'$, for $\lambda$ a generator of $\Bbb F_p^\times$.
Hence $\Bbb F_p[x,x']^{GL_2(p)}$ is the fixed point subring of $\Bbb F_p[a,b]$
under the map which fixes $b$ and sends $a$ to $\lambda a$. \qed\par
\proclaim \Corollary 5\cddot4. If $\cohz * {C_p\oplus C_p}=\Bbb Z[x,x']\otimes
\Lambda[w]$ and $GL_2(p)$ acts in the usual way on $C_p\oplus C_p$, then
$${\cohz * {C_p\oplus C_p}}^{SL_2(p)}=\Bbb Z[a,b]\otimes \Lambda[w]
\qquad\hbox{ and }$$ $$
{\cohz * {C_p\oplus C_p}}^{GL_2(p)}=\Bbb Z[a^{p-1},b]\otimes\Lambda[a^{p-2}w].$$
\par
\proof If $\lambda\in GL_2(p)$, the action of $\lambda$ sends $w$ to
${\rm det}(\lambda)w$.\qed\par
Let $G$ be an extension of $(C_p)^2$ by $SL_2(p)$ or any subgroup of $GL_2(p)$
containing $SL_2(p)$. The metacyclic group $P_1$ has a $(C_p)^2$ subgroup, but
this contains a characteristic subgroup of order $p$ (the subgroup of
$p$-th powers). It follows that the Sylow subgroups of $G$ are isomorphic to
$P_2$, and so $G$ must be the split extension. The groups $(C_p)^2\split
SL_2(p)$ and $(C_p)^2\split GL_2(p)$ seem to arise frequently
in cohomology calculations, so we
determine the $p$-part of their cohomology rings below.
\proclaim \Theorem 5\cddot5. Let $H\cong C_p\oplus C_p$, and let $Q$ be a
subgroup of $GL_2(p)$ containing $SL_2(p)$ and acting in the standard way
on $H$. Also let $C$ be a Sylow subgroup of $Q$, so $C\cong C_p$, and
$H\split C\cong P_2$ is a Sylow subgroup of $H\split Q$. Then the
image of $\cohz * {H\split Q}$ in $\cohz * {H\split C}$ is the subring
$$\{\xi\in {\cohz * {H\split C}}^{N(H\split C)} \mid
\res^{H\split C}_H(\xi)\in{\cohz * H}^Q\}.$$
Note also that $N_{H\split Q}(H\split C)=H\split N_Q(C)$,
that $N_Q(C)$ is isomorphic to
$C_{p-1}\oplus{\rm det}(Q)$, and $N_Q(C)/C_Q(C)$ is isomorphic to
$C_{{(p-1)}/2}$. \par
\proof The remarks at the end of the statement are standard group theory, and
were included to indicate how the main result can be applied. By considering
the restriction maps from $H\split Q$ to $H\split C$ and $H$
it may be seen that the image
is no bigger than claimed. Conversely, corollary~5\cddot4 tells us that given
an element $\xi$ of $\cohz * {H\split C}$ fixed
by the action of the normaliser we
need only check that for each $g\notin N(H\split C)$,
$$\cor^{H\split C}_H\res^{g(H\split C)g^{-1}}_Hc^*_g(\xi)=|C|\xi.$$
Since $H$ is normal
$$\cor^{H\split C}_H\res^{g(H\split C)g^{-1}}_Hc^*_g(\xi)=
\cor^{H\split C}_Hc_g^*\res^{H\split C}_H(\xi),$$
and this will be equal to $|C|\xi$ if $\res(\xi)$ is fixed by $Q$. \qed\par
\proclaim \Corollary 5\cddot6. If we take $P_2 \leq (C_p)^2\split SL_2(p)$
such that
the subgroup $\langle B,C\rangle$ is the normal subgroup, then the image of
$\cohz * {(C_p)^2\split SL_2(p)}$ in $\cohz * {P_2}$ has module generators
$$\zeta^i\alpha^j\beta^k\mu^\epsilon\qquad\hbox{for $j>0$, $\epsilon=0,1$ and
$i+2j\equiv k \quad(p-1)$,}$$
$$\zeta^{(p-1)i-j}\chi_j\quad\hbox{for $j<p-1$,}\qquad
\zeta^{i(p-1)}(\chi_{p-1}+\beta^{p-1}),$$
$$\zeta^i\alpha^j\nu \quad\hbox{where $i+2j\equiv -3 \quad(p-1)$,}\quad
\hbox{and}\quad (\zeta\beta)^i(\zeta^{p-1}+\beta^{p(p-1)})^j\mu^\epsilon.$$
\par
\proof First we determine the fixed points of $\coh {P_2}$ under the action of
its normaliser, which is $P_2\split C_{p-1}$. This group contains a chain of
normal subgroups $\langle C\rangle$, $\langle B,C\rangle$, $P_2$. Let $\lambda$
be a generator for $\Bbb F_p^\times$. The action by conjugation of a generator
of $C_{p-1}$ on $\langle B,C\rangle$ must be faithful, preserve $\langle
C\rangle$, and have determinant one, so without loss of generality may be taken
to be $C\mapsto C^\lambda$, $B\mapsto B^{\lambda^{-1}}$. Also this element acts
on the quotient group $P_2/\langle C\rangle=\langle \widebar A, \widebar
B\rangle$. The action on $\langle C\rangle$ is the determinant of the action on
$\langle \widebar A, \widebar B\rangle$, so the element must send $\widebar A$
to $\widebar A^{\lambda^2}$. The induced action on $\coh {P_2}$ is
$$\hbox{$\eqalign{\alpha &\mapsto \lambda^2\alpha\cr
\beta &\mapsto \lambda^{-1}\beta}\qquad
\eqalign{\mu &\mapsto \mu\cr
\nu &\mapsto \lambda^3\nu}\qquad
\eqalign{\zeta &\mapsto \lambda^p\zeta\cr
\chi_i &\mapsto \lambda^i\chi_i}$}$$
Since each monomial is an eigenvector, it is easy to identify the fixed
subspace, which has the following module generators.
$$\zeta^i\alpha^j\beta^k\mu^\epsilon\qquad\hbox{for $j>0$, $\epsilon=0,1$ and
$i+2j\equiv k \quad(p-1)$,}$$
$$\zeta^{(p-1)i-j}\chi_j,\qquad\zeta^i\alpha^j\nu\quad\hbox{for $i+2j+3\equiv
0\quad(p-1)$.}$$
The kernel of the restriction to $\langle B,C\rangle$ is the ideal generated by
$\alpha,\allowbreak\chi_2,\ldots,\chi_{p-2},\allowbreak\chi_{p-1}+\beta^{p-1},
\allowbreak p\zeta$, and $\nu$, and this ideal together with the subring
generated by $\zeta,\allowbreak \beta$, and $\mu$ spans $\coh {P_2}$ as a
module. Identifying elements of this subring with their images under
restriction to $\langle B, C\rangle$, it follows from corollary~5\cddot4 that
the fixed points of $\coh {\langle B, C\rangle}$ under the action of
$N(\langle B, C\rangle)$ are generated by $\zeta\beta$,
$\zeta^{p-1}+\beta^{p(p-1)}$, and $\mu$. The result follows. \qed\par
\beginsubsection The Held Group. Now we shall consider the cohomology of the
Held group~$\held$, a sporadic simple group of order
$2^{10}.3^3.5^2.7^3.17$, [He].
The Sylow 3- and 7-subgroups are isomorphic to $P_2$, so the results of
section~2 will be relevant. By the end of this section we shall have determined
$\co * \held {\Bbb Z[{1\over 2}]}$ by considering separately each of the
other primes that divide~$|\held|$. We list our results for the primes in order
of increasing difficulty. In the sequel we shall assume any
results concerning $\held$ and other groups that are stated explicitly in the
Atlas [Co], but we shall prove everything else that we require in a series of
lemmata, one for each $p$. Much of the content of these is in Held's original
paper [He].
\beginsubsection The prime 17. If $G$ has a cyclic Sylow subgroup $C$, it
follows from Swan's theorem (or directly from the stable element criterion if
$C$ has order $p$) that ${\cohz * G}_p$ is a polynomial algebra on a generator
of degree $2|N_G(C):C_G(C)|$. More recently Thomas has shown that in this case
the generator may be taken to be a Chern class of a representation (see [Th3],
together with the observation that $c.(-\rho)c.(\rho)=1$ so the Chern subring
is generated by classes of actual representations). In $\held$ there are two
classes of elements of order 17, so $|N_\held(C_{17}):C_\held(C_{17})|$
has order 8, and ${\cohz * \held}_{17}$ is a polynomial algebra on a generator
of degree 16. \par
\beginsubsection The prime 5.
\proclaim \Lemma 5\cddot7. The group
$\held$ has a unique conjugacy class of elements of order 5, and its Sylow
5-subgroup is $5A^2$, that is it is isomorphic to $C_5\oplus C_5$. The
centraliser of $5A^2$ is $5A^2$, and $N(5A^2)\cong (C_5)^2\split C_4A_4$.
The group $\aut((C_5)^2)$ contains a unique class of $C_4A_4$ subgroups.
\par
\proof The assertions concerning $\held$ are stated explicitly in the Atlas,
except for the claim that $5A^2$ is its own centraliser. Note that
$C(5A)$ has order 300, and of course $C(5A^2)\subset C(5A)$. $\held$ has only
two classes of maximal subgroups of order divisible by 25, $N(5A^2)$, and
another isomorphic to $S_4(4)\split C_2$. The only maximal subgroups of
$S_4(4)\split C_2$ with order divisible by 25 are isomorphic to
$(A_5\times A_5)\split(C_2\times C_2)$. In this group the subgroup
$C_5\times A_5$ has order 300 and centralises an element of order 5, so must be
$C(5A^2)$, but $C_5\times C_5$ is its own centraliser in $C_5\times A_5$.
\par
It remains to check the assertions concerning $GL_2(5)$. A Sylow 2-subgroup of
$GL_2(5)$, for example the group of diagonal and `anti-diagonal' matrices,
is non-abelian of order 32. It follows that any abelian subgroup of $GL_2(5)$
of order 16 must contain the centre. $PGL_2(5)\cong S_5$, so the result follows
since $S_5$ has a unique class of $A_4$ subgroups. \qed\par
\proclaim \Theorem 5\cddot8. ${\cohz * \held}_5$ is generated by elements
$\alpha$, $\beta$, $\gamma$ and $\chi$,
where $\alpha$ has degree~16, $\beta$ and $\gamma$ have degree~24
and $\chi$ has degree~15, subject only to the relations that all generators
have order 5 and $\gamma^2=3(\beta^2+\gamma^3)$.
The cohomology operation $\delta_p\pone\pi_*$ sends $\chi$ to
a multiple of $\gamma$.
\par
\proof It follows from corollary~5\cddot2 that ${\cohz * \held}_5$ is
isomorphic to the fixed points of $\cohz * {5A^2}$ under the action of
$N(5A^2)$. By lemma~5\cddot7 we may
take as generators for $N(5A^2)/5A^2$ considered as a
subgroup of $GL_2(p)$ the matrices
$$M_1=\pmatrix{2&0\cr 0&3}\quad
M_2=\pmatrix{2&0\cr 0&2}\quad
M_3=\pmatrix{0&1\cr -1&0}\quad
M_4=\pmatrix{1&1\cr 2&-2}.$$
Note that these define a composition series for this group, in the sense that
the subgroup generated by the first $i$ of them is normal in the subgroup
generated by the first $i+1$ of them. Let $\cohz * {5A^2}$ be generated by
elements $\delta$ and $\delta'$ of degree 2 and $\epsilon$ of degree 3, such
that the above matrices describe the action of $C_4A_4$ on the basis
$\delta,\delta'$ for $H^2$. Thus $M_1$ sends $\delta$ to $2\delta$, $\delta'$ to
$3\delta'$, and fixes $\epsilon $. The fixed point subring under the action of
$M_1$ is easily seen to be generated by $\delta^4,\allowbreak
\delta'^4,\allowbreak \delta\delta'$, and $\epsilon$. Since $M_2$ normalises
$M_1$, $M_2$ acts on this fixed point subring, and in fact it fixes $\delta^4$
and $\delta'^4$ and sends $\delta\delta'$ to $-\delta\delta'$ and $\epsilon$ to
$-\epsilon$. The fixed point subring is generated by $\delta^4$, $\delta'^4$,
$\delta^2\delta'^2$ and $\delta\delta'\epsilon$. The action of $M_3$ fixes
$\delta^2\delta'^2$, exchanges $\delta^4$ and $\delta'^4$ and sends
$\delta\delta'\epsilon$ to minus itself, so has fixed point subring generated
by $\delta^4+\delta'^4$, $\delta^2\delta'^2$, and $(\delta^5\delta'-
\delta'^5\delta)\epsilon$. This is a polynomial algebra on two generators of
degree~8 tensored with an exterior algebra on one generator of degree~15.
\par
The action of $M_4$ clearly fixes the exterior generator, but has no fixed
point in degree~8. To complete the proof it suffices to show that whenever
$C_3$ acts non-trivially on $\Bbb F_5[x,x']$, then the fixed point subring is
generated by elements $a$, $b$, $c$, of degrees two, three and three
respectively (where $x$ is of degree one), subject only to $c^2=3(b^2+a^3)$.
For this, extend the action to one on $\Bbb F_{25}[x,x']$, where it may be
diagonalised. Note that if $y$ is an eigenvector in degree one with eigenvalue
$\omega$, then $\bar y$ (defined by conjugating the coefficients of $y$ as an
expression in $x$ and $x'$) is an eigenvector of value $\bar\omega$.  The fixed
point subring over $\Bbb F_{25}$ is generated by $y\bar y$, $y^3$ and $\bar
y^3$. An element of this ring is in $\Bbb F_5[x,x']$ iff it is invariant
under
conjugation. Elements invariant under conjugation are the sum of a symmetric
polynomial over $\Bbb F_5$ in $y$ and $\bar y$, and $\sqrt 3$ times an
antisymmetric polynomial over $\Bbb F_5$ in $y$ and $\bar y$. The symmetric
polynomials in $y$ and $\bar y$ expressible in terms of $y\bar y$, $y^3$ and
$\bar y^3$ form a polynomial ring on $y\bar y$ and $y^3+\bar y^3$, while the
antisymmetric polynomials form a free module for the symmetric ones with basis
$y^3-\bar y^3$. Our three generators are $a=y\bar y$, $b=y^3+\bar y^3$ and
$c=\sqrt 3(y^3-\bar y^3)$. \qed \par
   The author wishes to thank Dr.~D.~J.~Benson who found an error in an earlier
statement of theorem~5\cddot8.
\par
\beginsubsection The prime 3.
\proclaim Lemma~5\cddot9. Let $P$ be a Sylow 3-subgroup of $\held$. Then
$P\cong P_2$, and $\held$ contains two (self-inverse) conjugacy classes of
elements of order 3, $3A$ and $3B$. $N(P)$ is isomorphic to $P\split D_8$ and
$D_8$ acts faithfully. The centre of $P$ is $3A$. There are two classes of
$C_3\oplus C_3$ subgroups of $\held$, $3A^2$ which is contained in a unique
Sylow 3-subgroup, and $3A3B$, which has
$$N(3A3B)\cong C_3\nonsplit (S_3\times S_4)\qquad
C(3A3B)\cong C_3\nonsplit (C_3\times V_4).$$
There is a unique class of $D_8$ subgroups of $GL_2(3)$.
\par
\proof The order of $C(3B)$ is not divisible by 27, so $Z(P)$ must be $3A$. One of
the maximal subgroups of $\held$ is $N(3A)\cong C_3\nonsplit S_7$, and we shall
work largely within this group. For example, $N(P)\leq N(3A)$, and we see that
$$N(P)=C_3\nonsplit N_{S_7}(C_3\times C_3)\cong P\split D_8.$$
If we write $3X$ for the conjugacy class in $\held$ represented by the 3-cycles
in $S_7$, and $3Y$ for that represented by products of two disjoint 3-cycles,
then
$$ C(3A,3X)\cong C_3.(C_3\times V_4)\qquad C(3A,3Y)\cong C_3\times C_3$$
We deduce that there are two classes of $C_3\oplus C_3$ subgroup in $\held$,
one of which must be $3A^2$, and the other $3A3B$ containing a unique $3A$
subgroup. Now
$$|N(3A,3B):C(3A,3B)|\leq 12\quad\hbox{and}\quad |N(3A^2):C(3A^2)|\leq 48,$$
but $N(3A,3Y)\cap N(3A)\cong C_3\nonsplit(S_3\times S_3).C_2$, which already has
order 24 times that of $C(3A,3Y)$. Therefore $3Y=3A$, and $N(3A^2)$ is either
$C_3\nonsplit(S_3\times S_3).C_2$ or a group containing this with index two. In
either case $N(3A^2)$ contains a unique Sylow subgroup. Also $3X=3B$, and
$N(3A,3B)$ is as claimed. \par
$\aut(\Bbb F_9,+)$ contains the subgroup generated by $\Bbb F_9^\times$ and the
Fr\"obenius automorphism, which is semidihedral of order 16, and contains a
unique $D_8$ subgroup. Since the Sylow 2-subgroups of $GL_2(3)$ have order 16
it follows that all $D_8$ subgroups of $GL_2(3)$ are conjugate.
\par\ \hfill\qed\par
\proclaim \Theorem 5\cddot10. ${\cohz * \held}_3$ is isomorphic to the subring
of $\cohz * {P_2}$ generated by $\chi_2,\allowbreak
\alpha^2+\beta^2,\allowbreak \alpha^2\beta^2,\allowbreak
\zeta(\alpha\nu-\beta\mu)$, and $\zeta^2$.
\par
\proof Identifying $\aut(P_2)/P_2$ with $GL_2(p)$, it follows from
lemma~5\cddot9 that we may take $N(P_2)$ as being generated by
$$M_1=\pmatrix{1&0\cr 0&-1}\qquad M_2=\pmatrix{-1&0\cr 0&1}\qquad
M_3=\pmatrix{0&1\cr -1&0}.$$
(Of course we only need $M_1$ and $M_3$). The action of $M_1$ fixes
$\beta,\allowbreak \nu$, and $\chi_2$, and multiplies the other generators by
$-1$. The action of $M_2$ fixes $\alpha,\allowbreak\mu$, and $\chi_2$ and
multiplies the other generators by $-1$. It follows that the fixed point
subring under $M_1$ and $M_2$ is spanned by the elements
$$\zeta^{2i},\qquad\zeta^{2i}\chi_2,\qquad
\zeta^{2i+\epsilon}\alpha^{2j+\epsilon}\nu,$$
$$\zeta^i\alpha^j\beta^k\nu^\epsilon\qquad\hbox{where $\epsilon+j\equiv i\equiv
k\quad (2)$, and $j=0$ or $k+\epsilon\leq 3$.}$$
The action of $M_3$ on the generators is given by
$$\hbox{$\eqalign{\alpha &\mapsto -\beta\cr
\beta &\mapsto \alpha}\qquad
\eqalign{\mu &\mapsto \nu\cr
\nu &\mapsto -\mu}\qquad
\eqalign{\zeta &\mapsto -\zeta\cr
\chi_2 &\mapsto \chi_2}$}$$
and a simple calculation shows that the fixed points under the action of $D_8$
are spanned by
$$\zeta^{2i}\chi_2,\qquad\zeta^{2i+1}(\alpha^{2j+1}\nu-\beta^{2j+1}\mu),$$
$$\zeta^{2i}(\alpha^{2j}+\beta^{2j}),\qquad\zeta^{2i}\alpha^{2j}\beta^2,$$
and further that this is the subring generated by the elements of the
statement. To complete the proof we must show that all these elements are
stable. Without loss of generality we may take $\langle A,C\rangle$ and
$\langle B,C\rangle$ to be the $3A3B$ subgroups of $P_2$, and by
corollary~5\cddot2 it suffices to check that for $H=\langle B,C\rangle$ and
$\xi$ any of the generators of the statement then
$$\cor^{P_2}_Hc_g^*\res^{P_2}_H(\xi)=3\xi=\cor^{P_2}_H\res^{P_2}_H(\xi)
\qquad\hbox{for all $g\in N_\held(H)$.}$$
If we take $\cohz * H \cong \Bbb Z[\beta',\gamma]\otimes\Lambda[\delta]$, then
the action of an element of $N_\held(H)$ sends $\beta'$ to $\lambda\beta'$,
$\gamma$ to $\lambda'\gamma+\lambda''\beta'$ and $\delta$ to
$\lambda\lambda'\delta$, where $\lambda,\lambda'\in \Bbb F_3^\times$
and $\lambda'' \in \Bbb F_3$. The restriction from $P_2$ to $H$ is given by
$$\hbox{$\eqalign{\alpha &\mapsto 0\cr
\beta &\mapsto \beta'}\qquad
\eqalign{\mu &\mapsto \delta\cr
\nu &\mapsto 0}\qquad
\eqalign{\zeta &\mapsto \gamma^3-\beta^2\gamma\cr
\chi_2 &\mapsto -\beta^2}$.}$$
It can be checked that the images of the generators of the statement are fixed
by $N_\held(H)$, so the generators satisfy the stability condition.\qed\par
\beginsubsection The prime 7.
\proclaim \Lemma 5\cddot11. Let $P$ be a Sylow 7-subgroup of $\held$. then
$P\cong P_2$, and $\held $ contains five conjugacy classes of element of order
seven, $7A,7B$ which are inverse to each other, $7C$ which is self-inverse, and
$7D,7E$ which are inverse to each other. The central elements are $7C$, and in
fact $N(7C)=N(P)\cong P\split (S_3\times C_3)$, with $S_3\times C_3$ acting
faithfully on $P$. There is a unique class of such subgroups of $GL_2(7)$.
There are three conjugacy classes of $C_7\times C_7$ subgroups of $\held$, with
normalisers
$$N(7C,7AB)\cong P\split(C_3\times C_3)\quad
N(7C^2)\cong (C_7)^2\split SL_2(7)\quad
N(7C,7DE)\cong P:(C_2\times C_3).$$
$7C^2$ is its own centraliser. Some other subgroups are
$$(C_7)^2\split C_6\cong N(7DE)\leq N(P)\qquad
N(7AB)\cong C_7 \split C_3\times L_2(7),$$ $$
(C_7)^2\split (C_3\times C_3)\cong N(7C,7AB)\cap N(7AB)\leq N(P).$$
\par
\proof The normaliser of a Sylow 3-subgroup of $GL_2(7)$ is $C_6\wr C_2$, which
contains a unique conjugacy class of $C_3\wr C_2\cong S_3\times C_3$ subgroup.
The Atlas tells us that $N(7C^2)$ is as claimed above, but $L_2(7)$ is not a
subgroup of $GL_2(7)$, so the only possibility is that $7C^2$ is
self-centralising. The Atlas tells us that $N(7C)$ is
$P\split (S_3\times C_3)$ as claimed, and since no element of $S_3\times C_3$
can centralise a $7C^2$ subgroup of $P$ we see that $S_3\times C_3$ acts
faithfully on $P$. We already know that such an action is essentially unique,
and so we fix generators for this group considered as acting on $C_7^{\bar A}
\oplus C_7^{\bar B} =P_2/\langle C\rangle$ as below.
$$\pmatrix{\lambda & 0\cr 0&\lambda'}\quad\hbox{where $\lambda^3=\lambda'^3
=1$}\qquad\pmatrix{0& 1\cr 1& 0}.$$
Direct computation shows that there are four conjugacy classes of $C_7$
subgroup in $P$ under the action of $N(P)$, which are tabulated in figure
5\cddot1 together with their centralisers in $N(P)$. For each of them except
$\langle C\rangle$ the centraliser has index three in the normaliser.
\par
\midinsert
$$\vcenter{\halign{\hfil#\hfil&\qquad\hfil#\hfil\cr
$C_7$ subgroup&Centraliser in $N(P)$\cr
$\langle C\rangle$&$P.C_3$\cr
$\langle  AC^i \rangle\sim \langle BC^i  \rangle $&$(C_7)^2\split C_3$\cr
$\langle AB^{2^i}C^j  \rangle$&$(C_7)^2\split C_2$\cr
$\langle AB^{-2^i}C^j  \rangle$&$(C_7)^2$\cr}}$$
\figure 5/1/Table of $C_7$ subgroups of $N(P)$/
\endinsert
$\held$ has only two classes of maximal subgroup whose order divides by $7^3$,
$N(7C^2)$ and $N(7C)$. We see that $\langle A,C  \rangle$ cannot be $7C^2$
because its normaliser has order divisible by~9. Also its normaliser must be
contained in $N(7C)$, so is $P.(C_3\times C_3)$. $\langle  A \rangle$ cannot
be $7C$ by the above, but neither can it be $7DE$ since the order of $C(7DE)$
is coprime to 3, hence $\langle  A \rangle$ must be $7AB$. There are $7C^2$
subgroups, but $\langle C,AB^{2^i}  \rangle $ cannot be $7C^2$ because it has
elements centralised by an involution. It follows that
$\langle AB^{2^i}  \rangle $ must be $7DE$ and $\langle AB^{-2^i}  \rangle $
is $7C$. Comparing orders tells us that $N(7DE)$ is contained in $N(P)$ and is
as claimed. The normaliser of $7AB$ is a maximal subgroup listed in the Atlas.
The intersection of $N(7C,7AB)$ and $N(7AB)$ has index seven in $N(7C,7AB)$,
because this group permutes its $7AB$ subgroups transitively, so is as claimed.
\qed \par
\proclaim \Theorem 5\cddot12. With notation as in lemma~5\cddot11, ${\cohz *
{N(P)}}_7$ is mapped via restriction to the subring of $\cohz * P$ generated by
$$\alpha^3+\beta^3,\quad\alpha^3\beta^3,\quad\chi_6,\quad
\alpha^5\beta\mu-\alpha^2\beta^4\mu,\quad\zeta\alpha\mu\quad\zeta\chi_5,\quad
\zeta(\alpha^5\beta^2-\alpha^2\beta^5),\quad\zeta^2\alpha\beta,$$
$$\zeta^2(\alpha^2\nu-\beta^2\mu),\quad\zeta^2\chi_4,\quad
\zeta^3\chi_3,\quad\zeta^3(\alpha^3-\beta^3),\quad\zeta^4\chi_2,\quad
\zeta^5(\alpha^2\nu+\beta^2\mu),\quad\zeta^6.
$$
\par
\proof As in theorem 5\cddot10 we take a basis for ${\cohz * P}\otimes\Bbb F_7$
consisting of the monomials
$$\zeta^i,\quad\zeta^i\chi_j,\quad\zeta^i\alpha^j\nu,\quad
\zeta^i\alpha^j\beta^k\mu^\epsilon\quad\hbox{where $j=0$ or $k+\epsilon\leq
6$}.$$
All monomials are eigenvectors for the action of the $C_3\times C_3$
subgroup of $N(P)/P$, and it may be checked that the fixed monomials are
$$\zeta^{3i},\quad\zeta^{3i-j}\chi_j,\quad \zeta^{3i+2}\alpha^{3j+2}\nu$$
$$\zeta^i\alpha^j\beta^k\mu^\epsilon\quad\hbox{where $i+j+\epsilon\equiv 0$ and
$i+k\equiv\epsilon$ modulo 3.}$$
The action of the involution represented by the matrix in the statement of
lemma~5\cddot11 on the original generators is
$$\hbox{$\eqalign{\alpha &\mapsto \beta\cr
\beta &\mapsto \alpha}\qquad
\eqalign{\mu &\mapsto -\nu\cr
\nu &\mapsto -\mu}\qquad
\eqalign{\zeta &\mapsto -\zeta\cr
\chi_i &\mapsto (-1)^i\chi_i}$}$$
and the fixed points under the whole of $N(P)/P$ are seen to be spanned by
$$\zeta^{6i},\quad\zeta^{6i-j}\chi_j,\quad
\zeta^{6i+2}(\alpha^{6j+2}\nu-\beta^{6j+2}\mu),\quad
\zeta^{6i}(\alpha^{3j}+\beta^{3j}),$$
$$\zeta^{6i-1}(\alpha^{6j-1}\nu+\beta^{6j-1}\mu),\quad
\zeta^i\alpha^j\beta^k\mu^\epsilon\quad\hbox{where $i\equiv2j-\epsilon\equiv
j+k\quad (6)$,}$$
$$\zeta^i(\alpha^j\beta^k\mu^\epsilon+\alpha^{j-3}\beta^{k+3}\mu^\epsilon)\quad
\hbox{where $i\equiv2j-\epsilon\equiv j+k+3\quad (6)$,}$$
$$\zeta^i(\alpha^j\beta^k\mu^\epsilon-\alpha^{j-3}\beta^{k+3}\mu^\epsilon)\quad
\hbox{where $i\equiv2j-\epsilon+3\equiv j+k\quad (6)$.}$$
More calculation shows that the fifteen elements of the statement generate this
subring.\qed\par
In a sense this almost completes the calculation of ${\cohz * \held}_7$, which
may be described as the subring of stable elements of ${\cohz * {N(P)}}_7$.
We know from corollary~5\cddot2 that only conjugates $P^g$ of $P$ such that
$P\cap P^g$ has order 49 can give any extra conditions, and by lemma~5\cddot11
the only such intersections are of type $7C^2$. It follows as in
theorem~5\cddot5 that $\xi\in{\cohz * {N(P)}}_7$ is in the image of $\cohz *
\held$ iff
$$\cor^P_{\langle AB^{-1},C\rangle}\phi^*\res^P_{\langle AB^{-1},C\rangle}(\xi)
=p\xi,\eqno(*)$$
where $\phi$ is some fixed element of $\aut({\langle A,C\rangle})$ of order 7
not fixing $\langle C \rangle$. (It is not obvious, but true, that all such
$\phi$ will give the same condition.)
\par
The description of ${\cohz * \held}_7$ given by condition $(*)$ is not very
explicit, since the condition is difficult to check. For example, this
description is not much use for calculating the Poincar\'e series of
${\cohz * \held}_7\otimes \Bbb F_7$. The following lemma gives a more useful
description of ${\cohz * \held}_7$, which will allow us to completely determine
${\cohz * \held}_7$ in theorem~5\cddot14.
\proclaim \Lemma 5\cddot13. Let $K$ be the subgroup $\langle AB^{-1},C\rangle$
of $P$, let $\res=\res^P_K$, and define $\zeta',\epsilon,\delta$,
elements of $\cohz * K$ by
$$\zeta'=\res(\zeta)\quad\epsilon=\res(\alpha)=-\res(\beta)\quad
\delta=\res(\mu)=-\res(\nu).$$
Now let $S$ be the subring of $\cohz * K $ generated by $\zeta'\epsilon$,
$\zeta'^6+\epsilon^{42}$ and $\delta$. Then ${\cohz * \held}_7$ is isomorphic
to $(\res)^{-1}(S)\cap{\cohz * {N(P)}}_7$.
\par
\proof From corollary 5\cddot2 it follows that an element $\xi$ of $\cohz *
{N(P)}_7$ is in the image of the restriction from $\cohz * \held_7$ iff
$$c_g^*\res^P_{g^{-1}Pg\cap P}(\xi)=\res^P_{P\cap gPg^{-1}}(\xi)$$
for all $g$ such that $P^g\cap P$ has order $p^2$. This can only happen if
$P^g\cap P$ is of type $7C^2$ (see lemma~5\cddot11). Since $N(P)$ acts
transitively on such subgroups we may restrict to the case when $P^g\cap P$ is
the subgroup $K=\langle AB^{-1},C\rangle$ and $g$ normalises $K$.
Lemma~5\cddot11 tells us that $N(7C^2)\equiv (C_7)^2\split SL_2(7)$, and the
ring $S$ described in the statement is (by corollary~5\cddot4) the fixed point
subring of $\cohz * K$ under the action of $SL_2(7)$. \qed\par
\proclaim \Theorem~5\cddot14. With notation as in lemma~5\cddot11, $\cohz *
\held_7$ is mapped via restriction to the subring of $\cohz * P$ generated by
the following elements.
$$\alpha^3+\beta^3,\quad\chi_6-\alpha^3\beta^3,\quad\zeta\alpha\mu,\quad
\zeta\chi_5,\quad\zeta^2\alpha\beta,\quad \zeta^2(\alpha^2\nu-\beta^2\mu),$$
$$\zeta^2\chi_4,\quad\zeta^3\chi_3,\quad\zeta^3(\alpha^3-\beta^3),\quad
\zeta^4\chi_2,\quad \zeta^5(\alpha^2\nu+\beta^2\mu),\quad
\zeta^6-\alpha^{39}\beta^3.$$
\par
\proof The twelve elements of the statement with $\alpha^3\beta^3$,
$(\alpha^5\beta-\alpha^2\beta^4)\mu$, and
$\zeta(\alpha^5\beta^2-\alpha^2\beta^5)$ generate the subring of $\cod * P$
described in theorem~5\cddot12. Also the images of the twelve elements under
restriction from $P$ to the subgroup $K$ lie in the subring $S$ (see
lemma~5\cddot13), so these elements are contained in $\cohz * \held_7$.
It remains to show that no other generators are required.
It may be verified that the image of the subring of this statement under
restriction to $K$ is the subring $S'$ of $S$ generated by
$\zeta'^6+\epsilon^{42}$, $\zeta'^2\epsilon^2$, $\zeta'^3\epsilon^3$,
$\zeta'\epsilon\delta$ and $\zeta'^2\epsilon^2\delta$. The elements
$\alpha^3\beta^3$, $(\alpha^5\beta-\alpha^2\beta^4)\mu$, and
$\zeta(\alpha^5\beta^2-\alpha^2\beta^5)$ restrict to $-\epsilon^6$,
$-2\epsilon^6\delta$, and $2\zeta'\epsilon^7$ respectively. We claim now that
the intersection of $S$ and $\res^P_K({\cohz * {N(P)}})$ is the subring $S'$ of
$S$ generated by the elements of the statement. Equivalently we claim that the
intersection of $S$ and the ring generated by $S'$, $\epsilon^6$,
$\epsilon^6\delta$ and $\zeta'\epsilon^7$ is $S'$. In even degrees $S\setminus
S'$ is the set of elements of $S$ divisible exactly once by $\epsilon$, while
elements of even degree in the ring generated by $S'$, $\epsilon^6$ and
$\zeta'\epsilon^7$ will be divisible by $\epsilon$ either at least twice or not
at all. A similar argument works in odd degrees, since the odd degree part of
$S\setminus S'$ is the set of elements not divisible by $\epsilon$.
\par
Now we have shown that the subring of $\cohz * P$ generated by the elements of
the statement has the same image under restriction to $K$ as $\cohz * \held_7$,
so we may assume that any remaining elements of $\cohz * \held_7$ are in the
kernel of $\res^P_K$. It will suffice to show that the kernel of $\res^P_K$ as
a map from $\cohz * {N(P)}_7$ is contained in the subring generated by the
elements of the statement. We claim that this kernel is the ideal of $\cohz *
{N(P)}_7$ generated by $\alpha^3+\beta^3$, $\zeta^5(\alpha^2\nu+\beta^2\mu)$,
$\zeta^4\chi_2$, $\zeta^3\chi_3$, $\zeta^2\chi_4$, $\zeta\chi_5$ and
$\chi_6-\alpha^3\beta^3$. We easily reduce to the case of elements of the
kernel of the form $\zeta^iP(\alpha,\beta)$ or
$\zeta^i(P(\alpha,\beta)+\lambda\alpha^r\nu)$. Notice that if $n\geq 7$ and
$$\zeta^j(\sum^6_{i=0}\lambda_i\alpha^{n-i}\beta^i+\lambda_7\beta^n) \in
\cohz * {N(P)}_7,\qquad\hbox{then}\quad$$
$$\zeta^{6k+j}(\sum^6_{i=0}
\lambda_i\alpha^{6m+n-i}\beta^i+\lambda_7\beta^{6m+n}) \in \cohz * {N(P)}_7,$$
so the remaining cases reduce to $\zeta^iP(\alpha,\beta)$ and
$\zeta^i(P(\alpha,\beta)\mu+\lambda\alpha^r\nu))$ where the degree of $P$ is at
most 13 and $i\leq 5$. It can be checked by lengthy calculation that all such
elements are in the ideal of $\cohz * {N(P)}_7$ generated by $\alpha^3+\beta^3$
and $\zeta^5(\alpha^2\nu+\beta^2\mu)$. To check that the kernel is contained in
the subring of the statement it suffices to find an expression for each product
of one of the generators of the kernel (as an ideal) and one of
$\alpha^3\beta^3$, $(\alpha^5\beta-\alpha^2\beta^4)\mu$ and
$\zeta(\alpha^5\beta^2-\alpha^2\beta^5)$ in terms of the twelve elements of the
statement. This too is a routine calculation. \qed\par

\def\state#1{\par\smallskip \item{#1.}}
\def\sta#1{#1}
\def\po{\hbox{$\rm p\kern -.06em\circ\kern -.06em$}}
\def\pc{{\rm pc}}
\def\es{{\rm S}}
\def\gl{{\rm GL}}
 
\def\lcm{{\rm\scriptstyle L.C.M.}}
 
\beginsection {6.} Yagita's Invariant. \par
\mark{\ }
The following short section concerns an invariant related to free actions of
groups on products of spheres. In [Ya] Yagita defined an integer invariant
$\po(G)$ for a prime $p$ and a finite group $G$ which must divide $n$ if $G$
acts freely on a product $(\es^{n-1})^k$ with trivial action on
$H^*((\es^{n-1})^k;\Bbb Z)$.
\beginsubsection The Invariant. Let $u$ be a generator for $\cohz 2 {C_p}$. To
each inclusion $i$ of $C_p$ in $G$ Yagita associates an integer $m(i)$, which
is the largest $m$ such that $i^*({\cohz * G})$ is contained in the subring of
$\cohz 2 {C_p}$ generated by $u^m$. Then $\po(G)$ is defined to be twice the
lowest common multiple of all such $m(i)$. \par
We recall that given a CW complex $X$ with a free cellular action of $G$ (such
a space is referred to as a free $G$ complex), there is a Cartan-Leray spectral
sequence with $E_2$ page $H^i({\rm B}G;H^j(X))$ converging to a
filtration of $H^{i+j}(X/G)$. Now let $X$ be a finite free $G$
complex, with cohomology ring isomorphic to an exterior algebra on $k$
generators with trivial $G$ action. An inclusion $i: C_p\rightarrowtail G$
induces a map of Cartan-Leray spectral sequences, where the map on $E_2$ pages
is essentially the restriction, and the map on $E_\infty$ pages is induced by
the projection from $X/C_p$ to $X/G$. By considering this map
of spectral sequences Yagita shows that $m(i)$ must divide $n$.
\par
While exploring the possibility of defining a similar invariant
in terms of the images of $\ch G$, the author discovered two errors in
[Ya], and wrote the following letter to Yagita. At the time of writing
there seemed to be no connection between the new invariant and actions on
spheres, but the author has recently found a connection which is explained
in lemma~6\cddot 1. The preprint referred to in the letter contained an
incorrect proof of corollary~3\cddot 7, and a determination of the order of
$\cohz n {G(a,1)}$ which disagrees with Yagita's theorem~2\cddot4. We
explain the method of calculation after the letter.
\par
\noindent\sl
Dear Prof. Yagita,
\par
I am writing to you to report some work that I have done which was inspired by
your paper `On the dimension of spheres whose product admits a free action by a
non-abelian group', Quart. J. Math. Oxford 36 (1985) 117--127. I have also
found two mistakes in this paper. I find that theorem~2\cddot4 is incorrect
for odd primes $p$ for the groups $G_2(a,1)$, but that your other results
that rely on this theorem are correct.
Secondly, your proof of lemma~1\cddot7
is incorrect; in this case I have been unable to find a proof or a
counterexample. Before explaining my criticisms, I shall define my adaptation,
$\pc (G)$, of your invariant. My invariant is easier to calculate than yours,
and gives an upper bound for yours, but does not seem to bear any direct
relationship to actions on spheres. I shall consider only odd primes $p$.
\par
You define, for each $i:\Bbb Z_{/p}\rightarrowtail G$, $m(i)$ to be the
greatest integer $m$ such that $i^*(R)\subset \Bbb Z [u^m]/(pu^m)$, where $R$
is $H^*({\rm B}G;\Bbb Z)$. I define $n(i)$ similarly,
except that I take $R$ to be the
subring of $H^*({\rm B}G;\Bbb Z)$ generated by Chern classes
of representations of $G$. I then define
$$\pc (G)=2\lcm\{n(i) \mid i:\Bbb Z_{/p}\rightarrowtail G \}.$$
(Incidentally, your stated definition of $\po (G)$ is $\lcm \{2m(i)\}$, which
gives $\po (G)=1$ if $p$ does not divide the order of $G$, and gives easy
counterexamples to lemma~1\cddot7. I shall use the definition $2\lcm
\{m(i)\}$.)
\par
\mark{Yagita's Invariant}
Various properties of $\pc (G)$, such as \sta{1} to \sta{4} below, may be proven
using the proofs you use for $\po (G)$:
\state{1} If $G$ is abelian, then $\pc (G)=2$. \par
\state{2} If $N\rightarrowtail G\twoheadrightarrow Q$, then
$\pc (G)\mid \lcm \{\pc(Q),\pc_N(G)\}$. \par
\state{3} If $H < G$, then $\pc(H)\mid\pc_H(G)$. \par
\state{4} If $H < G$, then $\pc_H(G)\mid\pc(G)$. \par
\state{5} $\po(G)\mid\pc(G)$, so my invariant gives an upper bound for yours. A
lower bound would of course be more useful!
It follows from \sta{1} and \sta{2} that
\state{6} $\pc(G) = \pc_{G'}(G)$, where $G'$ is the derived subgroup of $G$.
Similarly, it can be shown that
\state{6}$'$ $\po(G) = \po_{G'}(G)$.\par
If $G$ is a minimal non-abelian $p$-group or an extra-special $p$-group, then
$G'$ is isomorphic to $\Bbb Z_{/p}$, so \sta{6}$'$ could be used to simplify your
calculations for these groups. \par
Using \sta{6} it may be shown that
\state{7} If $G$ is a minimal non-abelian $p$-group, then $\pc(G)=2p$, and
hence
\state{8} If $G$ is a $p$-group, then $\pc(G)= 2$ if and only if $G$ is
abelian. \par
The regular representation of $G$ restricts to a $\Bbb Z_{/p}$ subgroup as
$|G:\Bbb Z_{/p}|$ copies of the regular representation of $\Bbb Z_{/p}$, so
its total Chern class restricts to $\Bbb Z_{/p}$ as
$$(\prod_{i=0}^{p-1}(1-iu))^{|G:\Bbb Z_{/p}|}=(1-u^{p-1})^{|G:\Bbb Z_{/p}|}, $$
and so
\state{9} $\pc(G)\mid 2(p-1)p^{n-1}$, where $p^n$ is the $p$-part of $|G|$.
\par
Using also representations of $G$ induced up from irreducible representations
of $\Bbb Z_{/p}$ I can show that $\pc(G)\mid 2|G|$, but the calculation is
longer. \par
 
I first noticed the flaw in lemma~1\cddot 7 while trying to prove the analogous
result for $\pc(G)$. Let $N\triangleleft G$, and let $x\in H^*({\rm B}N)$. Then
letting $c_t^*$ be the map $\coh N\mapright{}\coh N$ induced by
$n\mapsto n^t$, we have:
$$\res^G_N {\cal N}^G_N(x)=\prod_{t\in G/N} c_t^*(x)$$
(compare this with the formula for $\res^G_N\cor^G_N$), so the equation in your
proof of lemma~1\cddot7 is incorrect in general. I have been unable
either to prove
lemma~1\cddot7 or to find a counterexample. However, the analogous statement
for $\pc(G)$ is false:
\par
Let the cyclic group of order $p^n-1$ act on $(\Bbb Z_{/p})^n$ so as to
permute the non-zero elements transitively (this can be done -- it is the action
by multiplication of the multiplicative group of the field of order $p^n$ on
its additive group). Now let $N=(\Bbb Z_{/p})^n$, $Q=\Bbb Z_{/p^n-1}$, and $G$
the (unique) extension of $N$ by $Q$ corresponding to the above action.
$G$ has $p^n-1$ 1-dimensional representations, which restrict trivially to any
subgroup of order $p$, and one $(p^n-1)$-dimensional representation, whose total
Chern class restricts to any $\Bbb Z_{/p}$ subgroup as $(1-u^{p-1})^{p^{n-1}}
=1-u^{(p-1)p^{n-1}}$. Hence $\pc(G) = 2(p-1)p^{n-1}$, which does not divide
$(p^n-1)\pc(N)$.
\par
To calculate $\po(G)$, I followed a suggestion of the late Prof.~J.~F.~Adams:
\par \noindent
Firstly extend the action of $Q$ on $H^*({\rm B}N;\Bbb Z)=
\Bbb Z_{/p}[x_1,\ldots,x_n]$ to an action of $Q$ on $\Bbb F_{p^n}
[x_1,\ldots,x_n]$, calculate the fixed subring under this action, then
calculate the subring of this fixed by the Galois group of $\Bbb F_{p^n}$ over
$\Bbb Z_{/p}$. Using this method it may be shown that such groups $G$ do not
give counterexamples to $\po(G)\mid (|Q|,\po(N))$.
\par
Your theorem~2\cddot4, in the case $G_1$ (the metacyclic groups) is a special
case of a result of C.~T.~C.~Wall from `Resolutions for extensions of groups',
Proc. Camb. Phil. Soc. 57 (1961) 251--5. The case $G_2(1,1)$ is studied in
G.~Lewis `The integral cohomology rings of groups of order $p^3$' Trans. Amer.
Math. Soc. 132 (1968) 501--29. In this spectral sequence $d_2$ is not trivial;
for example $d_2:E_2^{1,4}\rightarrow E_2^{3,3}$ is not trivial. In my
preprint, which I have enclosed, I give sufficient information
(during the proof of theorem~3\cddot 6) to calculate the
order of $H^n({\rm B}G_2(a,1);\Bbb Z)$ for $a > 1$ (I use the name
$M(a+2)$ for your $G_2(a,1)$). I find that the order of
$H^n({\rm B}G_2(a,1);\Bbb Z)$ is $0$, $p^{a+1}$, $p^2$,
$p^{a+3}$ for $n=1$, 2, 3, 4 respectively. For $i+j \leq 4$, I
calculate that the addititve structure of $E_2$ for the spectral sequence you
consider is
$$\bigspec{1}{\relax\Bbb Z_{/p^a}\oplus\Bbb Z_{/p} \cr
\relax\Bbb Z_{/p}&\relax\Bbb Z_{/p} \cr
\relax\Bbb Z_{/p^a}&\relax\Bbb Z_{/p}&\relax\Bbb Z_{/p}\cr
-&-&-&-\cr
\relax\Bbb Z&-&\relax\Bbb Z_{/p}&-&\relax\Bbb Z_{/p}}$$
\par
I believe, therefore, that in the spectral sequence for $\langle A,C\rangle
\rightarrowtail G_2(a,1)\twoheadrightarrow \Bbb Z_{/p}$, where $a > 1$,
there is a non-trivial $d_n:E_n^{i,4-i}\rightarrow E_n^{i+n,5-i-n}$ for
some $i$ and $n$. I think that the values of $i$ and $n$ are 1 and 2
respectively.
\par
I have no plans to publish this work, but I shall include it in my Ph.~D.
thesis. \rm
\beginsubsection Remarks. Since receiving the above letter, Yagita has
retracted the non-metacyclic cases of his theorem~2\cddot4, and has sent the
author a proof of lemma~1\cddot7 under the strong condition that the action of
$G$ on $\coh N$ is nilpotent. The groups $G(a,1)$ considered by Yagita may
be presented as follows.
$$\langle A,B,C\mid A^{p^a}=B^p=C^p=[A,C]=[B,C]=1 \quad [A,B]=C\rangle$$
Yagita's theorem~2\cddot4 claimed that the spectral sequence with integer
coefficients for $G(a, 1)$ expressed as an extension of $\langle A,C
\rangle$ by $C_p$ collapses. We may consider $G(a, 1)$ as a central
extension of $\langle A^p\rangle $ by $P_2$. Since $G(a, 1)$ is nilpotent
of class two the corresponding circle group is isomorphic to $\sone \times
P_2$, which has cohomology given by the K\"unneth theorem.
$$\cohz * {\sone \times P_2} \cong \Bbb Z[\tau]\otimes \cohz * {P_2}$$
The class $\tau p^{a -1}+\beta$ in $\cohz 2 {\sone \times P_2}$ is the
first Chern class of a bundle over ${\rm B}(\sone\times P_2)$ with total space
${\rm B}G(a, 1)$. We could of course choose any other non-zero element
of $\cohz 2 {P_2}$ instead of $\beta$.
It is easy now to calculate the order of $\cohz n {G(a, 1)}$ using the
spectral sequence for this $\sone $ bundle.
\par
A linear action of $G$ on $\es^{2n-1}$ is an action induced by an $n$
dimensional complex representation of $G$. Linear actions on products of
spheres are products of linear actions. Linear actions induce the trivial
action on the cohomology of the spheres because $\gl_n(\Bbb C)$ is connected.
An element $g$ of $G$ will act fixed point freely in the linear action given by
a product of representations $\rho_1\times\ldots\times\rho_k$ if and only if
there is a $j$ such that the restriction of $\rho_j$ to $\langle g\rangle$ does
not contain the trivial representation of $\langle g\rangle$. We obtain the
following lemma.
\par
\proclaim \Lemma 6\cddot 1. If $G$ has a linear action on a product of $2n-1$
dimensional spheres such that every element of $G$ of order $p$ acts fixed
point freely, then $\pc(G)$ divides $2n$. \par
\proof
It suffices to check that for each $i:C_p\rightarrowtail G$ then $n(i)$ (as
defined in the above letter) divides
$n$. Let the action be specified by maps $\rho_1,\ldots,\rho_k$ from $G$ to
$\gl_n(\Bbb C)$, and let $C$ be the image of $C_p$ in $G$. If $C$ acts fixed
point freely, then there is a $\rho=\rho_j$ such that $\res^G_C(\rho)$ does not
contain the trivial representation. Then
$$ c.(\res^G_C(\rho))=\prod^n_{l=1}(1+\lambda_l u)$$
where $u$ is a generator for $\cohz 2 C$ and $\lambda_l$ is non-zero. But then
$c_n(\res^G_C(\rho))$ is a non-zero multiple of $u^n$.
\qed

\def\ess{\{1,\ldots,l\}}
\beginsection{7.} The Davis Complex.\par
The following chapter, which is largely expository although we hope to
obtain new results in this area eventually, concerns a construction due
to M.~Davis [Da], which assigns to a Coxeter group $G$ a simplicial complex
$D(G)$, which is contractible and has a simplicial $G$-action with finite
stabilisers.  We give a simplified account of a special case of Davis'
construction, and an elementary account of an example due to Bestvina [Be] of a
group whose cohomological dimension over $\Bbb F_2$ is greater than its
cohomological dimension over $\Bbb Q$.  None of these results are original,
although our proofs differ from those of [Da] and [Be].  Finally we suggest an
application to the construction of a 4-manifold with interesting homological
properties.
\par \mark{\ }
\proclaim \Definition 7\cddot1. A group $G$ is said to be a graph product if it
is generated by a collection $G_1,\ldots,G_l$ of finite subgroups, such that
for each $i\neq j$ the subgroup generated by $G_i$ and $G_j$ is either the
free product $G_i*G_j$ or the direct product $G_i\times G_j$.  In the special
case when each $G_i$ is a group of order two, $G$ is a right angled Coxeter
group.  A subgroup of $G$ generated by a subset of the $G_i$ shall be called a
special subgroup.  Note that the trivial subgroup is defined to be special
(generated by the empty set of $G_i$'s).
\par
\proclaim \Definition 7\cddot2. If $G$ is a graph product, then its Davis
complex $D(G)$ is the simplicial complex associated to the poset whose elements
are the left cosets of the finite special subgroups of $G$, ordered by
inclusion.
\par
Write $G(S)$ for the special subgroup of $G$ generated by the set of
$G_i$ for $i$ an element of $S$.  We may represent any $n$-simplex of $D(G)$ in
the form $(g,S_0,S_1,\ldots,S_n)$, where $g$ is an element of $G$ and
$S_0,\ldots,S_n$ is a strictly increasing chain of subsets of $\ess$ such that
$G(S_n)$ is finite.  Another such symbol $(h,T_0,\ldots,T_n)$ represents the
same simplex if and only if $T_i=S_i$ for all $i$ and $gG(S_0)=hG(S_0)$.  The
stabiliser of the simplex $(g,S_0,\ldots,S_n)$ is the subgroup $gG(S_0)g^{-1}$
of $G$, which is finite.
\par
We shall not prove that the Davis complex is contractible, but we shall prove
the following easier result.
\proclaim \prop 7\cddot3 (Davis).  The Davis complex of a graph product is
acyclic.  \par
\proof If for all pairs $i$, $j$ in $\ess$ the subgroup of $G$ generated by
$G_i$ and $G_j$ is $G_i\times G_j$, then $G$ is itself the direct product of
all of the $G_i$'s, and hence finite, and so $D(G)$ is a cone with vertex the
unique left coset of $G$ itself.  If not then without loss of generality we may
assume that $G_1$ and $G_l$ generate $G_1*G_l$.  Define subgroups $H$, $K$ and
$L$ of $G$ by
$$ H=G(\{1,\ldots,l-1\}),\quad K=G(\{2,\ldots,l\}),\quad
L=G(\{2,\ldots,l-1\}).$$
Then $H\cap K=L$, and $G=H*_LK$.  Consider now the subcomplex $X$ of $D(G)$
with simplices those $(g,S_0,\ldots,S_n)$ such that $l\notin S_n$.  By fixing a
transversal $T$ to $H$ and expressing $g=th$ for some $t\in T$ and $h\in H$, we
see that $X$ is isomorphic to a disjoint union of copies of $D(H)$, indexed by
$T$.  Note however that the action of $H$ on the $D(H)$ corresponding to $t$ is
twisted by conjugation by $t$.  Since inductively $D(H)$ may be assumed to be
acyclic we see that $X$ has no homology except in degree zero, and that
$H_0(X)$ is isomorphic as a $\Bbb ZG$-module to $\Bbb ZG/H$ (the permutation
module with basis the left cosets of $H$).  Similarly, if we define $Y$ and $Z$
by  \mark{The Davis Complex}
$$Y=\{(g,S_0,\ldots,S_n)|1\notin S_n\},\quad
Z=\{(g,S_0,\ldots,S_n)|1,l\notin S_n\},$$
then these subcomplexes have only degree zero homology, and as $\Bbb
ZG$-modules
$$H_0(Y)\cong \Bbb ZG/K,\quad H_0(Z)\cong\Bbb ZG/L.$$
Any simplex $(g,S_0,\ldots,S_n)$ of $D(G)$ must be contained in either $X$ or
$Y$ since $S_n$ cannot contain both 1 and $l$, and the intersection of $X$ and
$Y$ is $Z$.  The Mayer-Vietoris sequence for $D(G)$ implies that $H_i(D(G))=0$
for $i>1$, and gives
$$0\mapright{}H_1(D(G))\mapright{}\Bbb ZG/L\mapright{\eta}\Bbb ZG/H
\oplus \Bbb ZG/K\mapright{}H_0(D(G))\mapright{}0.$$
The map $\eta$ is often studied in the proof of the Mayer-Vietoris theorem for
group cohomology (see [Br]), and it is known that $\eta$ is injective and has
cokernel $\Bbb Z$. \qed\par
\proclaim \Corollary 7\cddot4. If $H$ is a torsion-free subgroup of a graph
product $G$ (for example the kernel of the map from $G$ to the direct product
of the $G_i$'s), then the chain complex for $D(G)$ is an $H$-free resolution
for $\Bbb Z$.  The chain complex for $D(G)$ with rational coefficients is a
$\Bbb QG$-projective resolution for $\Bbb Q$.  \par
\proof The chain complex for $D(G)$ is a resolution for the trivial module by
Proposition~7\cddot3.  The modules that occur in this resolution are
permutation modules with basis $G/K$ for various finite subgroups $K$.  If $H$
is torsion-free, then it cannot meet such a $K$, so acts freely.  The module
$RG/K$ is projective provided that the order of $K$ is invertible in $R$.
\qed \par
Before stating any further theorems it is useful to introduce another
simplicial complex $K(G)$ associated to the graph product $G$.  $K(G)$ has as
simplices the subsets $S$ of $\{1,\ldots,l\}$ such that $G(S)$ is finite.
Equivalently, $K(G)$ has vertices $\{1,\ldots,l\}$, edges the pairs $i$ and $j$
such that $G_i$ and $G_j$ commute, and is the full simplicial complex on its
1-skeleton.  (Recall that a simplicial complex is said to be full if whenever
its 1-skeleton contains a complete graph on $n$ vertices then the complex
contains an $(n-1)$-simplex whose boundary is that graph.)  The barycentric
subdivision of any simplicial complex is full.  Given a (finite) full
simplicial complex $K$, and a finite group $G_v$ for each vertex $v$ of $K$, we
may define a graph product $G$ with $K(G)=K$, by taking the free product of all
the $G_v$'s and adding the relations that $G_v$ and $G_w$ commute if $(v,w)$ is
an edge of $K$.
\proclaim \prop 7\cddot5.  The Davis complex $D(G)$ may be obtained by taking a
free $G$-orbit of copies of the cone on the barycentric subdivision of $K(G)$,
and identifying together parts of their boundaries as described below.
\par
\proof  By the boundary of a cone, we mean the simplices not containing the
cone point.  Simplices of $K(G)$ correspond to non-empty subsets $S$ of
$\{1,\ldots,l\}$ such that $G(S)$ is finite, so simplices of its barycentric
subdivision correspond to ascending chains $(S_0,\ldots,S_n)$ of non-empty
subsets such that $G(S_n)$ is finite.  We may add a cone point to this by
allowing our chains to contain the empty set.  Now we see that simplices of
$G\times CK(G)'$ are representable uniquely by symbols $(g,S_0,\ldots,S_n)$,
where $S_0,\ldots,S_n$ is a chain of subsets with $G(S_n)$ finite.  Putting the
equivalence relation on these symbols described after Definition~7\cddot2 is
equivalent to identifying parts of the boundaries of the cones.
\qed \par
{\it From now onwards, we shall only consider the case when $G$ is a
right-angled Coxeter group}\par
\proclaim \Theorem 7\cddot6 (Bestvina).  Let $X$ be the topological space
obtained by attaching a disc $D^2$ to a circle $\sone$ using the map
$z\mapsto z^n$ from $\sone$ to $\sone$.  Let $G$ be any right-angled Coxeter
group such that $K(G)\cong X$, and let $H$ be a torsion-free subgroup of $G$ of
finite index.  Then $H$ has cohomological dimension three, but $H^3(H;M)$ has
exponent dividing $n$ for all modules $M$.  Also $H^3(H;\Bbb ZH)$ is cyclic of
order $n$.  \par
\proof  By Corollary 7\cddot4, the Davis complex provides a finite free
resolution of $\Bbb Z$ over $\Bbb ZH$ of length three.  Using
Proposition~7\cddot3 we shall study this resolution.  The rest of this proof is
very similar to a \lq bare hands' proof that $H^2(X)$ is cyclic of order $n$.
Firstly, orient $C(K(G)')$ (the cone on the barycentric subdivision of $K(G)$)
by coherently orienting the cone on a disc having the same 2-simplices as
$K(G)'$, and then identifying the sides of the cone.  The third stage ($C_3$)
of the resolution consists of a direct sum of copies of $\Bbb ZG$, one for each
2-simplex of $K(G)'$.  The second stage ($C_2$) of the resolution consists of a
direct sum of three sorts of pieces.  Firstly, copies of $\Bbb ZG$
corresponding to interior 2-simplices of $C(K(G)')$ containing the cone point.
The map from $C_3$ to each of these components will be non-zero on exactly two
of $C_3$'s summands, and will be the identity on one and minus the identity on
the other.  Secondly, copies of $\Bbb ZG$ corresponding to 2-simplices of
$C(K(G)')$ containing the cone point and an edge of the disc.  The map from
$C_3$ to each of these components will be the identity map on $n$ of $C_3$'s
summands, and zero on the rest.  Finally, for each 2-simplex $\sigma$, a copy
of $\Bbb ZG/L$, where $L$ is the minimal vertex of $\sigma$ (necessarily a
subgroup of order two).  The map from $C_3$ to this component will be
projection on one of $C_3$'s summands, and zero on the rest.
\par
Now apply $\hom_{\Bbb ZH}(\_,M)$ to this resolution, and consider the
\lq adjoint map' from $C^2$ to $C^3$.  $C^3$ consists of a direct sum of copies
of $N=\hom_{\Bbb ZH}(\Bbb ZG,M)$, and $C^2$ consists of a direct sum of three
types of module, the first and second types being isomorphic to $N$.
The map from a summand of the first type to $C^3$ is of the form
$$m\mapsto (0,\ldots,0,m,0,\ldots,0,-m,0,\ldots,0),$$
and since the interior of the disc is connected the images of the first type
summands suffice to identify all the copies of $N$ that comprise $C^3$.  The
map from a second type summand to this quotient $N$ is multiplication
by $n$, and we see that $H^3(H;M)$ has exponent dividing $n$.
\par
Now specialise to the case when $M=\Bbb ZH$, so that $N\cong \Bbb ZG$.  Already
by quotienting $C^3$ by the image of first and second type summands of $C^2$ we
have shown that $H^3$ is a quotient of $RG$, where $R=\Bbb Z/(n)$.  A summand
of $C^3$ of the third type has the form $\Bbb ZG/\langle g_i\rangle$, where the
$g_i$ are the Coxeter generators for $G$, and the effect of quotienting by
these is that
$$H^3(H;\Bbb ZH)\cong RG/(\sum_i(1+g_i)RG).$$
Since $G$ is presented by the $g_i$ subject only to relations involving even
powers of each $g_i$ we see that the subgroup $G'$ of the group of units of
$RG$ generated by $-g_1,\ldots,-g_l$ is isomorphic to $G$ and $RG'=RG$.  Now
$H^3$ is just $RG'$ modulo its augmentation ideal, or $R$.
\qed\par
Next we shall look for conditions for $D(G)$ to be a combinatorial manifold.
\proclaim \Lemma 7\cddot7.  The link of a simplex of $D(G)$ is isomorphic to a
suspension of the link of a simplex of $K(G)'$, except that the 0-simplex
$(g,\emptyset)$ has link isomorphic to $K(G)'$.
\par
\proof  The link of the simplex $(g,S_0,S_1,\ldots,S_n)$ is the join of various
pieces, which we shall examine separately.  Firstly there is the complex
corresponding to the poset of special cosets strictly contained in $gG(S_0)$.
This is the barycentric subdivision of the Coxeter complex of $G(S_0)$, that is
the subdivision of the $m$ dimensional analogue of the octahedron (where $S_0$
contains $m+1$ elements), so is homeomorphic to $S^m$ except when
$S_0=\emptyset$,
when it is empty.  Secondly, there is the complex corresponding to the poset of
special cosets lying strictly between $gG(S_i)$ and $gG(S_{i+1})$.  If
$S_{i+1}\setminus S_i$ contains $m$ elements then this complex is the
barycentric subdivision  of the surface of an $(m-1)$-simplex (or empty if
$m=1$)---this is because faces of an $(m-1)$-simplex correspond to non-empty
subsets of its vertices.  Thirdly, there is the complex corresponding to the
poset of special cosets strictly containing $gG(S_n)$.  This is isomorphic to
the link in $K(G)'$ of any $|S_n|-1$-simplex with $S_n$ as its minimal vertex,
except when $S_n=\emptyset$, in which case it is the whole of $K(G)'$.
\qed \par
\proclaim \Corollary 7\cddot8.
$D(G)$ is a combinatorial $(n+1)$-manifold if and only if $K(G)$ is
homeomorphic to $S^n$.
\par
\proof  Keep track of the dimensions of the spheres that arise in the proof of
the preceeding lemma. \qed \par
I.~M.~Chiswell has used the Davis complex to provide a formula for the Euler
characteristic of a graph product [Ch], which we state here only in the case
when $G$ is a right-angled Coxeter group.
\proclaim \Theorem 7\cddot9 (Chiswell).  If $G$ is a right-angled Coxeter
group, then
$$\chi(G)=1-{1\over 2}\sum_{i\geq 0}{n_i\over2^i},$$
where $n_i$ is the number of $i$-simplices in $K(G)$.
\qed\par
P.~H.~Kropholler has suggested combining Corollary 7\cddot8 and
Theorem~7\cddot9 to try to exhibit a closed aspherical 4-manifold of negative
Euler characteristic as follows.  Given any full triangulation of the three
sphere, form the corresponding right-angled Coxeter group $G$, and let $H$ be a
torsion-free subgroup of finite index.  $H$ acts freely and cocompactly
on $D(G)$, so $D(G)/H$ is a closed manifold.  The strong version of Davis'
theorem (namely that $D(G)$ is contractible) implies that $D(G)/H$ is
aspherical, and its Euler characteristic is $|G:H|$ times that of $G$.
Unfortunately this seems not to work.  For example:
\proclaim \prop 7\cddot10.  If $K$ is any triangulation of $S^3$ having $m_i$
$i$-simplices, and $G$ is the Coxeter group corresponding to the barycentric
subdivision of $K$, then $\chi(G)$ is positive.
\par
\proof  Find $n_i$ (the number of $i$-simplices of $K'$ in terms of the $m_j$,
and substitute in to the formula of Theorem~7\cddot9.  It is helpful to note
that the $n_i$ (and the $m_i$) satisfy the following relations:
$$n_0-n_1+n_2-n_3=0\qquad 2n_3=n_2.$$
(The first of these comes from $\chi(S^3)=0$, and the second from the fact that
each 2-simplex in a 3-manifold bounds exactly two 3-simplices.
\qed \par
In the case of triangulations of $S^3$ arising from convex polytopes, the
theory of toric varieties may be brought to bear on this question.
G.~K.~Sankaran has pointed out that given a triangulation of $S^3$, the
signature of the cup product as a quadratic form on $H^4$ of the corresponding
toric variety is a positive multiple of the expression occuring in
Theorem~7\cddot9.
\par

\def\book#1/#2/#3/#4/#5/{\item{#1} #2, {\it #3,} #4, {\oldstyle #5}.
\par\smallskip}
\def\paper#1/#2/#3/#4/#5/(#6) #7--#8/{\item{#1} #2, #3, {\it #4,} {\bf #5}
({\oldstyle#6}) {\oldstyle #7}--{\oldstyle#8}.\par\smallskip}
\def\prepaper#1/#2/#3/#4/#5/#6/{\item{#1} #2, #3, {\it #4} {\bf #5} {#6}.
\par\smallskip}
\frenchspacing
 
\beginsection References.\par
\mark{\ }
\paper Al1/K.~AlZubaidy/Metacyclic $p$-groups and Chern classes/Illinois J.
Math./26/(1982) 423--431/
\paper Al2/K.~AlZubaidy/Rank 2 $p$-groups, $p>3$, and Chern classes/Pacific J.
Math./103/(1982) 259--267/
\paper Ar/S.~Araki/Steenrod reduced powers in the spectral sequences associated
with a fibering/Mem. Fac. Sci., Kyusyu/11/(1957) 15--97/
\paper At/M.~F.~Atiyah/Characters and the cohomology of finite groups/Publ.
Math. IHES/9/(1961) 23--64/
\prepaper Be/M.~Bestvina/The virtual cohomological dimension of Coxeter
groups///preprint/
 
\paper Bl/N.~Blackburn/Generalisations of certain elementary theorems on
$p$-groups/Proc. L.M.S./11/(1961) 1--22/

\book Br/K.~S.~Brown/Cohomology of Groups/Springer Verlag/1982/
\book Bu/W.~Burnside/Theory of Finite Groups/C.U.P./1897/
\book CE/H.~Cartan and S.~Eilenberg/Homological Algebra/P.U.P./1956/
\prepaper Ch/I.~M.~Chiswell/The Euler Characteristic of a graph
product///preprint/
\book Co/J.~H.~Conway {\it et al.}/Atlas of Finite Simple Groups/O.U.P./1985/
\paper Da/M.~Davis/Groups generated by reflections and aspherical manifolds not
covered by Euclidean space/Ann. of Math./117/(1983) 293--324/
\paper Di/L.~E.~Dickson/A fundamental system of invariants of the general
modular linear group with a solution of the form problem/Trans.
Amer. Math. Soc./12/(1911) 75--98/
\paper Die/T.~Diethelm/The mod $p$ cohomology rings of the nonabelian split
metacyclic $p$-groups/Arch. Math./44/(1985) 29--38/
\paper Ev1/L.~Evens/A generalisation of the transfer map in the cohomology of
groups/Trans. Amer. Math. Soc./108/(1963) 54--65/
\paper Ev2/L.~Evens/On the Chern classes of representations of finite
groups/Trans. Amer. Math. Soc./115/(1963) 180--93/
\book Go/D.~Gorenstein/Finite Groups/Chelsea/1980/
\paper Ha/P.~Hall/A contribution to the theory of groups of prime power
order/Proc. London Math. Soc./36/(1933) 29--95/
\paper He/D.~Held/The simple groups related to $M_{24}$/J. of Alg./13/(1969)
253--96/
 
\paper Hi/G.~Hirsch/Quelques propri\'et\'es des produits de Steenrod/
C.~R.~Acad. Sci. Paris/241/(1955) 923--25/
\prepaper Hu1/J.~Huebschmann/Chern classes for metacyclic groups///to appear/
\prepaper Hu2/J.~Huebschmann/Perturbation theory and free resolutions for
nilpotent groups of class 2/J. of Algebra,//to appear/
\prepaper Hu3/J.~Huebschmann/Cohomology of nilpotent groups of class 2/J. of
Algebra,//to appear/
\mark{References}
\prepaper La/D.~S.~Larson/The integral cohomology rings of split metacyclic
groups/unpublished report,//Univ. of Minnesota ({\oldstyle1987})/
\prepaper Lea/I.~J.~Leary/The integral and mod-$p$ cohomology rings of some
$p$-groups/unpublished essay,//({\oldstyle1989})/
\paper Lew/G.~Lewis/Integral cohomology rings of groups of order $p^3$/
Trans. Am. Math. Soc./132/(1968) 501--29/
\paper Ma/J.~P.~May/Matric Massey products/J. of Alg./12/(1969) 533--68/
\book McC/J.~McCleary/User's Guide to Spectral Sequences/Publish or
Perish/1985/
\prepaper Mo/B.~Moselle/Calculations in the cohomology of finite
groups/unpublished essay,//({\oldstyle1989})/
\paper Re/D.~L.~Rector/Modular characters and $K$-theory with coefficients in a
finite field/J. Pure and Appl. Alg./4/(1974) 137--58/
\paper St/N.~E.~Steenrod/Products of cocycles and extensions of mappings/Ann.
Math./48/(1947) 290--320/
\paper Sw/R.~G.~Swan/The $p$-period of a finite group/Illinois J. Math./4/(1960)
341--6/
\paper TY/M.~Tezuka and N.~Yagita/Cohomology of finite groups and the
Brown-Peterson cohomology/Springer LNM/1370/(1989) 396--408/
\paper TY2/M.~Tezuka and N.~Yagita/Cohomology of finite groups and
Brown-Peterson cohomology II/Springer LNM/1418/(1990) 57--69/
\paper Th1/C.~B.~Thomas/Chern classes and metacyclic $p$-groups/
Mathematika/18/(1971) 196--200/
\paper Th2/C.~B.~Thomas/Riemann-Roch formulae for group
representations/Mathematika/18/(1971) 196--200/
\book Th3/C.~B.~Thomas/Characteristic Classes and the Cohomology of Finite
Groups/C.U.P./1986/
\paper Th4/C.~B.~Thomas/Characteristic classes and 2-modular representations of
some sporadic simple groups/Contemp. Math./96/(1989) 303--18/
\paper Va/R.~Vasquez/Nota sobre los cuadrados de Steenrod en la sucesion
espectral de un espacio fibrado/Bol. Soc. MAt. Mexicana/2/(1957) 1--8/
\paper Wa/C.~T.~C.~Wall/Resolutions for extensions of groups/Proc. Cam. Phil.
Soc./57/(1961) 251--5/
\prepaper We/P.~J.~Webb/A split exact sequence for Mackey functors///to appear/
\prepaper Wei/E.~A.~Weiss/unpublished thesis/Bonner Mathematische Schriften,//
({\oldstyle1969})/
\paper Ya/N.~Yagita/On the dimension of spheres whose.../Quart. J. Math.
Oxford/36/(1985) 117--27/
\prepaper Ya2/N.~Yagita/Letter to the author///dated {\oldstyle 15}th January
{\oldstyle1990}/
 
\end